\documentclass[a4paper, 12pt]{amsart}
\usepackage[english]{babel}
\usepackage{amsfonts,amsmath,amssymb,amsthm,graphicx,afterpage,url}
\input{xy}
\xyoption{all}
\usepackage{epsfig}
\usepackage{overpic}
\usepackage{pgf,tikz}
\usetikzlibrary{arrows}

\usepackage{textcase}

\newcommand{\cycles}[1]{}
\newcommand{\algor}[1]{}
\newcommand{\invadraw}[1]{}

\newcommand{\compile}[1]{}

\def\lk{\mbox{lk}}

\def\sgn{\mathop{\fam0 sgn}}

\def\lk{\mathop{\fam0 lk}}

\def\R{{\mathbb R}} \def\Z{{\mathbb Z}}  

\DeclareMathOperator{\wu}{wu}

\theoremstyle{plain}
\newtheorem{theorem}{Theorem}[section]
    \newtheorem{lemma}[theorem]{Lemma}

    \newtheorem{conjecture}[theorem]{Conjecture}

\theoremstyle{definition}

\newtheorem{problem}[theorem]{Problem}
\newtheorem{remark}[theorem]{Remark}
\newtheorem{assertion}[theorem]{Assertion}

\newtheoremstyle{mydefinition}
  {3pt}
  {3pt}
  {\normalfont}
  {\parindent}
  {\bfseries}
  {.}
  { }
  {}

\theoremstyle{mydefinition}
\newtheorem{pr}[theorem]{Problem}

\graphicspath{
{./old_figs/}
{./graphs_in_space_figs/}
}
\usepackage{float}
\makeatletter

\@addtoreset{figure}{section}
\makeatother

\usepackage{hyperref}
\hypersetup{
   colorlinks   = true, 
   urlcolor     = blue, 
   linkcolor    = blue, 
   citecolor   = red 
} 

\begin{document}

\newpage
\title{Linking invariants of spatial graphs}

\author{E. Alkin, Yu. Khromin, A. Skopenkov} 

\thanks{\emph{E. Alkin, A. Skopenkov:} Moscow Institute of Physics and Technology.
\newline
\emph{Yu. Khromin:} Physical and Mathematical School №146 (Perm).
\newline
\emph{A. Skopenkov:} Independent University of Moscow, \url{https://users.mccme.ru/skopenko}.
\newline
This is a research project presented at the 2026 Summer Conference of the Tournament of Towns, see \url{https://www.turgor.ru/en/lktg/index.php} .
We are grateful to E. Flapan, A. Miroshnikov, A. Mizev, O. Nikitenko and R. Nikkuni for useful discussions, and to O. Nikitenko for preparation of some figures. 
Some figures are prepared using the paper \cite{Ta95}.} 

\date{}


\begin{abstract}
We recall definitions of linking numbers and Wu--Simon numbers for spatial graphs. 
We expose a `converse' to the Conway--Gordon--Sachs theorem (i.e. description of linking functions for embeddings $K_6\to\R^3$), 
and some results on Wu--Simon numbers. 
We conjecture and discuss a generalization of the Conway--Gordon--Sachs theorem to multiple linking. 
The exposition is based on plane diagrams, so no knowledge of spatial geometry is required. 
\end{abstract}

\maketitle
\tableofcontents

Main definitions are given in \S\ref{s:lkpd}.  
Main results are stated in \S\S \ref{s:main},\ref{s:maini} using the notion of a \emph{linking number} defined (via projections) in \S\ref{s:lkn}. 
Proofs of the main results are sketched in \S\S \ref{s:real},\ref{s:sol},\ref{s:wu}. 

\smallskip
\textbf{Remark} (on the style of this text).
\emph{In this text we expose a theory as a sequence of problems,} see e.g. \cite{HC19}, \cite[Introduction, Learning by doing problems]{Sk21m} and the references therein. 
Most problems are useful theoretical facts. 
So this text could in principle be read even without solving problems.
If a mathematical statement is formulated as a problem, then the objective is to prove this statement.
Open-ended questions are called {\bf riddles}; here one must come up with a clear wording, and a proof. 
\emph{If a problem is named `theorem' (`lemma', `corollary', etc.), then this statement is considered to be more important.}
Usually we \emph{formulate} beautiful or important statements \emph{before} giving a sequence of results (lemmas, propositions, etc.) which constitute its \emph{proof}.
We give hints on that after the statements but we do not want to deprive you of the pleasure of finding the right moment when you finally are ready to prove the statement.
In general, if you are stuck on a certain problem, try looking at the next ones; they may turn out to be helpful.
Problems marked by star and remarks are not used in the sequel; although problems with star are not necessarily complicated, they can be omitted during the first round of problem solving.
Important definitions are highlighted in \textbf{bold} for easy navigation.


\section{Links and plane diagrams}\label{s:lkpd}

\begin{figure}[H]\centering

\begin{minipage}{0.24\textwidth}
\centering
\includegraphics[scale=1.2]{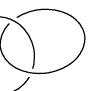}
\par (a)
\end{minipage}
\hfill
\begin{minipage}{0.24\textwidth}
\centering
\includegraphics[scale=1.2]{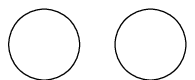}
\par (b)
\end{minipage}
\hfill
\begin{minipage}{0.24\textwidth}
\centering
\includegraphics[scale=0.6]{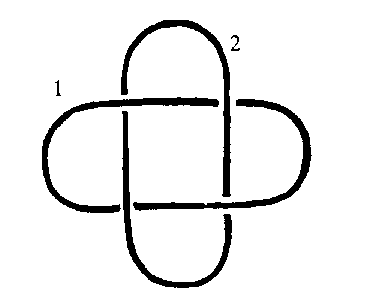}
\par (c)
\end{minipage}
\hfill
\begin{minipage}{0.24\textwidth}
\centering
\includegraphics[scale=0.4]{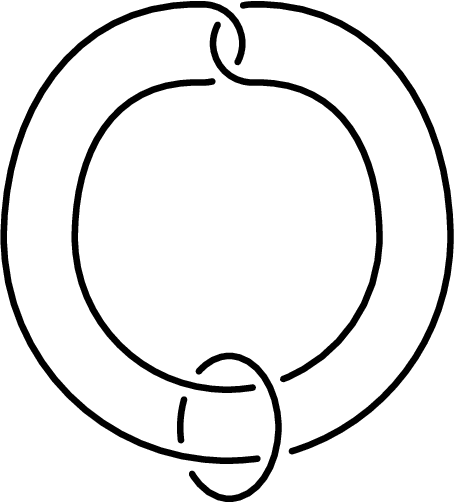}
\par (d)
\end{minipage}

\vspace{0.2cm}
\caption{The Hopf link (a), the trivial link (b) and another two links (c, d)}
\label{f:hopf}
\end{figure}

You can imagine a {\it knot} as a thin elastic string whose ends have been glued together; see the components in Figure~\ref{f:hopf}. 
Rigorously speaking, a {\bf knot} is a spatial closed non-self-intersecting polygonal line.
A {\bf link} (Figure~\ref{f:hopf}) is an ordered set of pairwise disjoint knots, which are called the {\it components} of the link.

Links are usually represented by their `nice' plane projections called {\it plane diagrams} (Figure~\ref{f:hopf}; 
by projection we understand shadow).
Imagine laying down  the rope on a table and carefully recording how it crosses itself (i.e. which part lies on top of the other).
A {\it plane diagram of a link} is its \emph{generic} projection onto a plane,  
together with an indication 
at every crossing, which segment `passes under' and which `passes over' (Figure~\ref{f:hopf}). 
The projections of the same link on different planes can look quite dissimilar.

In this text we work rigorously not with links and their plane diagrams, but with \emph{abstract plane diagrams}. 
Here `abstract' means `not necessarily constructed from any link', and is mostly omitted.  

Rigorously, a set of closed polygonal lines in the plane is said to be {\it generic} if 

(i) their vertices are pairwise distinct,
 
(ii) any two segments of the polygonal lines either are disjoint, or intersect at their common endpoint, or intersect at their common interior point, and
 
(iii) no three interiors of these segments have a common point.


Common interior points as in (ii) are called {\bf crossings}. 

An {\bf (abstract) plane diagram} consists of 

$\bullet$ a generic ordered set of some closed polygonal lines in the plane, 
 
$\bullet$ at every crossing, an indication which segment `passes under' and which `passes over' (Figure~\ref{f:hopf}). 

Each of the closed polygonal lines is called a {\bf component} of the plane diagram. 

\begin{assertion}\label{p:diag} * For any one-component abstract plane diagram there is a knot projected to this diagram.
(Such a knot need not be unique; see though \cite[Assertion 1.1.c]{Sk20u}.)
\end{assertion}

\section{Linking number modulo 2}\label{s:lkn2}

The notion of a \emph{linking number} of two closed curves in 3-space was invented by K. Gauss in his studies of electromagnetism \cite{RN11}. 
Speaking informally, the linking number is the number of turns one curve makes around the other. 
It is not clear how to define this number (even informally) for complicated pairs of curves. 
Let us provide a rigorous definition in the language of  
plane diagrams. 
 
The {\bf linking number modulo 2}, $\lk_2\in\Z_2$, of a two-component plane diagram is the number modulo 2 of the crossings at which the first component passes over the second component.

\begin{figure}[H]\centering
\begin{minipage}{0.2\textwidth}
\centering
\includegraphics[scale=1.2]{3-7.borr-link.eps}
\par (b)
\end{minipage}
\quad
\begin{minipage}{0.2\textwidth}
\centering
\includegraphics[scale=1.2]{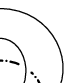}
\par (w)
\end{minipage}
\quad
\begin{minipage}{0.2\textwidth}
\centering
\includegraphics[scale=1.2]{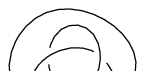}
\par (t)
\end{minipage}
\vspace{0.3cm}
\caption{The Borromean rings (b), the Whitehead link (w), and the trefoil knot (t). 
\newline
Every two of the Borromean rings are unlinked, but all the three of them are linked 
(i.e. cannot be deformed in three pairwised disjoint balls, cf. \S\ref{s:isot})}
\label{f:borwhitre}
\end{figure}

\begin{pr}\label{a:lk2} (a,b,c) For your choice of the order of the components, find the linking numbers modulo 2 of the plane diagrams in Figure~\ref{f:hopf}, of the two-component subdiagrams 
of Borromean rings in Figure~\ref{f:borwhitre}.b, and of every pair of disjoint cycles of length 3 in Figure~\ref{f:ex-graph-emb}, right.   
\end{pr}

\begin{figure}[H]\centering
\includegraphics[scale=0.6]{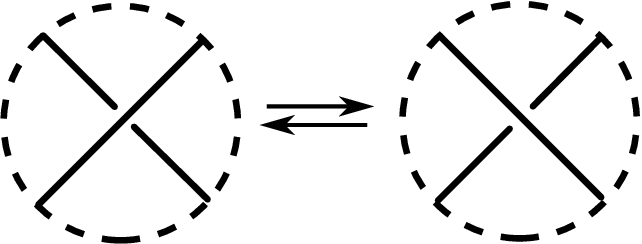}
\vspace{0.2cm}
\caption{Crossing change.
\emph{We show parts bounded by dashed circles of two plane diagrams. 
The diagrams are identical outside the parts.
No other segments of the plane diagrams except for the pictured ones intersect the parts.
Same for Figures~\ref{f:cross-ori},~\ref{f:triplic},~\ref{f:reid},~\ref{f:reid-1-5},~\ref{f:reid-6-7}.}}
\label{f:cross}
\end{figure}

A {\bf crossing change} is a change from overcrossing to undercrossing or vise versa, see Figure~\ref{f:cross}.

\begin{assertion}\label{con-lk2} For any two-component plane diagrams $D_+$ and $D_-$ which differ by the crossing change 
of a crossing of different components, we have $\lk\phantom{}_2D_+-\lk\phantom{}_2D_- = 1\in\Z_2$. 
\end{assertion}

\begin{assertion}\label{a:lk2-pr}
Switching the components 
preserves the linking number modulo 2.
\end{assertion}

In the proof (and below) you can  
use without proof the following 
{\it Parity lemma:} 
any two generic closed polygonal lines in the plane 
intersect at an even number of points.
For a discussion and a proof see \S1.3 `Intersection number for polygonal lines in the plane' of \cite{Sk18, Sk}. 

The \textit{inverse} of a polygonal line $l = L_1 L_2 \ldots L_m$ is the polygonal line $\overline l := L_m L_{m-1} \ldots L_1$.
The \textit{concatenation} of polygonal lines $l = L_1\ldots L_m$ and $p = L_mP_2\ldots P_n$ is the polygonal line $lp := L_1 \ldots L_m P_2 \ldots P_n$.  

\begin{lemma}\label{l:add2} Let $l_1,l_2,l_3$ be polygonal lines joining point $A$ to point $B$, and disjoint with a polygonal line $p$.
Then
$$\lk\phantom{}_2(l_1\overline{l_2}, p) + \lk\phantom{}_2(l_2\overline{l_3}, p) = \lk\phantom{}_2(l_1\overline{l_3}, p).$$
\end{lemma}

\section{Linking number}\label{s:lkn}

\begin{figure}[H]\centering
\includegraphics[scale=0.5]{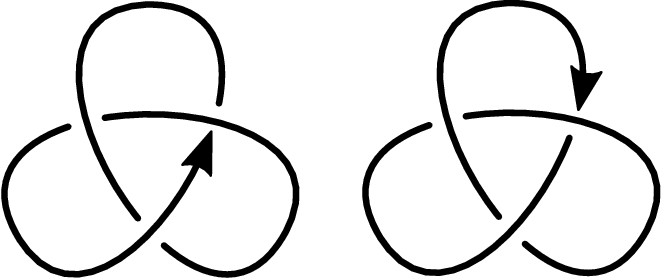}
\caption{Two trefoil knots with the opposite orientations}
\label{f:knot-ori}
\end{figure}

One knows what an oriented polygonal line is, so one knows what an \emph{oriented link} is (Figure~\ref{f:knot-ori}), and what an \emph{oriented 
plane diagram} is.

\begin{figure}[H]\centering
\includegraphics[scale=0.8]{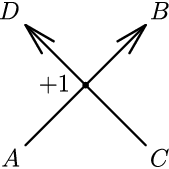}
\qquad 
\includegraphics[scale=0.8]{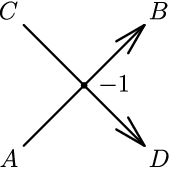}
\caption{The sign of intersecting pair of segments}
\label{f:sign}
\end{figure}

Let $(\overrightarrow{AB},\overrightarrow{CD})$ be an ordered pair of intersecting vectors (oriented segments) in the plane. 
Define {\bf the sign} of the pair to be $+1$ if $ABC$ is oriented clockwise and to be $-1$ otherwise (Figure~\ref{f:sign}).

The {\bf linking number} $\lk$ of a two-component oriented plane diagram 
is the sum of signs at all those crossings at which the first component passes above the second component.
Here at every crossing the {\it first} (the {\it second}) vector is the oriented edge of the first (the second) component.

\begin{pr}\label{e:lk} (a,b,c) For your choice of the order of the components and their orientations, find the linking numbers of the plane diagrams in Figure~\ref{f:hopf}, of the two-component subdiagrams of Borromean rings in Figure~\ref{f:borwhitre}.b, and of every pair of cycles of length 3 in Figure~\ref{f:ex-graph-emb}, right.  

(d) For any integer $n$ there is an oriented two-component plane diagram whose linking number is $n$. 

(e) For any triple of integers $k,l,m$ there is an oriented three-component plane diagram whose pairwise linking numbers are $k,l,m$. 
\end{pr}

\begin{figure}[H]\centering
\includegraphics[scale=0.35]{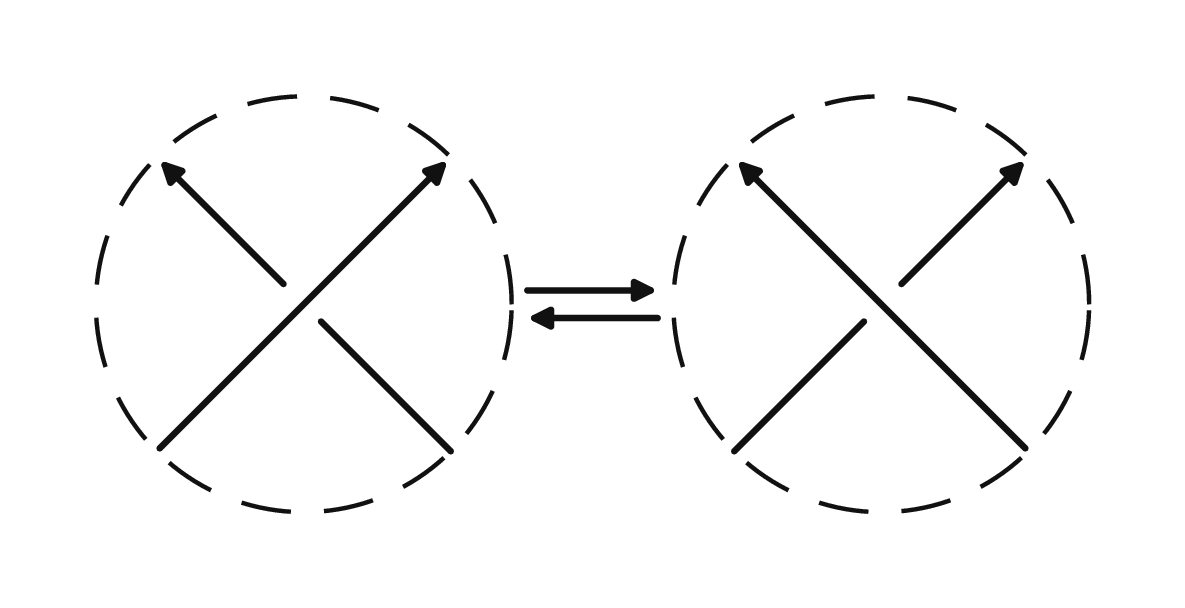}
\caption{Oriented crossing change}
\label{f:cross-ori}
\end{figure}

\begin{assertion}\label{con-lk} For any oriented two-component plane diagrams $D_+$ and $D_-$ (left and right in Figure~\ref{f:cross-ori}, respectively) 
which differ by the oriented crossing change of a crossing of different components, we have $\lk D_+-\lk D_- = 1$. 
\end{assertion}


\begin{assertion}\label{p:lk-pr} 
(a) Reversing the orientation of either of the components 
negates the linking number.

(b) Switching the components preserves the linking number.
\end{assertion}


For (b) the following lemma is useful. 

\begin{lemma}\label{l:half} The linking number equals $(n_{12}-n_{21})/2$, where
 
$\bullet$ $n_{12}$ is the sum of signs at all those crossings at which the first component passes above the second component, and

$\bullet$ $n_{21}$ is the sum of signs at all those crossings at which the second component passes above the first component.
\end{lemma}

In the proof (and below) you can  
use without proof the following {\it Triviality lemma:} for any two generic closed oriented polygonal lines in the plane,  
the sum of the signs of their intersection points is zero.
For a discussion and a proof see \S1.3 `Intersection number for polygonal lines in the plane' of \cite{Sk18}, \cite{Sk}.

\begin{lemma}\label{l:add} Let $l_1,l_2,l_3$ be polygonal lines joining point $A$ to point $B$, and disjoint with a closed polygonal line $p$.
Then
$$\lk(l_1\overline{l_2}, p) + \lk(l_2\overline{l_3}, p) = \lk(l_1\overline{l_3}, p).$$
\end{lemma}

\section{Main results on sets of linking numbers}\label{s:main}

\begin{figure}[H]
\centering
\includegraphics[scale=1.5]{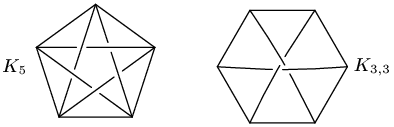}\qquad
\includegraphics[scale=0.8]{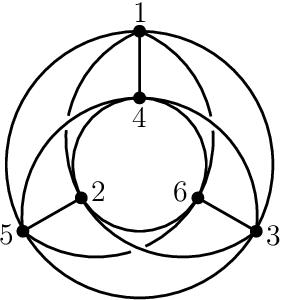}
\caption{Embeddings of graphs 
in 3-space $\R^3$}
\label{f:ex-graph-emb}
\end{figure}

An \emph{embedding} $K\to\R^3$ of a graph $K$ in 3-space $\R^3$ is
a drawing without self-intersections (Figure~\ref{f:ex-graph-emb}).
We consider those embeddings of graphs
for which edges are represented by polygonal lines.
A simple rigorous definition of such embeddings is given in Remark~\ref{rem:def-emb}.
But you do not need it to solve the problems, because all results can be 
restated and proved in terms of 
\emph{plane diagrams}, see Remark~\ref{r:emvsdi}. 


Denote by

$\bullet$ $[n]$ the set $\{1, 2, \ldots, n\}$;

$\bullet$ $K_n$ the complete graph with the vertex set $[n]$;   

$\bullet$ $K_{m,n}$ the complete bipartite graph with parts $[m]$ and $[n]'$ (we denote by $A'$ a copy of $A$).  

\begin{theorem}[Conway--Gordon--Sachs]\label{t:cgs}
For any embedding $K_6\to\R^3$
the sum of the linking numbers over all unordered pairs of disjoint cycles of length 3 in $K_6$ is odd.
\end{theorem}

This result is not part of this project, but serves only as motivation for Theorems~\ref{t:div3}, \ref{t:odd}, 
Assertion~\ref{t:div3m} and Conjecture~\ref{t:div3mp}. 
 
\emph{Of this and the following sections, only Theorems~\ref{t:div3},~\ref{t:odd},~\ref{p:k5nontr}, 
 Assertion~\ref{t:div3m} and Problem~\ref{a:crossing} are intended for solving within this project.
Some of them might be not so easy to solve without hints given in later sections.}

\begin{remark}[restatement in terms of plane diagrams]\label{r:emvsdi}
The \emph{plane diagram of an embedding of a graph in $\R^3$}
is defined analogously to the plane diagram of a link (using a notion of \emph{generic graph drawing in the plane} defined in \S\ref{s:isot}).
An \emph{(abstract) plane diagram of a graph} is defined in \S\ref{s:isot}, analogously to an (abstract) plane diagram (defined in \S\ref{s:lkpd}). 
Theorem~\ref{t:cgs} is restated as follows:  
\emph{for any plane diagram of $K_6$, the sum of the linking numbers over all unordered pairs of disjoint triangles in $K_6$ is odd.} 
Analogously one restates all those results of this text that are stated in terms of embeddings of graphs in 3-space.   
\end{remark}
 
\begin{theorem}\label{t:div3}  
For any graph $K$ there is an embedding $K\to\R^3$ such that the linking number of any two disjoint cycles of length 3 is divisible by 3. 
\end{theorem} 

This is proved by R. Nikkuni in 2026. 
For motivation see the end of this section. 

For graphs with less than 6 vertices Theorem~\ref{t:div3}, Assertion~\ref{t:div3m} and Conjecture~\ref{t:div3mp} are not substantive. 
In Theorem~\ref{t:div3} (and presumably in Conjecture~\ref{c:borr}) mod 3 can be replaced by mod $m$ for any odd integer $m$. 
But not by mod 2 (because of Theorem~\ref{t:cgs}).

\begin{assertion}\label{t:div3m}  
The same as Theorem~\ref{t:div3} with `is divisible by 3' replaced by `is non-zero'. 
\end{assertion}


\begin{conjecture}\label{t:div3mp}  
For any graph $K$ there is an embedding $K\to\R^3$ such that for any three pairwise disjoint cycles of length 3, 
the greatest common divisor of their pairwise linking numbers is 1. 
\end{conjecture}

The analogue without `of length 3' is wrong by Negami's theorem \cite{Ne91, PS05} 
(and by \cite[Theorem 3.1]{FD05} which gives a specific constant for the corresponding particular case of \cite{Ne91}).  

\begin{problem}\label{p:k6} For any integer $k$ there is an embedding $K_6\to\R^3$ such that the linking number of
 
$\bullet$ one (unordered) pair of disjoint $3$-cycles is $2k+1$, and
 
$\bullet$ any other pair of disjoint $3$-cycles is zero.
\end{problem}

This elementary result appeared in recent research \cite[Proposition 1.2]{KS20} \cite[end of \S2.2]{AT21}, although it could have been known earlier.
 
\begin{theorem}\label{t:odd}
Consider any $10$ integers $n_{123,456},n_{124,356},\ldots$, corresponding to the $10$ unordered partitions of the set $[6]$ into two $3$-element subsets.
If the sum of these integers is odd, then there is an embedding $K_6\to\R^3$ such that
the linking number of every pair $\{ijk,pqr\}$ of disjoint $3$-cycles is $n_{ijk,pqr}$.
Here we assume that $i<j<k$ and $p<q<r$, and take the orientations $(ijk)$ and $(pqr)$ on the cycles. 
\end{theorem}

This is proved in \cite{Mi}; this follows from Lemmas~\ref{l:move} and \ref{l:expr} below, cf. Theorem~\ref{t:odd1}. 

\begin{problem}[open]\label{op:set_lk} 
For every $n$ describe all ordered sets 
realizable as sets of linking numbers 
of an embedding $K_n\to\R^3$.
\end{problem}
 
\emph{General hints.} Describe complicated constructions by simple modifications of simple constructions! 
Most of the presented geometric problems are reduced to purely algebraic problems. 

\begin{figure}[H]\centering
\includegraphics[scale=.1]{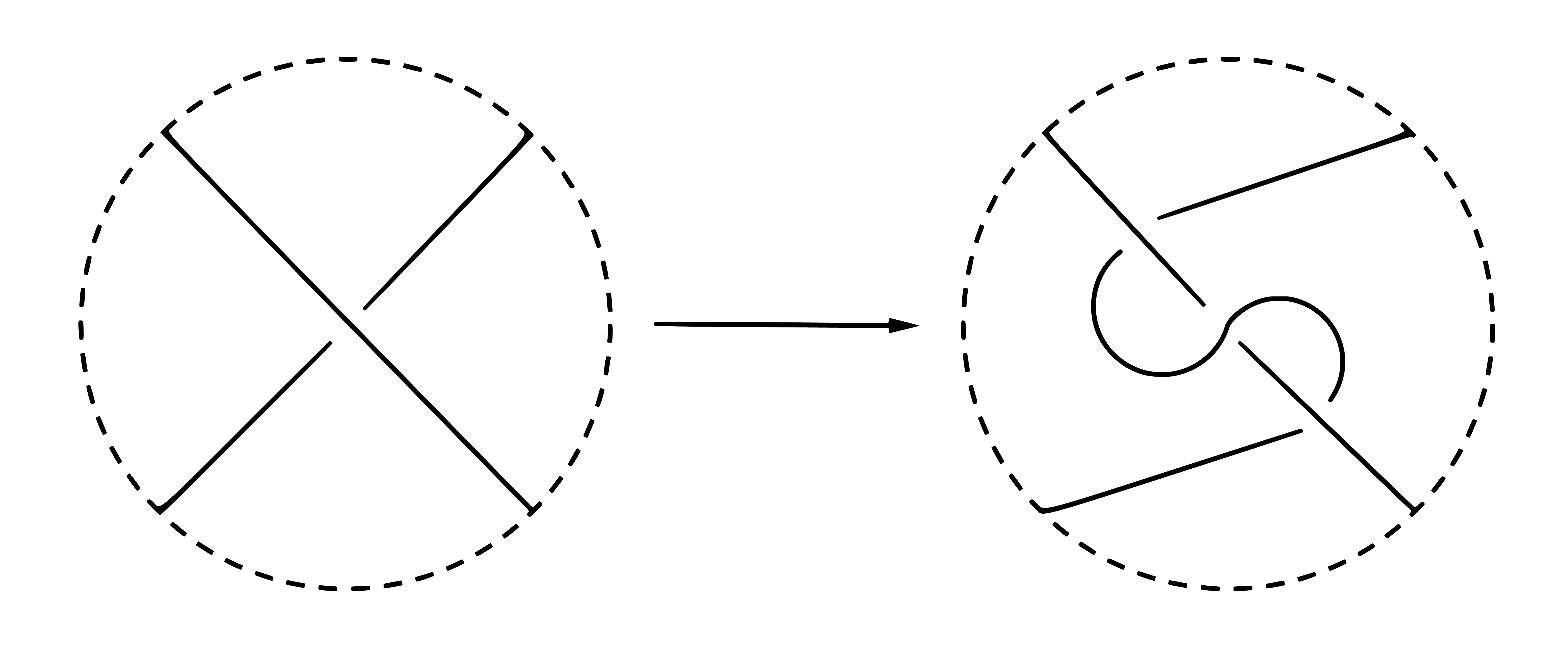}
\caption{Triplication of a crossing} 
\label{f:triplic}
\end{figure}

\smallskip
\emph{Hint to Theorem~\ref{t:div3}:} use Figure \ref{f:triplic} and Lemma \ref{l:half}.

\textbf{A generalization of the Conway--Gordon--Sachs Theorem \ref{t:cgs} to multiple linking.}

This subsection is not a part of this project but serves only as  motivation for Theorem~\ref{t:div3}. 
 
\begin{conjecture}\label{c:borr} 
Let $K_{11}\to\R^3$ be an embedding such that the linking number of any two disjoint cycles of length 3 is divisible by 3.  
Then there is a 
`mod 3 Borromean' triple of pairwise disjoint cycles of length 3 in $K_{11}$, i.e. a triple whose mod 3 triple linking number is non-zero. 
\end{conjecture}

\textbf{Remark.} 
For Borromean rings, see Figure \ref{f:borwhitre}.b and e.g. \cite{Zi10, Sk24}, \cite[\S4]{Sk} together with the references therein. 
The integer \emph{(Milnor) triple linking number} is defined for an oriented 3-component link whose pairwise linking numbers are zeros, see e.g.
\cite[\S2, Theorem~(3)]{MM01}, \cite{Sk24}, \cite[\S4.9]{Sk}. 
Analogously one defines the mod $m$ \emph{triple linking number} for an oriented 3-component link whose pairwise linking numbers are divisible by $m$.
Thus for Conjecture~\ref{c:borr} to be substantive, one needs Theorem~\ref{t:div3}. 
By Assertion~\ref{t:div3m} and Conjecture~\ref{t:div3mp} one cannot omit the `mod 3 pairwise triviality' assumption in Conjecture~\ref{c:borr}, 
even if in the conclusion one replaces `mod 3' by `mod $m$ for some $m>1$'. 

For large enough $n$ and any embedding $K_n\to\R^3$ there are three pairwise disjoint cycles (not necessarily of length 3) in $K_n$ whose images form (a link isotopic to) Borromean rings. 
This holds by Negami's theorem \cite{Ne91, PS05}. 
However, the analogue of Conjecture~\ref{c:borr} for $K_n$ with large enough $n$ instead of $K_{11}$, is not trivial.  
 
For different but related problems on intrinsic $n$-linking see e.g. survey \cite{Na20} and the references therein. 
 
A prospective proof and a counterexample for $K_{10}$ instead of $K_{11}$ is given in \cite{AS}. 


\begin{conjecture}\label{c:borrb} Let  $K_{6r}\to\R^3$  be an embedding such that for every $s=2,\ldots,r-1$ 
the $s$-fold linking number of any $s$ pairwise disjoint cycles of length 3 is divisible by $r!+1$.  
Then there is a `Brunnian mod $(r!+1)$' $r$-tuple of  pairwise disjoint cycles of length 3 in $K_{6r}$, 
i.e. an $r$-tuple whose $r$-fold linking number mod $(r!+1)$  is non-zero. 
\end{conjecture}

 
\textbf{Remark.} 
(a) The conjecture is a \emph{linking version of the topological Tverberg conjecture}, see e.g. survey \cite{Sk16}. 
So the conjecture could be non-trivial. 
We would not be surprised if the conjecture is true for $r$ being a prime power, and false otherwise. 
For a (dis)proof, configuration spaces might be useful. 


(b) For Conjecture~\ref{c:borrb} to make sense, one needs the following analogue of Theorem~\ref{t:div3}. 

Conjecture. \emph{For every $r\ge3$ there is an embedding $K_{6r}\to\R^3$ such that for every $s=2,\ldots,r-1$ the $s$-fold linking number of any $s$ pairwise disjoint cycles of length 3 is divisible by $r!+1$.}  

The above proof of Theorem~\ref{t:div3} can be interpreted as making connected sum with the Hopf link, for any pair of edges contributing to the linking number (as defined in \S\ref{s:lkn}). 
One can try to prove the above conjecture for $r=3$ by making connected sum with the Borromean link,for any triple of edges contributing to the triple linking number. 
For this, one needs calculation of the triple linking number via plane diagrams. 

\section{Main results on isotopy and homology of spatial graphs}\label{s:maini}

\begin{figure}[H]\centering
\includegraphics[scale=0.8]{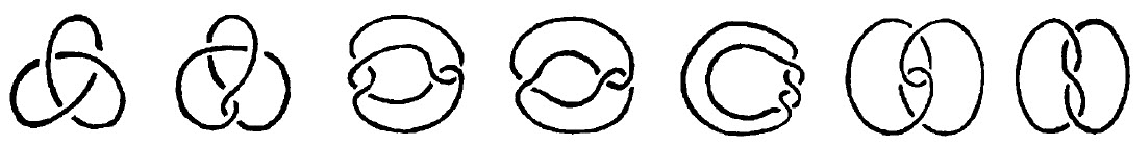}
\caption{Knots isotopic to the trefoil knot} 
\label{f:trefoil}
\end{figure}

By an {\it isotopy} of a link we mean its continuous deformation in 3-space as a set of thin elastic strings; throughout the deformation, no self-intersections (of one string with itself or of different strings) are allowed.
See Figure~\ref{f:trefoil} and the reduction to plane diagrams  in \S\ref{s:isot}. 

\begin{pr}\label{a:crossing} Please prove these assertions at an intuitive level (because you do not have a rigorous definition of an isotopy). 

(a) All the knots (represented by their plane diagrams) in 
Figure~\ref{f:trefoil} are isotopic to each other.

(b) The Whitehead link (Figure~\ref{f:borwhitre}.w) is isotopic to the link from Figure~\ref{f:hopf}.d.

(c) Two plane diagrams of a link are obtained by projections to different planes.
Then these diagrams are  obtained by projections to one plane from some two (different but) isotopic links. 
\end{pr}

\emph{Hint to (a).} `Probably the best way ... 
is to make a model of the trefoil knot by using a rope or shoelace, 
and then move it around from one position to the other' \cite[\S2]{Pr95}. 
 
\begin{theorem}\label{p:k5nontr}  (a) There are non-isotopic embeddings $K_5\to\R^3$ whose restrictions to any simple cycle in $K_5$ are isotopic. 

(b) (chirality \cite{Si86}) No embedding $K_5\to\R^3$ is isotopic to its mirror image. 

(a',b') The same for $K_{3,3}$ instead of $K_5$. 
\end{theorem}

The embeddings from (a,a') are distinguished by the integer \emph{Wu--Simon number}, and \emph{bipartite Wu--Simon number},  defined in \S\ref{s:wu}. 
The same invariants are used to prove (b,b'). 

\begin{remark}[embeddings modulo knots]\label{r:emk}
\emph{Isotopy} classification (readily calculable, see \cite[Remark 1.2]{Sk16c}) of knots is so hard that it is unlikely to be obtained.
If embeddings of a graph $K$ to 3-space are isotopic, then their restrictions to any simple cycle in $K$ are isotopic.
So isotopy classification of embeddings of graphs in 3-space is unlikely to be obtained 
(although \emph{connected sums} with non-isotopic knots can potentially yield isotopic graphs).
Therefore it is interesting to define and classify `embeddings modulo knots'.
There are several formalizations of this idea: \emph{homology} (defined in \S\ref{s:homol}), \emph{link homotopy}, \emph{almost isotopy}, \emph{non-ambient isotopy}; see a survey \cite{Sk'}.
\end{remark}


\begin{theorem}\label{t:tani} Two embeddings of a graph (not necessarily connected) in 3-space are homologous if and only if 
their restrictions to any subgraph homeomorphic to $K_5,K_{3,3}$ or $K_3\sqcup K_3$ are homologous. 
Or, equivalently, if all the linking numbers, the Wu--Simon numbers, and the  bipartite Wu--Simon numbers (of the restrictions) are equal.
\end{theorem}

This theorem is not part of this project. 
See Remark \ref{r:conf}. 
 
\begin{problem}[open; cf. Problem \ref{op:set_lk}]\label{op:set_lkwu} 
(a) For every $n$ find sets of integers realizable as sets of linking numbers and Wu--Simon numbers of an embeddings $K_n\to\R^3$ (start with $n=6$). 

(b) Analogous question adding bipartite Wu--Simon numbers. 
\end{problem}

\section{Realization of linking numbers: to the proof of Theorem~\ref{t:odd}}\label{s:real}

Let us reformulate Theorem~\ref{t:odd}. 
Denote by $X$ the set of unordered partitions of $[6]$ 
into two 3-element subsets. 
A \textbf{linking function} $l_f:X\to\Z$ of an embedding $f:K_6\to\R^3$ assigns to any partition $x=\{\{i,j,k\},\{p,q,r\}\}\in X$,  where $i<j<k$ and $p<q<r$, the linking number of $f(ijk)$ and $f(pqr)$. 
(The set of all functions $X \to \Z$ can be regarded as the set of all vectors with integer coordinates, indexed by partitions from~$X$.)
 
\begin{theorem}[cf. Theorem~\ref{t:odd}]\label{t:odd1} 
For any function $\varphi: X \to \Z$ with an odd sum of the values there is an embedding $f$ with $l_f=\varphi$. 
\end{theorem}

\begin{figure}[H]\centering
\includegraphics[scale=0.6]{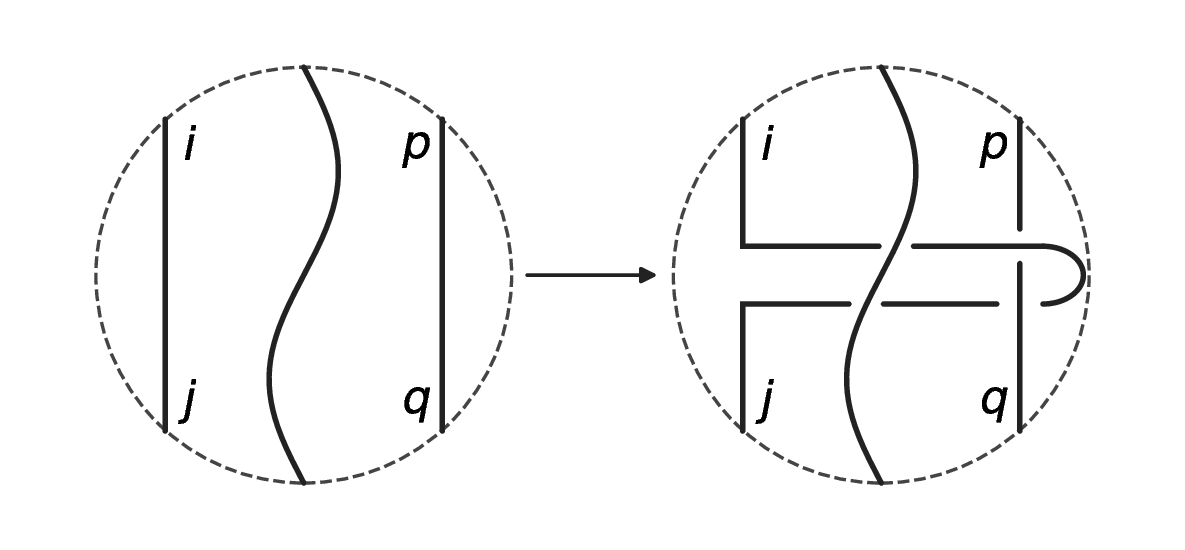}\qquad
\caption{Finger move of one edge across the other}
\label{f:finger_gra3spa}
\end{figure}

\begin{pr}\label{p:move} (a) How does 
$l_f$ change under the transformation of an embedding $f:K_6\to\R^3$ shown in Figure \ref{f:finger_gra3spa}? 

(b) (riddle) Same question for an embedding $f:K_n\to\R^3$. 
\end{pr}

Having solved Problem \ref{p:move}.b, you can prove Assertion~\ref{t:div3m}.

An \emph{interesting} pair is an ordered pair $(ij,pq)$ of disjoint oriented edges of $K_6$ such that $i<j$ and $p<q$. 
Take an interesting pair $(ij,pq)$. 
Take $a,b\in [6]$ such that $\{i,j,p,q,a,b\}=[6]$. 
Define a function $v_{ij,pq} : X\to\Z$ 
by\footnote{Although the formulas use a choice of  $a$ and $b$ from $\{a,b\}$, they are independent of the choice.}
$$v_{ij,pq}(\{i,j,a\},\{p,q,b\}) := \sgn(a-i)\sgn(a-j)\sgn(b-p)\sgn(b-q) \in \{ +1, -1 \},$$
$$v_{ij,pq}(\{i,j,b\},\{p,q,a\}) := \sgn(b-i)\sgn(b-j)\sgn(a-p)\sgn(a-q) \in \{ +1, -1 \},$$
and $v_{ij,pq}(x) := 0$ for every other $x\in X$.

E.g. $v_{15,23}(\{2,3,5\},\{1,4,6\})=0$ and 
$$v_{15,23}(\{2,3,6\},\{1,4,5\}) = v_{15,23}(\{1,5,4\},\{2,3,6\}) = 1\cdot(-1)\cdot1\cdot1 = -1.$$

\begin{pr}\label{p:intepair}  How does the function $v_{ij,pq}$ depends on 
switching the edges $ij,pq$? 
\end{pr}


\begin{lemma}\label{l:move} 
For any embedding $f:K_6\to\R^3$ and interesting pair $(ij,pq)$  there are embeddings $f_+,f_-:K_6\to\R^3$ such that $l_{f_\pm}-l_{f}=\pm v_{ij,pq}$. 
\end{lemma}

A function $v: X \to \Z$ is said to be \emph{realizable} if $v$ is equal to some integer linear combination of functions $v_{ij,pq}$, 
over interesting pairs $(ij,pq)$. 

\begin{lemma}[\cite{Mi}]\label{l:expr} Any function $X \to \Z$ with an even sum of values is realizable.
\end{lemma}

The lemma follows from Assertions~\ref{a:coord}.ab. 

For any $x\in X$ define the \emph{indicator function} 
$$\chi_x:X\to\Z\quad\text{by}\quad \chi_x(y)=\begin{cases}1 &\text{if }x=y \\ 0&\text{otherwise}\end{cases}.$$

\begin{assertion}\label{a:coord} 
(a) For any two partitions from $X$, either the sum or the difference of their indicator  functions is realizable.

(b) Twice the indicator function of any partition from $X$ is realizable.  
\end{assertion}


\textbf{Remark.} Here is a different more conceptual definition (not used in this project). 
An \emph{interesting} pair is an ordered pair of disjoint oriented edges of $K_6$. 
Take an interesting pair $(ij,pq)$. 
Take $a,b\in [6]$ such that $\{i,j,p,q,a,b\}=[6]$. 
Define a function $v_{ij,pq} : X\to\Z$ by 
$$v_{ij,pq}(\{i,j,a\},\{p,q,b\}) := \sgn(i,j,a)\sgn(p,q,b),\quad v_{ij,pq}(\{i,j,b\},\{p,q,a\}) := \sgn(i,j,b)\sgn(p,q,a),$$
and $v_{ij,pq}(x) := 0$ for every other $x\in X$.
Here the \emph{sign} of a vector of 3 integers is 
$$\sgn(k,l,m) := \sgn(k-l)\sgn(l-m)\sgn(m-k)\in\{+1,-1\}.$$ 
This is the sign of a permutation arranging $k,l,m$ in increasing order. 

E.g. $v_{15,23}(\{2,3,5\},\{1,4,6\})=0$ and 
$$v_{15,23}(\{2,3,6\},\{1,4,5\}) = v_{15,23}(\{1,4,5\},\{2,3,6\}) = \sgn(1,5,4)\sgn(2,3,6) = (-1)\cdot1=-1.$$

\section{Some answers, hints and solutions to \S\ref{s:lkpd}-\S\ref{s:real}}\label{s:sol}

\begin{figure}[h]\centering
\includegraphics[scale=0.5]{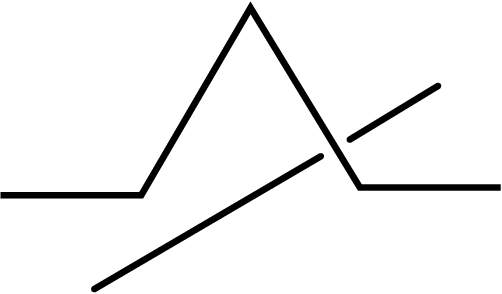}
\caption{The bridge over a crossing}
\label{f:bridge}
\end{figure}

{\bf \ref{p:diag}.} See Figure \ref{f:bridge}.
On the upper edge of every crossing of the plane diagram choose two points, close to the intersection, 
and lying on the opposite sides of the crossing.
Replace the line segment between the two chosen points by a `bridge' rising above the plane diagram, which joins these two points. 
Replace all crossings by the corresponding bridges. 
We obtain a knot whose projection is the union of segments of given plane diagram.

\smallskip
{\bf \ref{a:lk2}.} \textit{Answers:} 
\qquad (a) 1 for the Hopf link and 0 for other links; \qquad (b) 0;  

(c) 1 for the pair of cycles $(123)$ and $(456)$; 0 for other pairs.
 
\smallskip
{\bf \ref{a:lk2-pr}.} Take a plane diagram of a link. 
By the Parity lemma (stated after Assertion~\ref{a:lk2-pr})
the number of crossings where the first component passes above the second one has the same
parity as the number of crossings where the second component passes above the first one.
This completes the proof.

\smallskip
{\bf \ref{e:lk}.} \textit{Answers:}
 \qquad (a) 1; 0; 2; 0;  \qquad(b) 0. 

(c) $+1$ for the pair of cycles $(123)$ and $(456)$; 0 for other pairs.
 
\smallskip
{\bf \ref{p:lk-pr}.} (a) Reversing the orientation of either of the components negates the sign of every crossing.

(b) Switching the components negates the sign of every crossing. 
Thus (b) follows by Lemma~\ref{l:half}.

\smallskip
{\bf \ref{l:half}.} This follows because by the Triviality lemma (stated after the Lemma) $n_{12}+n_{21}=0$.  
 
\smallskip
{\bf \ref{t:div3}.} We follow a text by R. Nikkuni, with his kind permission. 

Replace each crossing by 
three crossings as shown in Figure~\ref{f:triplic}.
Take any two disjoint oriented 3-cycles 
in the graph. 
Then the sign of each crossing 
of these cycles is defined both in the original diagram and in the new one. 
For each crossing of the original diagram, the edge that passed \emph{over}, after the replacement, 
passes \emph{over} in two of the 
new crossings, and \emph{under} in the 
remaining new crossing.
The signs of 
the two new crossings coincide with the sign of the original crossing. 
The sign of the 
remaining new crossing is the opposite to the sign of the original crossing.
For the original diagram, we use the notation $n_{12}$ and $n_{21}$ from Lemma~\ref{l:half}, 
while for the new diagram we define $n'_{12}$ and $n'_{21}$ analogously.
For each original crossing that contributes $\pm1$ to the sum $n_{12}-n_{21}$, 
the corresponding 
three new crossings contribute $\pm(1-(-1)+1)=\pm3$ to the sum $n'_{12}-n'_{21}$.
Then, for any two disjoint cycles, the sum $n'_{12}-n'_{21}$ is divisible by $3$. 
Hence, by Lemma~\ref{l:half}, the obtained planar diagram is as required.

\smallskip
{\bf \ref{t:div3m}.} \emph{Hint.} Choose a pair of cycles with zero linking number. 
Modify the planar diagram so that the linking number for this pair of cycles becomes nonzero, and no new pairs with zero linking numbers appear. 
For this, apply Lemma~\ref{l:move} repeatedly. 
In the process, the linking numbers of the required pair and of some other pairs of cycles will change by $\pm N$ for an arbitrary preassigned integer $N$.

\smallskip
{\bf \ref{p:k6}.} This follows from Lemma~\ref{l:move}, Assertion~\ref{a:coord}.a, and the result of Problem~\ref{e:lk}.c.

\begin{figure}[h]\centering
\includegraphics[scale=0.25]{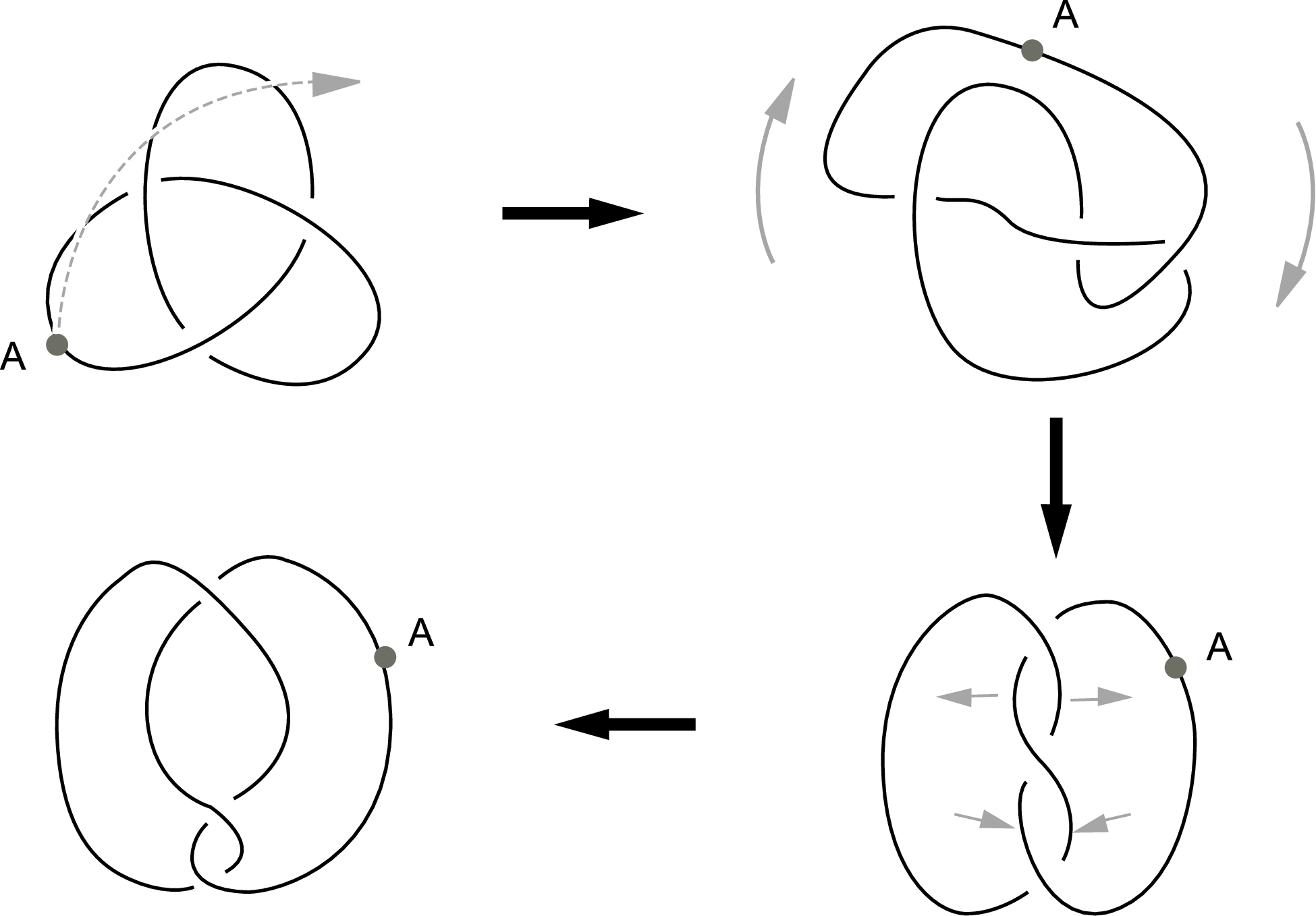}
\caption{Isotopy the trefoil knot}
\label{f:itref}
\end{figure}

\smallskip
{\bf \ref{a:crossing}.} (a) Figure \ref{f:itref} (prepared by D. Kroo) gives some hints concerning transformations of the trefoil knot. 

\smallskip
{\bf \ref{p:move}.} (a) See Lemma~\ref{l:move}. 
  
\smallskip
{\bf \ref{p:intepair}.} The function does not change.

\smallskip
{\bf \ref{l:move}.} We present a simplified exposition of the proof from \cite{KS20}.

We give a solution for the embedding $f_-$. 
Informally, we obtain $f_-$ by turning the edge $f(ij)$ by a $-1$ twist around the edge $f(pq)$, as in Figure~\ref{f:finger_gra3spa}.
We now give a rigorous construction.

On the edge $ij$, choose nearby points $I$ and $J$ with no crossings between them, 
so that the orientation of the segment $IJ$ coincides with that of the edge $ij$.
On the edge $pq$, similarly choose points $P$ and $Q$ so that the orientation of the segment $PQ$ coincides with that of the edge $pq$.
Using a finger move, replace the segment $IJ$ of the original diagram (i.e. the diagram of the embedding $f$) by a non-self-intersecting oriented polygonal line $l_1$ in general position, so that $l_1$ passes through the points $P$ and $Q$, as shown in Figure~\ref{f:finger_gra3spa}. 
In particular, at the crossing $P$ the polygonal line $l_1$ passes over, while at the crossing $Q$ and at all other crossings it passes under.

In the original diagram take any two disjoint cycles $x$ and $y$ of length $3$.
Denote by $x_-$ and $y_-$ the corresponding cycles in the obtained diagram (i.e. the diagram of the embedding $f_-$).
Denote by $\overline x$ the partition of the vertex set of the graph $K_6$ into two 3-element subsets, one of which consists of the vertices of the cycle $x$.
In what follows, after introducing orientations on the cycles $x$ and $y$, we shall take the consistent orientations on $x_-$ and $y_-$.

If the edge $ij$ is contained in neither of the cycles $x$ and $y$, then for any orientation on these cycles we have
$$\lk(x_-,y_-) = \lk(x,y).$$
Hence $l_{f_-}(\overline x)-l_f(\overline x)=0.$

If the edge $ij$ belongs to one of the cycles $x$ and $y$, then without loss of generality it belongs to $x$. 
Orient the cycle $x$ in the same direction as the edge $ij$. 
Regard $x$ as the first cycle.
Denote by $l_2$ the oriented polygonal line which is the complement of the oriented segment $IJ$ in the cycle $x$.

Now, if the cycle $y$ does not contain the edge $pq$, then for any choice of orientation on it, by Lemma~\ref{l:add} we obtain
$$\lk(x_-,y_-)=\lk(l_1l_2,y)=\lk(l_1JI,y)+\lk(IJl_2,y)=0+\lk(x,y).$$
Hence $l_{f_-}(\overline x)-l_f(\overline x)=0.$

If the cycle $y$ does contain the edge $pq$, then denote the vertices of the cycles $x$ and $y$ so that $x=ija$, $y=pqb$. 
Orient the cycle $y$ in the same direction as the edge $pq$.
Then the sign of the crossing at $P$ is $-1$.
So by Lemma~\ref{l:add} we obtain
$$\lk(x_-,y_-)=\lk(l_1l_2,y)=\lk(l_1JI,y)+\lk(IJl_2,y)=-1+\lk(x,y).$$
Hence $l_{f_-}(\overline x)-l_f(\overline x)=-\sgn(a-i)\sgn(a-j)\sgn(b-p)\sgn(b-q)$.

From all the above, the formula $l_{f_-}-l_f=-v_{ij,pq}$ follows.

The construction of the embedding $f_+$ is obtained from the above construction by changing the crossings at the points $P$ and $Q$ 
(overcrossing is changed to undercrossing and vice versa).
  
\smallskip
{\bf \ref{a:coord}.} We follow \cite{Mi}. 
 
(a) For any different partitions $x,y \in X$ there are $\alpha \in x,\beta\in y$ such that $|\alpha\cap\beta|=2$. 
Then there is an interesting pair $(ij,pq)$ such that $\{i,j\}=\alpha\cap\beta$ and $\{p,q\}=\overline\alpha\cap\overline\beta$. 
Hence $v_{ij,pq}(z)= \begin{cases}\pm1&\text{if }z=x\text{ or }z=y\\0&\text{otherwise}\end{cases}$.
So $\pm v_{ij,pq}=\chi_x\pm \chi_y$.

(b) Take any partition $x\in X$, and denote $s := \{\{1,2,3\},\{4,5,6\}\} \in X$. 
By (a), $\chi_x\pm\chi_s$ is realizable, so $2\chi_x\pm2\chi_s$ is realizable. 
Then $2\chi_s = v_{12,45}+v_{16,34}+v_{23,56}$, so $2\chi_s$ is realizable. 
Hence $(2\chi_x\pm 2\chi_s)\mp2\chi_s = 2\chi_x$ is also realizable.

(Alternatively, for any permutation $\sigma$ of $[6]$ we have $\pm2\chi_{\sigma(s)} = v_{\sigma(12,45)}+v_{\sigma(16,34)}+v_{\sigma(23,56)}$, where we extend the definition of $v_{ij,kl}$ by the Remark after Assertion \ref{a:coord}.) 


\section{Isotopy via plane diagrams}\label{s:isot}

\begin{figure}[H]\centering
\includegraphics[scale=0.45]{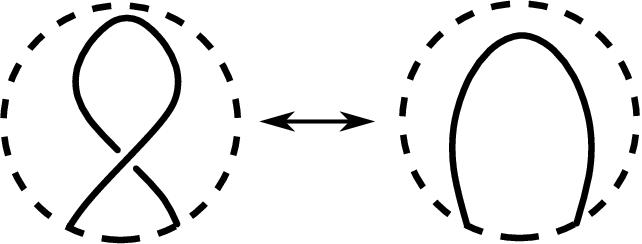}\qquad
\includegraphics[scale=0.45]{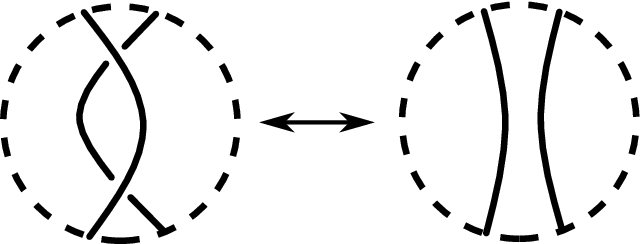}\qquad
\includegraphics[scale=0.45]{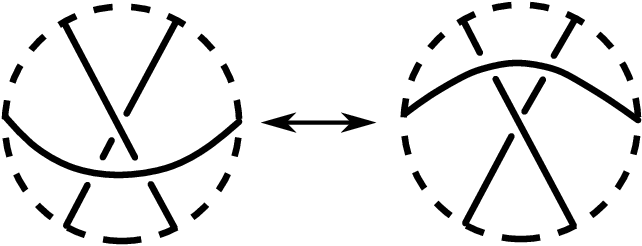}
\caption{Reidemeister moves}
\label{f:reid}
\end{figure}

Instead of links up to isotopy, in this text we study 
plane diagrams up to equivalence generated by {\bf Reidemeister moves} (Figure~\ref{f:reid}) and {\it plane isotopy}.
By a {\it plane isotopy} of a 
plane diagram we mean its continuous deformation in the plane, without introducing self-intersections. 
This is motivated by the following theorem (which is not part of this project). 

\begin{theorem}[Reidemeister; {\cite[\S1.7]{PS96}}]\label{l:reid} Two links are isotopic if and only if some plane diagram of the first link can be obtained from some plane diagram of the second one by Reidemeister moves and plane isotopy.
\end{theorem}

\begin{pr}\label{a:crossing-r} For the first pair of knots in Figure~\ref{f:trefoil} decompose your isotopy between them 
into Reidemeister moves shown in Figure~\ref{f:reid} and plane isotopy.
\end{pr}

\begin{lemma}\label{l:lk2} (a) The linking number modulo 2 \qquad 

(b) The linking number
 
 is preserved under Reidemeister moves.
\end{lemma}

This lemma is easily proved separately for each Reidemeister move.
For (b), check that the signs of all crossings do not change.

See more e.g. in \cite[\S\S 1-3]{Sk20u}. 

For an analogue of Theorem~\ref{l:reid} for graphs we need the following definitions.


Let $K$ be a graph with $V$ vertices. 
A (piecewise-linear, PL) {\bf graph drawing} $f$ of $K$ in $\R^d$, $d = 2, 3$ is 
 
$\bullet$ a set of $V$ points in $\R^d$ corresponding to the vertices, together with
 
$\bullet$ a set of (non-closed) polygonal lines in $\R^d$ corresponding to the edges; each polygonal line joins a pair from the set of points, which pair corresponds to an edge.

The point corresponding to a vertex $v$ is called the {\bf image} $f(v)$ of a vertex $v$ under $f$.
The polygonal line corresponding to the edge $\sigma$ (more precisely, the union of its segments), is called the {\bf image} $f(\sigma)$ of the edge $\sigma$ under the mapping $f$.
The {\bf image} of a set of edges is the union of the images of all edges in the set.

\begin{remark}[rigorous definition of an embedding]\label{rem:def-emb}
An \emph{embedding}\footnote{Standard definitions of an `embedding' uses the notion of a map not of the notion of a graph drawing.}
of $K$ in $\R^d$, $d = 2, 3$ is a graph drawing of $K$ in $\R^d$ such that

$\bullet$ the images of vertices are pairwise distinct,

$\bullet$ the images of edges are non-self-intersecting polygonal lines, and

$\bullet$ none of the images of edges intersects the interior of any other image.\footnote{Then, the images of any two edges either do not intersect, or intersect only at their common endpoint.}

One can prove that for any graph $K$ there is an embedding of $K$ in $\R^3$.
\end{remark}

A graph drawing in the plane is said to be \emph{generic} if it satisfies the conditions (ii), (iii) from~\S\ref{s:lkpd}, and the condition 

(i') the images of vertices are pairwise distinct, and any internal vertex of the image any edge is different from any other vertex (from the same or another edge).


Common interior points as in the condition (ii) are called {\bf crossings}.

An {\bf (abstract) plane diagram} of a graph $K$ consists of 

$\bullet$ a generic graph drawing of $K$ in the plane, 
 
$\bullet$ at every crossing, an indication which segment `passes under' and which `passes over' (Figure~\ref{f:ex-graph-emb}).\footnote{Rigorously speaking, the indication is a map $\Phi:\{(x,y)\in K\times K\ :\ x\ne y,\ fx=fy\}\to\{+1,-1\}$ such that $\Phi(x,y)=-\Phi(y,x)$, where $f:K\to\R^2$ is a map corresponding to the graph drawing.} 

\begin{figure}[H]\centering
\includegraphics[scale=1.]{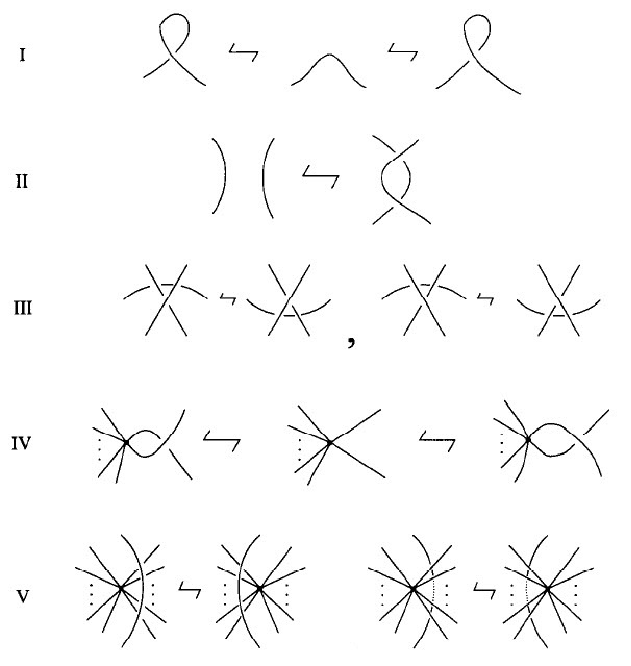}
\caption{Reidemeister moves for plane diagrams of graphs}
\label{f:reid-1-5}
\end{figure}

Analogously to Theorem~\ref{l:reid},
\emph{two embeddings of a graph in $\R^3$ are isotopic if and only if some plane diagram of the first embedding (see Remark~\ref{r:emvsdi}) can be obtained from some plane diagram of the second one by Reidemeister moves for plane diagrams of graphs (Figure~\ref{f:reid-1-5}) and plane isotopy.}

\section{Wu--Simon number}\label{s:wu}

Before studying this section please read \cite[\S\S 1.3, 1.4]{Sk} and proof of \cite[Lemma 1.4.3 and Remark 1.4.4b]{Sk}. 
This would help you invent \emph{the Wu--Simon number} defined below, or understand its definition 
(although the definition is formally independent of \cite[\S\S 1.3, 1.4]{Sk}). 
See the definition of a plane diagram of a graph at the end of \S\ref{s:isot}.

For a plane diagram $D$ of a graph $K$, we usually identify an edge $e$ of $K$ with its image under the generic graph drawing that is a part of $D$.


Let $D$ be a plane diagram of $K_5$.
Consider the orientations on segments of~$D$ induced by oriented cycles $12345$ and $14253$ in $K_5$.
At every crossing, consider the intersection sign, where the first vector passes above the second vector.
Color a crossing of non-adjacent edges~$e_1, e_2$ of~$K_5$ in \emph{blue}, if both~$e_1$, $e_2$ belong to the cycle $14253$, and in \emph{red}, otherwise.
The {\bf Wu--Simon number} $\wu D$ is defined \cite[\S4]{Ta94} to be the difference of
the sum of signs at all the red crossings, and 
the sum of signs at all the blue crossings.

 
By the van Kampen--Flores theorem for the plane \cite[Lemma 1.4.3]{Sk18},
any plane diagram of~$K_{5}$ has an odd number of crossings of non-adjacent edges, so the Wu--Simon number is odd. 

\begin{problem}\label{e:wu-k5}
(a) Find the Wu--Simon number of the plane diagram from Figure~\ref{f:ex-graph-emb} (left), for your choice of numeration of vertices.

(b) For any odd 
$n$ there is a plane diagram of $K_5$ whose Wu--Simon number is $n$. 
\end{problem}

\begin{assertion}\label{con-wu-k5} 
For any plane diagrams $D_+$ and $D_-$ of $K_5$ (left and right in Figure~\ref{f:cross-ori}, respectively) 
which differ by the oriented crossing change   of a crossing of non-adjacent edges, we have $\wu D_+-\wu D_- = \pm 2$.
\end{assertion}

\begin{assertion}\label{a:wu-k5-auto}* 
Let $D,D_{\sigma}$ be plane diagrams of $K_5$ such that $D_{\sigma}$ is obtained from $D$ by a renumeration $\sigma: [5] \to [5]$ of the vertices of $K_5$.
Then $\wu D_{\sigma} = \sgn\sigma \cdot \wu D.$
\end{assertion}


\emph{Hint.} It suffices to prove the statement for the particular case when $\sigma$ is a transposition.

\begin{lemma}\label{l:wu-k5} 
The Wu--Simon number (of a plane diagram of $K_5$) is preserved under Reidemeister moves for plane diagrams of graphs (Figure~\ref{f:reid-1-5}).
\end{lemma}


Let $D$ be a plane diagram of $K_{3,3}$.
Consider the orientations on segments of~$D$ induced by the orientations $ij'$, $i,j \in [3]$, on edges of $K_{3,3}$.
At every crossing, consider the intersection sign, where the first vector passes above the second vector.
Color a crossing of non-adjacent edges~$e_1, e_2$ of~$K_{3,3}$ in \emph{blue}, if either

$\bullet$ $e_1 = ij'$ and $e_2 = ji'$ for some distinct $i,j \in [3]$, or 

$\bullet$ exactly one edge of $e_1, e_2$ belongs to $\{ 11', 22', 33' \}$.  

Color a crossing of non-adjacent edges~$e_1, e_2$ of~$K_{3,3}$  in \emph{red}, otherwise.

The {\bf bipartite Wu--Simon number} $\wu D$ is defined \cite[\S4]{Ta94} to be the difference of
the sum of signs at all the red crossings, and 
the sum of signs at all the blue crossings.

By a version for $K_{3,3}$ of the van Kampen--Flores theorem for the plane \cite[Remark 1.4.4.a]{Sk18},
any plane diagram of~$K_{3,3}$ has an odd number of crossings of non-adjacent edges, so the bipartite Wu--Simon number is odd. 
 
 
 
 
 
\begin{problem}\label{e:wu-k33}
(a) Find the bipartite Wu--Simon number of the plane diagram from Figure~\ref{f:ex-graph-emb} (middle), for your choice of numeration of vertices.
 
(b) For any odd 
$n$ there is a plane diagram of $K_{3,3}$ whose bipartite Wu--Simon number is $n$. 
\end{problem}

\begin{assertion}\label{con-wu-k33} 
The same as Assertion~\ref{con-wu-k5} replacing `$K_5$' by `$K_{3,3}$'.
\end{assertion}
 
\begin{assertion}\label{a:wu-k33-auto}* 
Let $D$, $D_{\sigma}$ be plane diagrams of $K_{3,3}$ such that $D_{\sigma}$ is obtained from $D$ by a correct renumeration $\sigma: [3] \cup [3]' \to [3] \cup [3]' $ of the vertices of $K_{3,3}$ (here `correct' means that $\sigma$ is induced by an automorphism of the graph $K_{3,3}$).
Then $\wu D_{\sigma} = \pm \wu D$.
\end{assertion}
 
\begin{lemma}\label{l:wu-k33} 
The bipartite Wu--Simon number (of a plane diagram of $K_{3,3}$) is preserved under Reidemeister moves for plane diagrams of graphs (Figure~\ref{f:reid-1-5}).
\end{lemma}
 
\begin{figure}[h]
\centerline{\includegraphics[width=3.5cm]{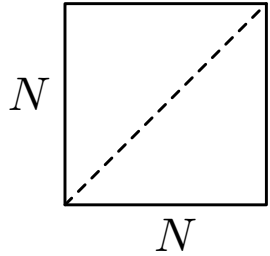}}
\caption{The deleted product}\label{f:dp}
\end{figure}

\begin{remark}[configuration spaces approach]\label{r:conf} 
In the 1950s Wu defined an invariant of higher-dimensional embeddings (of $k$-dimensional complexes to $\R^{2k+1}$). 
Let $K$ be a graph.   
The {\it deleted product} of $K$ is the product of~$K$ with itself, minus the diagonal (Figure \ref{f:dp}):
$\{(x,y)\in K\times K\ :\ x\ne y\}$. 
This is the configuration space of ordered pairs of distinct points in~$K$.
The {\it simplicial deleted product} of $K$ is
$$\widetilde K := \cup \{ \sigma\times\tau\ : \sigma, \tau \textrm{ are simplices of }K,\ \sigma\cap\tau = \emptyset\}.$$
Let $f:K\to\R^3$ be an embedding. 
The \emph{Wu invariant} 
$$\wu(f)\in H^2_{ss}(\widetilde K;\Z)$$   
assumes values in the \emph{group of cohomology classes of skew-symmetric 2-dimensional cocycles},  
see \cite[definition of $\mathcal{L}(f)$ in \S1]{Ta95}. 

Classical algebraic topology suggests the following standard way to obtain numerical invariants from a cohomological invariant. 
Fix orientations on edges of $K$. 
A \emph{2-dimensional chain} in $\widetilde K$ is an integer-valued function from the set of ordered pairs of disjoint 
edges of~$K$. 
A \emph{2-dimensional cycle}~$C$ in~$\widetilde K$ is 2-dimensional chain whose (homological) \emph{boundary} 
is zero. 
Examples of skew-symmetric 2-dimensional cycles are 
$$(*)\qquad K_3\times K'_3 + K'_3\times K_3,\quad \widetilde{K_5}\quad\text{and}\quad \widetilde{K_{3,3}}$$ 
(for certain orientations on the edges of the graphs $K_3 \sqcup K'_3$, $K_{5}$ and $K_{3,3}$, respectively). 
Cf. \cite{ADN+, MNS}. 
Denote by $H_2^{ss}(\widetilde K;\Z)$ the group of all skew-symmetric 2-dimensional cycle in $\widetilde K$. 

Define the \emph{scalar product}\footnote{This corresponds to $\frown$-product of the cohomology and homology groups of the quotient $\widetilde K$ with \emph{twisted coefficients}. 
Another name is \emph{the value of~$[A]$ at~$C$}. 
Cf. \cite[Remark 1.5.5c]{Sk18}.} 
$$\cdot:H^2_{ss}(\widetilde K;\Z)\times H_2^{ss}(\widetilde K;\Z)\to\Z\quad\text{by}\quad 
[A]\cdot C := \sum_{\{e_1,e_2\}} A(e_1,e_2)C(e_1,e_2),$$ 
where $A$ is a cocycle, $[A]$ its cohomology class, and the summation is taken over all unordered pairs of disjoint edges of $K$; for such an unordered pair the integers $A(e_1,e_2)$ and $C(e_1,e_2)$ are not well-defined, but their product is well-defined.  
Thus for a fixed~$C$ one defines an integer-valued invariant $\wu_C(f)\in\Z$ by $\wu_C(f):=\wu(f)\cdot C$. 

For cycles $Z$ in ($*$) the invariants $\wu_Z(f)$ are twice the linking number, the Wu-Simon number and the bipartite Wu-Simon number, respectively 
(the last two were independently discovered in \cite{Si86}). 
The invariants $\wu_C(f)$ are called \emph{reduced Wu invariants} or \emph{generalized Simon invariants} \cite[Definition 2.6]{FFN} 
 
The scalar product  $\cdot$ is clearly bilinear (this is essentially recovered in different terminology in \cite[Theorem 4.1]{FFN}). 
The scalar product  $\cdot$ is non-degenerate (this holds by \cite[Theorem 4.9]{Ta95}). 
I.e. for any homomorphism $\varphi:H^2_{ss}(\widetilde K;\Z)\to\Z$ there is $C\in H_2^{ss}(\widetilde K;\Z)$ such that $\varphi(\alpha)=\alpha\cdot C$ for any $\alpha\in H^2_{ss}(\widetilde K;\Z)$.  

\emph{Cycles $Z$, analogous to examples ($*$) for subgraphs of $K$ that are homeomorphic to
$K_3 \sqcup K'_3$, $K_{5}$, $K_{3,3}$, are rational generators in $H_2^{ss}(\widetilde K;\Z)$} \cite{Ni00}. 

Thus corresponding invariants $\wu_Z(f)$ are rational generators of all invariants $\wu_C(f)$
(this is essentially recovered in different terminology in \cite[Theorem 4.4]{FFN}). 

The paper \cite{FFN} used the reduced Wu invariants to obtain interesting \emph{chirality} and \emph{intrinsic chirality} results. 
The paper \cite{FFN} also contains some \emph{explicit calculations} of the invariants. 

Theorem \ref{t:tani} holds \cite{ST03} by the non-degeneracy of $\cdot$, the result of \cite{Ni00} quoted above, and the following result 
\cite[Main Theorem in \S1]{Ta95}: \emph{two embeddings of a graph in $\R^3$ are homologous if and only if their Wu invariants are equal}. 
\end{remark}

\textbf{Some answers, hints and solutions.} 

\begin{figure}[H]\centering
\includegraphics[scale=0.1]{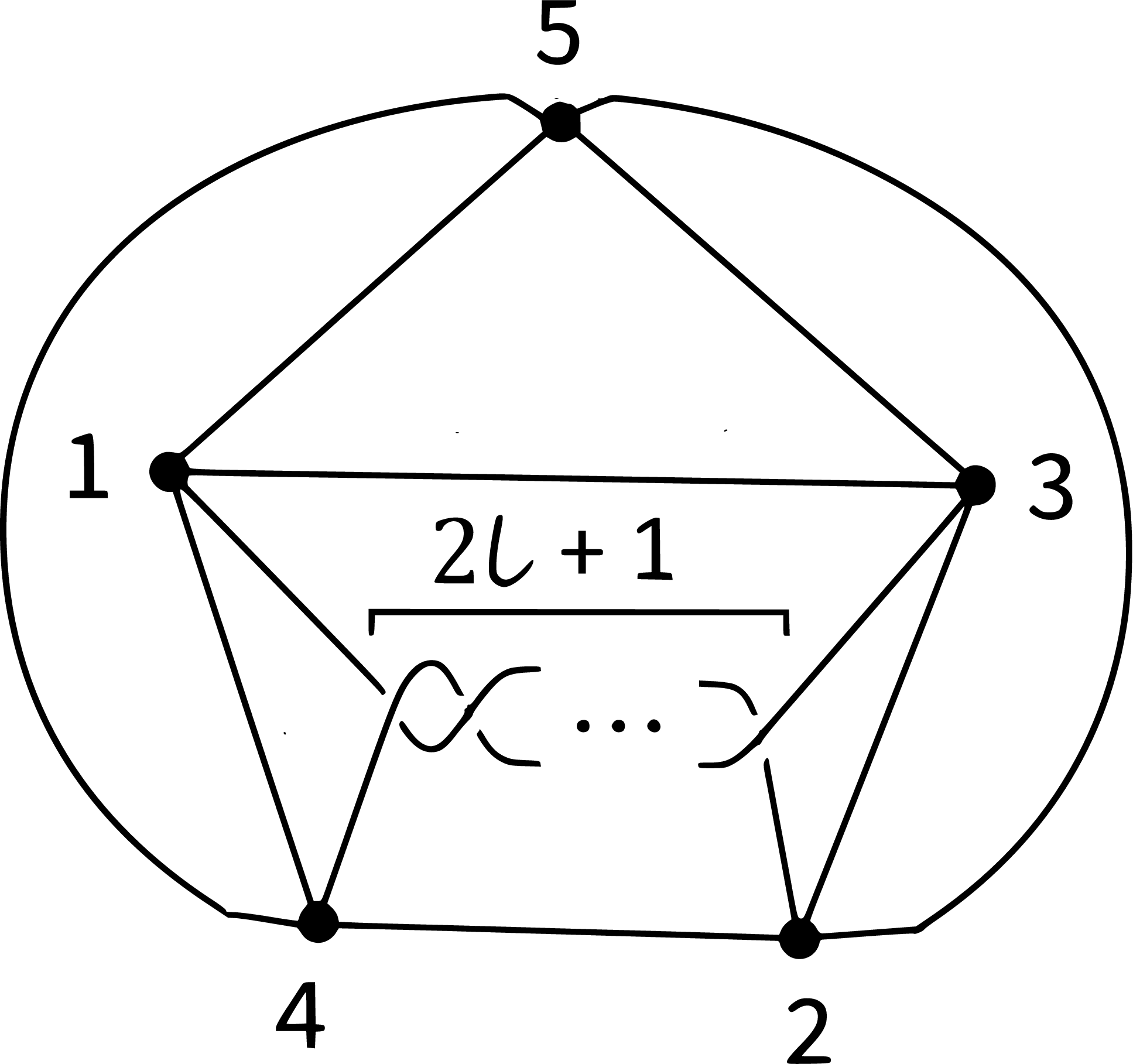}\qquad \qquad\includegraphics[scale=0.7]{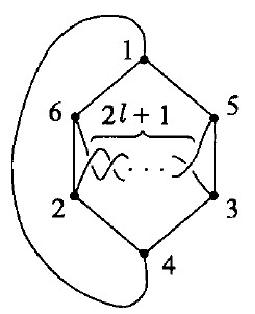}
\caption{Some plane diagrams of $K_5$ and of $K_{3,3}$}
\label{f:ex-k5}
\end{figure}

\smallskip
{\bf \ref{e:wu-k5}.} 
(a) \textit{Answer:} $\pm 5$; \qquad (b) As an example one can take the diagram from Figure~\ref{f:ex-k5}, left.

\smallskip
{\bf \ref{l:wu-k5}.}
We give a proof only for the second Reidemeister move. 
If the edges involved in this move are the same edge or are two adjacent edges, then the crossings involved do not affect the Wu--Simon number. 
If these edges are non-adjacent, then the crossings have the same colour but opposite signs, so their sum is zero, hence the Wu--Simon number does not change.

\smallskip
{\bf \ref{e:wu-k33}.} 
(a) \textit{Answer:} $\pm 3$; \qquad (b) As an example one can take the diagram from Figure~\ref{f:ex-k5}, right.

\smallskip
{\bf \ref{a:wu-k33-auto}.}
\emph{Hint.} 
This follows from the following two particular cases.

$\bullet$ If $\sigma(i) = i'$ and $\sigma(i') = i$ for any $i = 1, 2, 3$, then $\wu D_{\sigma} = \wu D$.
 
$\bullet$ If $\sigma$ is the transposition of some two vertices from one part of $K_{3,3}$, 
then $\wu D_{\sigma} = -\wu D$.

\section{Homology}\label{s:homol}



An \emph{extended diagram} of an embedding $K\to\R^3$ is the
union of a diagram of the embedding, and some oriented closed curves, every curve labelled by an edge of $K$
(different curves can be labelled by the same edge; the diagram and the union of the curves need not be contained in disjoint disks). 

\begin{figure}[H]\centering
\includegraphics[scale=1.2]{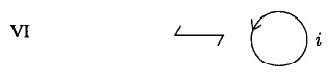}\\
\includegraphics[scale=1.2]{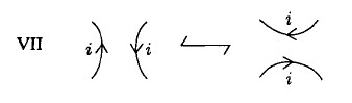}
\caption{Homology Reidemeister moves}
\label{f:reid-6-7}
\end{figure}



Two embeddings $f,g:K\to\R^3$ are \textbf{homologous} if some diagram of $f$ can be obtained from some diagram of $g$ by (isotopy) Reidemeister moves (Figure~\ref{f:reid-1-5}) and \emph{homology Reidemeister moves} (Figure~\ref{f:reid-6-7}) of extended diagrams 
(disappearance or appearance of an additional circle component labelled by any edge of $K$, and saddle move).
 
(Alternatively, in the definition of a \emph{concordance}, one is allowed to change open 2-disks in $K\times I$  by connected sum with tori.)


\begin{pr}\label{p:heqex} (a) Any knot is homologous to the unknot (trivial knot).

(b) The Borromean rings link is homologous to the trivial link. 

(c) Any two embeddings $K_4\to\R^3$ are homologous.
\end{pr}

\begin{pr}\label{p:hlk3} 
(a) There are non-homologous embeddings $K_3\sqcup K_3\to\R^3$.
(As an example one can take the standard and the Hopf links, Figure~\ref{f:hopf}.) 

(b) Two embeddings $K_3\sqcup K_3\to\R^3$ are homologous if and only if their linking numbers (defined in \S\ref{s:lkn}) are equal.
\end{pr}
 

\begin{pr}\label{p:hwu3} (a) There are non-homologous embeddings $K_5\to\R^3$.

(b) Two embeddings $K_5\to\R^3$ are homologous if and only if their Wu--Simon numbers (defined in \S\ref{s:wu}) are equal.

(a',b') Prove the analogues of (a,b) for $K_{3,3}$ and bipartite Wu--Simon number (defined in \S\ref{s:wu}). 
\end{pr}

{\it Books, surveys, and expository papers in this list are marked by the stars.}


\begin{thebibliography}{RSS95}

\UseRawInputEncoding

\newcommand{\aate}{\bibitem[AA38]{AA38} \emph{A. Adrian Albert}. Symmetric and alternate matrices in an arbitrary field, I. Trans. Amer. Math. Soc., (1938) 43(3):386--436.}

\newcommand{\abc}{\bibitem[ABC+]{ABC+} * \emph{M. Atiyah, A. Borel, G. J. Chaitin, D. Friedan, J. Glimm, J. J. Gray, M. W. Hirsch, S. MacLane, B. B. Mandelbrot, D. Ruelle, A. Schwarz, K. Uhlenbeck, R. Thom, E. Witten, C.  Zeeman.} Responses to ``Theoretical Mathematics: Toward a cultural synthesis of mathematics and theoretical physics'', by A. Jaffe and F. Quinn. Bull. Am. Math. Soc. 30 (1994) 178--207. arXiv:math/9404229.}

\newcommand{\abgmns}{\bibitem[ABM+]{ABM+} * \emph{E. Alkin, E. Bordacheva, A. Miroshnikov, O. Nikitenko, A. Skopenkov,} Invariants of almost embeddings of graphs in the plane: results and problems, arXiv:2408.06392.}

\newcommand{\abms}{\bibitem[ABM+]{ABM+} * \emph{Э. Алкин, Е. Бордачева, А. Мирошников, А. Скопенков,} Инварианты почти вложений графов в плоскость, arXiv:2410.09860.}

\newcommand{\adnt}{\bibitem[Ad93]{Ad93} * \emph{M. Adachi}. Embeddings and Immersions. Amer. Math.
Soc., 1993. (Transl. of Math. Monographs; V.~124).}

\newcommand{\adoe}{\bibitem[Ad18]{Ad18} {\it K. Adiprasito,} Combinatorial Lefschetz theorems beyond positivity, arXiv:1812.10454v4.}

\newcommand{\adnsv}{\bibitem[ADN+]{ADN+} * \emph{E. Alkin, S. Dzhenzher, O. Nikitenko, A. Skopenkov, A. Voropaev.}
Cycles in graphs and in hypergraphs: results and problems, arXiv:2308.05175.}

\newcommand{\agles}{\bibitem[AGL]{AGL86} Mathematical Economics,  ed. by A. Ambrosetti, F. Gori, R. Lucchetti,
Lect. Notes Math. 1330, Springer, 1986.}


\newcommand{\akzz}{\bibitem[Ak00]{Ak00} * \emph{П. М. Ахметьев.} Вложения компактов, стабильные
гомотопические группы сфер и теория особенностей, Успехи Мат. Наук.  2000. 55:3. C.~3-62.}

\newcommand{\akoe}{\bibitem[AK19]{AK19} \emph{S. Avvakumov, R. Karasev.} Envy-free division using mapping degree,
Mathematika, 67:1 (2020), 36--53. arXiv:1907.11183.}

\newcommand{\akto}{\bibitem[AK21]{AK21} \emph{G. Arone and V. Krushkal.}
Embedding obstructions in $\R^d$ from the Goodwillie-Weiss calculus and Whitney disks, 	Asian J. Math. 27 (2023), 135--186. arXiv:2101.10995. }

\newcommand{\akm}{\bibitem[AKM]{AKM} \emph{M. Abrahamsen, L. Kleist and T. Miltzow.}
Geometric Embeddability of Complexes is $\exists\mathbb R$-complete. arXiv:2108.02585.}

\newcommand{\aksoe}{\bibitem[AKS]{AKS} \emph{S. Avvakumov, R. Karasev and A. Skopenkov.} Stronger counterexamples to the topological Tverberg conjecture, Combinatorica, 43 (2023), 717--727. arXiv:1908.08731.}


\newcommand{\aksts}{\bibitem[AKS]{AKS} \emph{E. Alkin, Yu. Khromin, A. Skopenkov,} Linking invariants of spatial graphs, 
\linebreak 
\url{https://old.mccme.ru/circles/oim/graphs_in_space.pdf}}


\newcommand{\akuoe}{\bibitem[AKu19]{AKu19} \emph{S. Avvakumov, S. Kudrya.}
Vanishing of all equivariant obstructions and the mapping degree.
Discr. Comp. Geom., 66:3 (2021) 1202--1216. arXiv:1910.12628.}

\newcommand{\aktf}{\bibitem[Ak]{Ak} \emph{Д. Акимов,} Нечетность суммы чисел оборотов для почти вложения графа $K_{3,2}$.}

\newcommand{\alto}{\bibitem[Al22]{Al22} \emph{E. Alkin,}
Hardness of almost embedding simplicial complexes in $\R^d$, II. arXiv:2206.13486}

\newcommand{\amtf}{\bibitem[AM25]{AM25} \emph{E. Alkin, A. Miroshnikov,} On winding numbers of almost embeddings of $K_4$ in the plane, Math. Notes, to appear, arXiv:2501.15642.}

\newcommand{\amsw}{\bibitem[AMS+]{AMSW} \emph{S. Avvakumov, I. Mabillard, A. Skopenkov and U. Wagner.}
Eliminating Higher-Multiplicity Intersections, III. Codimension 2, Israel J. Math. 245 (2021) 501--534.  arxiv:1511.03501.}


\newcommand{\amse}{\bibitem[AMS]{AMS} * \emph{E. Alkin, A. Miroshnikov and A. Skopenkov,} 
Invariants of almost embeddings of graphs in the plane (Russian version), arXiv:2408.06392v2.}

\newcommand{\ams}{\bibitem[AMS]{AMS} * \emph{Э. Алкин, А. Мирошников, А. Скопенков,} Инварианты почти вложений графов в плоскость, arXiv:2410.09860v2.}

\newcommand{\amsp}{\bibitem[AMS']{AMS'} * \emph{Э. Алкин, А. Мирошников, А. Скопенков,} Число оборотов замкнутой ломаной вокруг точки, arXiv:2603.22351.}

\newcommand{\amspe}{\bibitem[AMS']{AMS'} * \emph{E. Alkin, A. Miroshnikov, A. Skopenkov,} The winding number of a closed curve around a point, arXiv:2603.22351.}


\newcommand{\anzt}{\bibitem[An03]{An03} * \emph{Д. В. Аносов.} Отображения окружности, векторные поля и их применения. М: МЦНМО, 2003.}

\newcommand{\anstf}{\bibitem[ANS]{ANS} \emph{E. Alkin, O. Nikitenko, A. Skopenkov,} Homotopy classification of closed polygonal lines: results and problems, arXiv:2508.16287.}

\newcommand{\arnf}{\bibitem[Ar95]{Ar95} * \emph{V. I. Arnold,}  Topological invariants of plane curves and caustics, University Lecture Series, Vol. 5, Amer. Math. Soc., Providence, RI, 1995.}

\newcommand{\arszo}{\bibitem[ARS01]{ARS01} \emph{P. Akhmetiev, D. Repov\v s and A. Skopenkov},
Embedding products of low-dimensional manifolds in $\R^m$, Topol. Appl. 113 (2001), 7--12.}

\newcommand{\arszt}{\bibitem[ARS02]{ARS02} \emph{P. Akhmetiev, D. Repovs and A. Skopenkov.} Obstructions to approximating maps of $n$-manifolds into $R^{2n}$ by embeddings, Topol. Appl., 123 (2002), 3--14.}

\newcommand{\asts}{\bibitem[AS25]{AS25} \emph{E. Alkin and A. Skopenkov,} Homotopy classification of closed polygonal lines, arXiv:2508.16287, Russian version.}

\newcommand{\asoed}{\bibitem[As]{As} \emph{A. Asanau,} \lowercase{A SIMPLE PROOF THAT CONNECTED SUM OF ORDERED
ORIENTED LINKS IS NOT WELL-DEFINED,} Math. Notes, to appear.}

\newcommand{\asoe}{\bibitem[As]{As} \emph{A. Asanau,} On the \lowercase{TRIPLE SELF-INTERSECTION NUMBER FOR GRAPHS IN THE PLANE,} unpublished, 2018.}

\newcommand{\avos}{\bibitem[Av14]{Av14} \emph{S. Avvakumov,} The classification of certain linked 3-manifolds in 6-space, Moscow Math. J., 16:1 (2016), 1--25. arXiv:1408.3918.}

\newcommand{\avose}{\bibitem[Av17]{Av17} \emph{S. Avvakumov,} The classification of linked 3-manifolds in 6-space, Algebraic \& Geometric Topology, 22:6 (2022) 2587--2630. arXiv:1704.06501.}



\newcommand{\bant}{\bibitem[Ba93]{Ba93} * \emph{T. Bartsch.} Topological methods for variational problems with
symmetries, Lecture Notes in Mathematics, 1560, Springer-Verlag, Berlin, 1993.}

\newcommand{\batt}{\bibitem[Ba23]{Ba23} * \emph{I. Barany.} Tverberg's theorem, a new proof. arXiv:2308.10105.}

\newcommand{\bbsn}{\bibitem[BB79]{BB} \emph{E.~G. Bajm{{\'o}}czy and I.~B{{\'a}}r{{\'a}}ny,}
\newblock On a common generalization of {B}orsuk's and {R}adon's theorem,
\newblock Acta Math.\ Acad.\ Sci.\ Hungar.\ 34:3 (1979), 347-350.}

\newcommand{\bbzos}{\bibitem[BBZ]{BBZ} * \emph{I.~B{{\'a}}r{{\'a}}ny, P.~V.~M. Blagojevi{{\'c}} and G.~M. Ziegler.} Tverberg's Theorem at 50: Extensions and Counterexamples, Notices of the Amer. Math. Soc., 63:7 (2016), 732--739.}


\newcommand{\bcm}{\bibitem[BCM]{BCM} * 13th Hilbert Problem on superpositions of functions, presented by A. Belov, A. Chilikov, I. Mitrofanov, S. Shaposhnikov and A. Skopenkov,
\url{http://www.turgor.ru/lktg/2016/5/index.htm}.}

\newcommand{\beet}{\bibitem[BE82]{BE82} * \emph{V.G. Boltyansky and V.A. Efremovich.} Intuitive Combinatorial Topology. Springer.}

\newcommand{\beetr}{\bibitem[BE82]{BE82} * \emph{В. Г. Болтянский и В. А. Ефремович.} Наглядная топология. М.:  Наука, 1982.}


\newcommand{\bfzn}{\bibitem[BF09]{BF09} \emph{K. Barnett, M. Farber}. Topology of Configuration Space of Two Particles on a Graph, I.  Algebr. Geom. Topol. 9 (2009) 593--624.	arXiv:0903.2180.}

\newcommand{\bfzof}{\bibitem[BFZ14]{BFZ14} \emph{P. V. M. Blagojevi{\'c}, F. Frick, and G. M. Ziegler,}
Tverberg plus constraints, Bull. Lond. Math. Soc. 46:5 (2014), 953-967, arXiv:1401.0690.}


\newcommand{\bfzos}{\bibitem[BFZ]{BFZ} \emph{P. V. M. Blagojevi{\'c}, F. Frick and G. M. Ziegler,}
Barycenters of Polytope Skeleta and Counterexamples to the Topological Tverberg Conjecture, via Constraints,
J. Eur. Math. Soc., 21:7 (2019) 2107-2116. arXiv:1510.07984.}


\newcommand{\bgso}{\bibitem[BG71]{BG71} J.C. Becker and H. H. Glover, {\it Note on the Embedding of Manifolds in Euclidean Space,} Proc. of the Amer. Math. Soc., 27:2 (1971) 405-410.}


\newcommand{\bgos}{\bibitem[BG16]{BG16} \emph{A. Bj\"orner and A. Goodarzi}, On Codimension one Embedding of Simplicial Complexes, in book: A Journey Through Discrete Mathematics, arXiv:1605.01240.}

\newcommand{\biet}{\bibitem[Bi83]{Bi83} * \emph{R. H. Bing.} The Geometric Topology of 3-Manifolds. Providence, R.~I. 1983. (Amer. Math. Soc. Colloq. Publ., 40).}

\newcommand{\bitz}{\bibitem[Bi20]{Bi20} * \emph{A. Bikeev.} Realizability of discs with ribbons on the M\"obius strip. Mat. Prosveschenie, 28 (2021), 150-158;
erratum to appear. arXiv:2010.15833.}

\newcommand{\bitzr}{\bibitem[Bi20]{Bi20} * \emph{А. Бикеев.} Реализуемость дисков с ленточками на ленте Мебиуса.
Мат. просвещение. Сер. 3. 28 (2021), 150--158.}

\newcommand{\bito}{\bibitem[Bi21]{Bi21} {\it A. I. Bikeev,}
Criteria for integer and modulo 2 embeddability of graphs to surfaces, arXiv:2012.12070v2.}


\newcommand{\bagos}{\bibitem[BG17]{BG17} \emph{S. Basu and S. Ghosh.} Equivariant maps related to the topological Tverberg conjecture, Homology, Homotopy and Applications 19:1 (2017) 155--170.}

\newcommand{\bkkmzof}{\bibitem[BKK]{BKK} \emph{M. Bestvina, M. Kapovich and B. Kleiner,}
Van Kampen's embedding obstruction for discrete groups, Invent. Math. 150 (2002) 219--235. arXiv:math/0010141.}

\newcommand{\bl}{\bibitem[BL]{BL} \url{https://en.wikipedia.org/wiki/Brunnian_link}}

\newcommand{\blf}{\bibitem[BL4]{BL4} Students form a 4-component Brunnian link,  \url{http://www.mccme.ru/circles/oim/foto2014/brunn4.png} (5Mb)}

\newcommand{\bmzf}{\bibitem[BM04]{BM04} \emph{Boyer, J. M. and Myrvold, W. J.} On the cutting edge: simplified $O(n)$ planarity by edge addition,  Journal of Graph Algorithms and Applications, 8:3 (2004) 241--273.}

\newcommand{\bm}{\bibitem[BM15]{BM15} \emph{I. Bogdanov and A. Matushkin.} Algebraic proofs of linear versions of the Conway--Gordon--Sachs theorem and the van Kampen--Flores theorem, arXiv:1508.03185.}


\newcommand{\bmzzn}{\bibitem[BMZ09]{BMZ09} \emph{P. V. M. Blagojevi{\'c}, B. Matschke, G. M. Ziegler,}
Optimal bounds for a colorful Tverberg-Vre\'cica type problem, Advances in Math., 226 (2011), 5198-5215, arXiv:0911.2692.}

\newcommand{\bmzof}{\bibitem[BMZ15]{BMZ15} \emph{P. V. M. Blagojevi{\'c}, B. Matschke, G. M. Ziegler,}
Optimal bounds for the colored Tverberg problem, J. Eur. Math. Soc.,  17:4 (2015) 739--754,
arXiv:0910.4987.}

\newcommand{\bpns}{\bibitem[BP97]{BP97} * \emph{R. Benedetti and C. Petronio.} Branched standard spines of 3-manifolds, Lecture Notes in Math. 1653, Springer-Verlag, Berlin-Heidelberg-New York, 1997.}

\newcommand{\brst}{\bibitem[Br72]{Br72} \emph{J. L. Bryant.} Approximating embeddings of polyhedra in codimension 3, Trans. Amer. Math. Soc., 170 (1972) 85--95.}

\newcommand{\brts}{\bibitem[Br68]{Br68} \emph{P. Bruegel,} 1568,
\url{https://en.wikipedia.org/wiki/The_Magpie_on_the_Gallows}.}


\newcommand{\bren}{\bibitem[Br82]{brown1982} * \emph{K.~S. Brown.} \newblock Cohomology of Groups. \newblock Springer-Verlag New York, 1982.}


\newcommand{\bssos}{\bibitem[BS17]{BS17} * \emph{I.~B\'{a}r\'{a}ny and P. Sober\'{o}n,} Tverberg's theorem is 50 years old: a survey, Bull. Amer. Math. Soc. (N.S.) 55:4 (2018), 459--492. arXiv:1712.06119.}

\newcommand{\bsto}{\bibitem[BS21]{BS21} * \emph{A. Buchaev and A. Skopenkov,} Simple proofs of estimations of Ramsey numbers and of discrepancy, Mat. Prosveschenie, to appear, arXiv:2107.13831.}

\newcommand{\brsnn}{\bibitem[BRS99]{BRS99} \emph{D. Repov\v s, N. Brodsky and A. B. Skopenkov.}
A classification of 3-thickenings of 2-polyhedra, Topol. Appl. 1999. 94. P.~307-314.}

\newcommand{\bsseo}{\bibitem[BSS]{BSS} \emph{I.~B\'{a}r\'{a}ny, S.~B. Shlosman, and A.~Sz{\H{u}}cs,}
\newblock On a topological generalization of a theorem of {T}verberg,
\newblock J.\ London Math.\ Soc.\ (II. Ser.) 23 (1981), 158--164.}

\newcommand{\btzs}{\bibitem[BT07]{BT07} \emph{A. Bj\"orner, M. Tancer}, Combinatorial Alexander Duality --- a Short and Elementary Proof, Discr. and Comp. Geom., 42 (2009) 586. arXiv:0710.1172.}

\newcommand{\buse}{\bibitem[Bu68]{Bu68} \emph{A. R. Butz,} Space filling curves and mathematical programming, Information and Control, 12:4 (1968) 314--330.}


\newcommand{\bz}{\bibitem[BZ16]{BZ16} * \emph{P. V. M. Blagojevi\'c and G. M. Ziegler,} Beyond the Borsuk-Ulam theorem: The topological Tverberg story, in: A Journey Through Discrete Mathematics, Eds. M. Loebl,
J. Ne\v set\v ril, R. Thomas, Springer, 2017, 273--341. arXiv:1605.07321v3.}



\newcommand{\cano}{\bibitem[Ca91]{Ca91} * \emph{D. de Caen}, The ranks of tournament matrices, Amer. Math. Monthly, 98:9 (1991) 829--831.}

\newcommand{\ca}{\bibitem[Ca]{Ca} \emph{J. Carmesin.} Embedding simply connected 2-complexes in 3-space, I-V, arXiv:1709.04642, arXiv:1709.04643, arXiv:1709.04645, arXiv:1709.04652, arXiv:1709.04659.}

\newcommand{\cfsz}{\bibitem[CF60]{CF60} \emph{P. E. Conner and E. E. Floyd}, Fixed points free involutions and equivariant maps, Bull. Amer. Math. Soc., 66 (1960) 416--441.}

\newcommand{\cfs}{\bibitem[CFS]{CFS} \emph{D. Crowley, S.C. Ferry, M. Skopenkov,} The rational classification of links of codimension $>2$, Forum Math. 26 (2014), 239--269. arXiv:1106.1455.}

\newcommand{\cget}{\bibitem[CG83]{CG83} \emph{J. H. Conway and C. M. A. Gordon},
Knots and links in spatial graphs, J. Graph Theory  7 (1983), 445--453.}

\newcommand{\cten}{\bibitem[Ch]{Ch} \emph{Chuang Tzu,} translated by H. A. Giles, Bernard Quaritch, London, 1889.}

\newcommand{\ctruku}{\bibitem[Ch]{Ch} \emph{Chuang Tzu,} translated to Russian by S. Kuchera, in: Ancient Chinese Philosophy, v. I, Mysl, Moscow, 1972.}


\newcommand{\chnn}{\bibitem[Ch99]{Ch99} * \emph{А. В. Чернавский,} Теорема Жордана.  Мат. Просвещение, 3 (1999), 142--157.}

\newcommand{\hcon}{\bibitem[HC19]{HC19} * \emph{C. Herbert Clemens.} Two-Dimensional Geometries. A Problem-Solving Approach, Amer. Math. Soc., 2019.}

\newcommand{\ckmoo}{\bibitem[CKMS]{CKMS} \emph{M. \v Cadek, M. Kr\v c\'al. J. Matou\v sek, F. Sergeraert,
L. Vok\v r\'inek, U. Wagner.} Computing all maps into a sphere, J. of the ACM, 61:3 (2014). arXiv:1105.6257.}


\newcommand{\ckmvwot}{\bibitem[CKM12+]{CKM12+} \emph{M. \v Cadek, M. Kr\v c\'al. J. Matou\v sek, L. Vok\v r\'inek, U. Wagner.} Polynomial-time computation of homotopy groups and Postnikov systems in fixed dimension, SIAM J. Comput., 43:5 (2014), 1728--1780. arXiv:1211.3093.}

\newcommand{\ckmvw}{\bibitem[CKM+]{CKM+} \emph{M. \v Cadek, M. Kr\v c\'al. J. Matou\v sek, L. Vok\v r\'inek, U. Wagner.} Extendability of continuous maps is undecidable, Discr. and Comp. Geom. 51 (2014) 24--66.
arXiv:1302.2370.}

\newcommand{\ckppt}{\bibitem[CKP+]{CKP+} \emph{E. Colin de Verdi\'ere, V. Kalu\v za, P. Pat\'ak, Z. Pat\'akov\'a and M. Tancer.} A direct proof of the strong Hanani-Tutte theorem on the projective plane. Journal of Graph Algorithms and Applications, 21:5 (2017) 939--981.}

\newcommand{\cksof}{\bibitem[CKS+]{CKS+} * New ways of weaving baskets, presented by G. Chelnokov, Yu. Kudryashov, A.Skopenkov and A. Sossinsky, \url{http://www.turgor.ru/lktg/2004/lines.en/index.htm}.}

\newcommand{\ckv}{\bibitem[CKV]{CKV} \emph{M.~{\v{C}}adek, M.~Kr\v{c}\'{a}l, and L.~Vok\v{r}\'{\i}nek.}
Algorithmic solvability of the lifting-extension problem, Discr. Comp. Geom. 57 (2017), 915--965. arXiv:1307.6444.}


\newcommand{\clr}{\bibitem[CLR]{CLR} * \emph{Т. Кормен, Ч. Лейзерсон, Р. Ривест.} Алгоритмы:
построение и анализ, МЦНМО, Москва, 1999.}

\newcommand{\clreng}{\bibitem[CLR]{CLR} * \emph{T. H. Cormen, C. E.Leiserson, R. L.Rivest, C. Stein.} Introduction to Algorithms, MIT Press, 2009.}

\newcommand{\crzfru}{\bibitem[CR]{CR} * \emph{Р. Курант, Дж. Роббинс,} Что такое математика. М.: МЦНМО, 2004.}

\newcommand{\crzfen}{\bibitem[CR]{CR} * \emph{R. Courant and H. Robbins,} What is Mathematics, Oxford Univ. Press.}

\newcommand{\crsne}{\bibitem[CRS98]{CRS98} * \emph{A. Cavicchioli, D. Repov\v s and A. B. Skopenkov.}
Open problems on graphs, arising from geometric topology, Topol. Appl. 84 (1998), 207--226.}

\newcommand{\crsos}{\bibitem[CRS']{CRS'} \emph{M. Cencelj, D. Repov\v s and M. Skopenkov,} Homotopy type of the complement of an immersion and classification of embeddings of tori, Russian Math. Surv. 62:5 (2007), 985--987,
arXiv:0803.4285.}

\newcommand{\crsot}{\bibitem[CRS]{CRS} \emph{M. Cencelj, D. Repov\v s and M. Skopenkov,}
Classification of knotted tori in the 2-metastable dimension, Mat. Sbornik, 203:11 (2012), 1654--1681.
arxiv:0811.2745.}

\newcommand{\csoo}{\bibitem[CS08]{CS08} \emph{D. Crowley and A. Skopenkov.} A classification of smooth embeddings of 4-manifolds in 7-space, II, Intern. J. Math., 22:6 (2011) 731-757, arxiv:0808.1795.}

\newcommand{\csos}{\bibitem[CS16]{CS16} \emph{D. Crowley and A. Skopenkov,} Embeddings of non-simply-connected 4-manifolds in 7-space. I. Classification modulo knots, Moscow Math. J., 21 (2021), 43--98. arXiv:1611.04738.}


\newcommand{\csoso}{\bibitem[CS16o]{CS16o} \emph{D. Crowley and A. Skopenkov,} Embeddings of non-simply-connected 4-manifolds in 7-space. II. On the smooth classification, Proc. A of the Royal Soc. of Edinburgh 152:1 (2022), 163--181. arXiv:1612.04776.}


\newcommand{\crsk}{\bibitem[CS]{CS} \emph{D. Crowley and A. Skopenkov,} Embeddings of non-simply-connected 4-manifolds in 7-space. III. Piecewise-linear classification. draft.}

\newcommand{\cutz}{\bibitem[Cu20]{Cu20} \emph{C. Culter,} Cantor sets are not tangent homogeneous,
Topol. Appl. 271 (2020) 1--9.}


\newcommand{\dies}{\bibitem[Di87]{Di} * \emph{T. tom Dieck,} Transformation groups, Studies in Mathematics, vol. 8, Walter de Gruyter, Berlin, 1987.}

\newcommand{\dize}{\bibitem[Di08]{Di08} * \emph{T. tom Dieck,} Algebraic topology, EMS Textbooks in Mathematics, 
EMS, Z\"urich, 2008.}

\newcommand{\dent}{\bibitem[De93]{De93}  \emph{T.K. Dey.} On counting triangulations in $d$-dimensions. Comput. Geom.  3:6 (1993) 315--325.}

\newcommand{\denf}{\bibitem[DE94]{DE94}  \emph{T.K. Dey and H. Edelsbrunner.} Counting triangle crossings and halving planes, Discrete Comput. Geom. 12 (1994), 281--289.}

\newcommand{\dgn}{\bibitem[DGN+]{DGN+} * S. Dzhenzher, T. Garaev, O. Nikitenko, A. Petukhov, A. Skopenkov, A. Voropaev, Low rank matrix completion and realization of graphs: results and problems, arXiv:2501.13935.}

\newcommand{\dgnr}{\bibitem[DGN+]{DGN+} * Минимизация ранга восполнением матриц, представляли А. Воропаев, Т. Гараев, С. Дженжер, О. Никитенко, А. Петухов и А. Скопенков, \url{https://www.mccme.ru/circles/oim/netflix_rus.pdf}.}

\newcommand{\dstt}{\bibitem[DS22]{DS22}  \emph{S. Dzhenzher and A. Skopenkov,} A quadratic estimation for the K\"uhnel conjecture on embeddings, arXiv:2208.04188.}

\newcommand{\botf}{\bibitem[Dz25]{Dz25} \emph{E. Dzhenzher,} Symmetric 1-cycles in the deleted product of a graph, Topol. Appl. (2025) 109277.}



\newcommand{\embo}{\bibitem[Eb]{Eb} * \url{http://www.map.mpim-bonn.mpg.de/Embeddings_of_manifolds_with_boundary:_classification}}

\newcommand{\embe}{\bibitem[Em]{Em} * \url{http://www.map.mpim-bonn.mpg.de/Embedding_(simple_definition)}}

\newcommand{\ers}{\bibitem[ERS]{ERS} * Invariants of graph drawings in the plane, presented by A. Enne, A. Ryabichev, A. Skopenkov and T. Zaitsev, \url{http://www.turgor.ru/lktg/2017/6/index.htm}}



\newcommand{\feto}{\bibitem[Fe21]{Fe21} \emph{M. Fedorov.} A description of values of Seifert form for punctured $n$-manifolds in $(2n-1)$-space, arXiv:2107.02541.}

\newcommand{\ffen}{\bibitem[FF89]{FF89} * \emph{А. Т. Фоменко и Д. Б. Фукс.} Курс гомотопической топологии. М.: Наука, 1989.}

\newcommand{\ffene}{\bibitem[FF89]{FF89} * \emph{A.T. Fomenko and D.B. Fuchs.} Homotopical Topology, Springer, 2016.}


\newcommand{\fhzo}{\bibitem[FH10]{FH10}  \emph{M. Farber, E. Hanbury}. Topology of Configuration Space of Two Particles on a Graph, II. Algebr. Geom. Topol. 10 (2010) 2203--2227. arXiv:1005.2300.}


\newcommand{\fkosc}{\bibitem[FK17]{FK17} \emph{R. Fulek, J. Kyn{\v{c}}l,} Counterexample to an Extension of the Hanani-Tutte Theorem on the Surface of Genus 4, Combinatorica, 39 (2019) 1267--1279, arXiv:1709.00508.}

\newcommand{\fkos}{\bibitem[FK17]{FK17} \emph{R. Fulek, J. Kyn{\v{c}}l,} Hanani-Tutte for approximating maps of graphs, arXiv:1705.05243.}

\newcommand{\fkon}{\bibitem[FK19]{FK19} \emph{R. Fulek, J. Kyn{\v{c}}l,}
$\Z_2$-genus of graphs and minimum rank of partial symmetric matrices,
35th Intern. Symp. on Comp. Geom. (SoCG 2019), Article No. 39; pp. 39:1--39:16, \linebreak
\url{https://drops.dagstuhl.de/opus/volltexte/2019/10443/pdf/LIPIcs-SoCG-2019-39.pdf}.
We refer to numbering in arXiv version: arXiv:1903.08637.}

\newcommand{\fktnf}{\bibitem[FKT]{FKT} \emph{M. H. Freedman, V. S. Krushkal and P. Teichner.} Van Kampen's
embedding obstruction is incomplete for 2-complexes in~$\R^4$, Math. Res. Letters. 1994. 1. P.~167-176.}

\newcommand{\fltf}{\bibitem[Fl34]{Fl34} \emph{A. Flores}, \"Uber $n$-dimensionale Komplexe die im $E^{2n+1}$ absolut selbstverschlungen sind, Ergeb. Math. Koll. 6 (1934) 4--7.}

\newcommand{\fnzn}{\bibitem[FN09]{FN09} \emph{T. Fleming and R. Nikkuni,} Homotopy on spatial graphs and the Sato-Levine invariant, Trans. Amer. Math. Soc. 361:4 (2009), 1885--1902.}

\newcommand{\fo}{\bibitem[Fo]{Fo} * \emph{L. Fortnow.} Time for Computer Science to Grow Up,  \url{https://people.cs.uchicago.edu/~fortnow/papers/growup.pdf}.}

\newcommand{\fozf}{\bibitem[Fo04]{Fo04} * \emph{R. Fokkink.} A forgotten mathematician, Eur. Math. Soc. Newsletter 52 (2004) 9--14.}


\newcommand{\fpstz}{\bibitem[FPS]{FPS} \emph{R. Fulek, M.J. Pelsmajer and M. Schaefer.}
Strong Hanani-Tutte for the Torus, arXiv:2009.01683.}

\newcommand{\frse}{\bibitem[Fr78]{Fr78} \emph{M. Freedman,} Quadruple points of 3-manifolds in $S^4$, Comment. Math. Helv. 53 (1978), 385-394.}

\newcommand{\fres}{\bibitem[FR86]{FR86} \emph{R. Fenn, D. Rolfsen.}
Spheres may link homotopically in 4-space, J. London Math. Soc. 34 (1986) 177-184.}

\newcommand{\frofea}{\bibitem[Fr15']{Fr15'} \emph{F. Frick}, Counterexamples to the topological Tverberg conjecture, arXiv:1502.00947v1.}


\newcommand{\frof}{\bibitem[Fr15]{Fr15} \emph{F. Frick}, Counterexamples to the topological Tverberg conjecture,
Oberwolfach reports, 12:1 (2015), 318--321. arXiv:1502.00947.}

\newcommand{\fros}{\bibitem[Fr17]{Fr17} \emph{F. Frick}, O\lowercase{N AFFINE TVERBERG-TYPE RESULTS WITHOUT CONTINUOUS GENERALIZATION}, arXiv:1702.05466}


\newcommand{\fstz}{\bibitem[FS20]{FS20} \emph{F. Frick and P. Sober\'on}, The topological Tverberg problem beyond prime powers, arXiv:2005.05251.}

\newcommand{\ftss}{\bibitem[FT77]{FT77} \emph{R. Fenn, P. Taylor,} Introducing doodles, pp. 37-43
in: Topology of Low-Dimensional Manifolds, Proceedings of the Second Sussex Conference, 1977,
Ed. R. Fenn, V. 722 of Lecture Notes in Math.}

\newcommand{\fvto}{\bibitem[FV21]{FV21} \emph{M. Filakovsk\'y, L. Vok\v r\'inek.} Computing homotopy classes for diagrams, Discr. Comp. Geom. 70 (2023), 866--920. arXiv:2104.10152.}

\newcommand{\fwz}{\bibitem[FWZ]{FWZ} \emph{M. Filakovsk\'y, U. Wagner, S. Zhechev.} Embeddability of simplicial complexes is undecidable. Oberwolfach reports, to appear.}

\newcommand{\fwztz}{\bibitem[FWZ]{FWZ} \emph{M. Filakovsk\'y, U. Wagner, S. Zhechev.} Embeddability of simplicial complexes is undecidable. Proceedings of the 2020 ACM-SIAM Symposium on Discrete Algorithms.}



\newcommand{\ga}{\bibitem[GA]{GA} * \url{https://en.wikipedia.org/wiki/Galactic_algorithm}}

\newcommand{\gatt}{\bibitem[Ga23]{Ga23} \emph{T. Garaev}, On drawing $K_5$ minus an edge in the plane, arXiv:2303.14503.}

\newcommand{\gdikrse}{\bibitem[GDI]{GDI} * {\it A. Chernov, A. Daynyak, A. Glibichuk, M. Ilyinskiy, A. Kupavskiy, A. Raigorodskiy and A. Skopenkov,} Elements of Discrete Mathematics As a Sequence of Problems (in Russian),
MCCME, Moscow, 2016. Update of a part: \url{http://www.mccme.ru/circles/oim/discrbook.pdf}}

\newcommand{\gdikrs}{\bibitem[GDI]{GDI} * {\it А.А. Глибичук, А.Б. Дайняк, Д.Г. Ильинский, А.Б. Купавский, А.М. Райгородский, А.Б. Скопенков, А.А. Чернов,} Элементы дискретной математики в задачах, М, МЦНМО, 2016.
Обновляемая версия части книги: \url{http://www.mccme.ru/circles/oim/discrbook.pdf}}

\newcommand{\giso}{\bibitem[Gi71]{Gi71} * {\it S. Gitler,} Immersion and Embedding of Manifolds,
Proc. Symp. Pure Math. 22, 87-96 (1971).}

\newcommand{\gkp}{\bibitem[GKP]{GKP} * {\it R. Graham, D. Knuth, and O. Patashnik,} Concrete Mathematics: A Foundation for Computer Science, Addison–Wesley, first published in 1989, \url{https://www.csie.ntu.edu.tw/~r97002/temp/Concrete\%20Mathematics\%202e.pdf}.}

\newcommand{\gmpptw}{\bibitem[GMP+]{GMP+} \emph{X. Goaoc, I. Mabillard, P. Pat\'ak, Z. Pat\'akov\'a, M. Tancer, U. Wagner}, On Generalized Heawood Inequalities for Manifolds: a van Kampen--Flores-type Nonembeddability Result,
Israel J. Math., 222(2) (2017) 841-866. arXiv:1610.09063.}


\newcommand{\gppot}{\bibitem[GPP+]{GPP+} \emph{X. Goaoc, P. Pat\'ak, Z. Pat\'akov\'a, M. Tancer, and U. Wagner.} Bounding Helly numbers via Betti numbers. In 31st International Symposium on Computational Geometry, volume 34
of LIPIcs. Leibniz Int. Proc. Inform., pp. 507-521. Schloss Dagstuhl. Leibniz-Zent. Inform., Wadern, 2015. Full version: arXiv:1310.4613.}

\newcommand{\group}{\bibitem[Gr]{Gr} * \url{https://en.wikipedia.org/wiki/Groupthink}}

\newcommand{\grsz}{\bibitem[Gr69]{Gr69} \emph{B. Gr\"unbaum.} Imbeddings of simplicial complexes. Comment. Math. Helv., 44:1, 502--513, 1969.}


\newcommand{\gres}{\bibitem[Gr86]{Gr86} * \emph{M. Gromov}, Partial Differential Relations,
Ergebnisse der Mathematik und ihrer Grenzgebiete (3), Springer Verlag, Berlin-New York, 1986.}

\newcommand{\groz}{\bibitem[Gr10]{Gr10} \emph{M. Gromov,}
\newblock Singularities, expanders and topology of maps. Part 2: From combinatorics to topology via algebraic isoperimetry, \newblock Geometric and Functional Analysis 20 (2010), no.~2, 416--526.}

\newcommand{\grsn}{\bibitem[GR79]{GR79} \emph{J. L. Gross	and R. H. Rosen}, A linear time planarity algorithm for 2-complexes, Journal of the ACM, 26:4 (1979), 611--617.}

\newcommand{\gs}{\bibitem[GS]{GS} \emph{М. Гортинский и О. Скрябин.} Критерий вложимости графов в плоскость вдоль прямой, препринт.}

\newcommand{\gssn}{\bibitem[GS79]{GS} \emph{P.~M. Gruber and R.~Schneider,} Problems in geometric convexity. In {\em Contributions to geometry (Proc. Geom. Sympos., Siegen, 1978)}, 255--278. Birkh{\"a}user, Basel-Boston, Mass., 1979.}

\newcommand{\gsnn}{\bibitem[GS99]{GS99} \emph{R. Gompf and A. Stipsicz,}
4-manifolds and Kirby calculus, GSM20, AMS, Providence, RI, 1999.}


\newcommand{\gszs}{\bibitem[GS06]{GS06} \emph{D. Goncalves and A. Skopenkov,} Embeddings of homology equivalent manifolds with boundary, Topol. Appl., 153:12 (2006) 2026-2034. arxiv:1207.1326.}

\newcommand{\gssoe}{\bibitem[GSS+]{GSS+} * Projections of skew lines, presented by A. Gaifullin, A. Shapovalov, A. Skopenkov and M. Skopenkov, \url{http://www.turgor.ru/lktg/2001/index.php}.}

\newcommand{\gtes}{\bibitem[GT87]{GT87} * \emph{J. L. Gross and T. W. Tucker.}
Topological graph theory. New York: Wiley-Interscience, 1987.}

\newcommand{\guzn}{\bibitem[Gu09]{Gu09} \emph{A. Gundert.} On the complexity of embeddable simplicial complexes. Diplomarbeit, Freie Universit\"at Berlin, 2009. 	arXiv:1812.08447.}


\newcommand{\ha}{\bibitem[Ha]{Ha} * \emph{F. Harary.} Graph theory.
Рус. пер.: Ф. Харари. Теория графов. М., Мир, 1973.}

\newcommand{\hats}{\bibitem[Ha37]{Ha37} \emph{W. Hantzsche,} Einlagerung von Mannigfaltigkeiten in euklidische R\" aume, Math. Zeitschrift, 43:1 (1937) 38--58.}

\newcommand{\hastk}{\bibitem[Ha62k]{Ha62k} {\em A.~Haefliger,}  Knotted $(4k-1)$-spheres in $6k$-space, Ann. of Math. 75 (1962) 452--466.}

\newcommand{\hastl}{\bibitem[Ha62l]{Ha62l} \emph{A. Haefliger,} Differentiable links, Topology, 1 (1962) 241--244.}

\newcommand{\hast}{\bibitem[Ha63]{Ha63} \emph{A.~Haefliger,} Plongements differentiables dans le domain stable, Comment. Math. Helv. 36 (1962-63) 155--176.}

\newcommand{\hassa}{\bibitem[Ha66A]{Ha66A} \textit{A. Haefliger}. Differential embeddings of~$S^n$ in $S^{n+q}$ for $q>2$. Ann. Math. (2), 83 (1966), 402--~436.}

\newcommand{\hass}{\bibitem[Ha66C]{Ha66C} \emph{A.~Haefliger,}  Enlacements de spheres en codimension superiure \`a 2, Comment. Math. Helv. 41 (1966-67) 51--72.}

\newcommand{\hase}{\bibitem[Ha68]{Ha68} \emph{A. Haefliger,} Knotted Spheres and Related Geometric Topic,
in Proc. Int. Congr. Math., Moscow, 1966 (Mir, Moscow, 1968), 437--445.}

\newcommand{\hasn}{\bibitem[Ha69]{Ha69} \emph{L.~S.~Harris,} Intersections and embeddings of polyhedra, Topology 8 (1969) 1--26.}

\newcommand{\hasf}{\bibitem[Ha74]{Ha74} * \emph{P. Halmos,} How to talk mathematics. Notices of the Amer. Math. Soc., 21 (1974) 155--158.}

\newcommand{\haef}{\bibitem[Ha84]{Ha84} \emph{N. Habegger,} Obstruction to embedding disks II: a proof of a conjecture by Hudson, Topol. Appl. 17 (1984).}

\newcommand{\haes}{\bibitem[Ha86]{Ha86} \emph{N. Habegger,} Knots and links in codimension greater than 2, Topology, 25:3 (1986) 253--260.}

\newcommand{\hogr}{\bibitem[HG]{HG} * \url{http://www.map.mpim-bonn.mpg.de/Homology_groups_(simplicial;_simple_definition)}}

\newcommand{\hifn}{\bibitem[Hi59]{Hi59} \emph{M. W. Hirsch.} Immersions of manifolds, Trans. Amer. Math. Soc. 93 (1959) 242--276.}

\newcommand{\hjsf}{\bibitem[HJ64]{HJ64} \emph{R. Halin and H. A. Jung.}
Karakterisierung der Komplexe der Ebene und der 2-Sph\"are, Arch. Math. 1964. 15. P.~466-469.}

\newcommand{\hkns}{\bibitem[HK96]{HK96} \emph{N. Habegger and U. Kaiser,} Homotopy classes of 2 disjoint $2p$-spheres in $\R^{3p+1}$, Topol. Appl. 71 (1996) 1--8.}

\newcommand{\hkne}{\bibitem[HK98]{HK98} \emph{N. Habegger and U. Kaiser,} Link homotopy in 2--metastable range, Topology 37:1 (1998) 75--94.}

\newcommand{\hmsnt}{\bibitem[HMS]{HMS93} * \emph{C. Hog-Angeloni, W. Metzler and A. J. Sieradski.}
Two-dimensional homotopy and combinatorial group theory. Cambridge: Cambridge Univ. Press, 1993. (London Math. Soc. Lecture Notes, 197).}

\newcommand{\ho}{\bibitem[Ho]{Ho} * The Hopf fibration, \url{https://www.youtube.com/watch?v=AKotMPGFJYk}}

\newcommand{\hozs}{\bibitem[Ho06]{Ho06} \emph{H. van der Holst,} Graphs and obstructions in four dimensions, J. Combin. Theory Ser. B 96:3 (2006), 388--404.}


\newcommand{\hpzn}{\bibitem[HP09]{HP09} \emph{H. van der Holst and R. Pendavingh,} On a graph property generalizing planarity and flatness, Combinatorica, 29 (2009) 337--361.}

\newcommand{\hssf}{\bibitem[HS64]{HS64} \emph{A. Haefliger and B. Steer,} Symmetry of linking coefficients, Comment. Math. Helv. 39 (1964) 259-270.}

\newcommand{\htsf}{\bibitem[HT74]{HT74} \emph{J. Hopcroft and R. E. Tarjan,} Efficient planarity testing, J. of the Association for Computing Machinery, 21:4 (1974) 549--568.}

\newcommand{\hufn}{\bibitem[Hu59]{hu59} * \emph{S. T. Hu,} Homotopy Theory, Academic Press, New York, 1959.}

\newcommand{\huss}{\bibitem[Hu66]{Hu66} * \emph{J.~F.~P.~Hudson,} Extending piecewise linear isotopies, Proc. London Math. Soc. (3) 16 (1966) 651--668.}

\newcommand{\husn}{\bibitem[Hu69]{Hu69} * \emph{J. F. P. Hudson.} Piecewise linear topology, W. A. Benjamin, Inc., New York-Amsterdam, 1969.}


\newcommand{\io}{\bibitem[Io]{Io} * \url{https://en.wikipedia.org/wiki/Category:Impossible_objects}}

\newcommand{\info}{\bibitem[IF]{IF} * \url{http://www.map.mpim-bonn.mpg.de/Intersection_form}}

\newcommand{\irsf}{\bibitem[Ir65]{Ir65} \emph{M.~C.~Irwin,} Embeddings of polyhedral manifolds, Ann. of Math. (2)
82 (1965) 1--14.}

\newcommand{\isot}{\bibitem[Is]{Is} * \url{http://www.map.mpim-bonn.mpg.de/Isotopy}}


\newcommand{\jqnt}{\bibitem[JQ93]{JQ93} * \emph{A. Jaffe, F. Quinn,} ``Theoretical mathematics'': Toward a cultural synthesis of mathematics and theoretical physics. Bull.Am.Math.Soc. 29 (1993) 1-13. arXiv:math/9307227.}

\newcommand{\jozt}{\bibitem[Jo02]{Jo02} \emph{C. M. Johnson.} An obstruction to embedding a simplicial $n$-complex into a $2n$-manifold, Topology Appl. 122:3 (2002) 581--591.}

\newcommand{\jvz}{\bibitem[JVZ]{JVZ} D. Joji\'c, S. T. Vre\'cica, R. T. \v Zivaljevi\' c,
Topology and combinatorics of 'unavoidable complexes', arXiv:1603.08472v1.}


\newcommand{\kalai}{\bibitem[Ka]{Ka} \emph{G. Kalai,} From Oberwolfach: The Topological Tverberg Conjecture is False, `Combinatorics and more' blog post, February 6, 2015, \url{gilkalai.wordpress.com}}

\newcommand{\kano}{\bibitem[Ka91]{Ka91} \emph{G. Kalai,} The diameter of graphs of convex polytopes and $f$-vector theory, Applied geometry and discrete mathematics, DIMACS Ser. Discrete Math. Theoret. Comput. Sci., vol. 4, Amer. Math. Soc., Providence, RI, 1991, pp. 387--411.}

\newcommand{\kefn}{\bibitem[Ke59]{Ke59} \emph{M. Kervaire.} An interpretation of G. Whitehead's generalization of H. Hopf's invariant, Ann. of Math. 62 (1959), 345--362.}

\newcommand{\kh}{\bibitem[Kh]{Kh} \emph{А.И. Храбров.} Руководство по чтению лекций
\url{http://vm.tstu.tver.ru/topics/pdf_tests/lection.pdf}}

\newcommand{\kho}{\bibitem[Kho]{Kho} \emph{N. Khoroshavkina.} A simple characterization of graphs of cutwidth 2, arXiv:1811.06716.}

\newcommand{\kkrot}{\bibitem[KKR]{KKR} \emph{K. Kawarabayashi, Y. Kobayashi and B. Reed.} The disjoint paths problem in quadratic time, J. of Comb. Theory, Ser. B, 102:2 (2012), 424--435.}

\newcommand{\kln}{\bibitem[KLN]{KLN} \emph{ J. Kratochv\'il, A. Lubiw and J. Ne\v set\v ril.} Noncrossing subgraphs in topological layouts, SIAM J. on Discr. Math. 4(2) (1991), 223--244.}

\newcommand{\kmsth}{\bibitem[KM63]{KM63} \emph{M. A. Kervaire and J. W. Milnor,} Groups of homotopy spheres. I,  Ann. of Math. (2) 77 (1963), 504-537.}

\newcommand{\kozeru}{\bibitem[Ko18]{Ko18} * \emph{Е. Колпаков.}
Доказательство теоремы Радона при помощи понижения размерности, Мат. Просвещение, 23 (2018), arXiv:1903.11055.}

\newcommand{\koze}{\bibitem[Ko18]{Ko18} * \emph{E. Kolpakov.}
A proof of Radon Theorem via lowering of dimension, Mat. Prosveschenie, 23 (2018), arXiv:1903.11055.}

\newcommand{\ko}{\bibitem[Ko]{Ko} \emph{E. Kolpakov.} A `converse' to the Constraint Lemma, arXiv:1903.08910.}

\newcommand{\koon}{\bibitem[Ko19]{Ko19} \emph{E. Kogan.} Linking of three triangles in 3-space, arXiv:1908.03865.}

\newcommand{\koto}{\bibitem[Ko21]{Ko21} \emph{E. Kogan.} On the rank of $\Z_2$-matrices with free entries on the diagonal, arXiv:2104.10668.}

\newcommand{\koee}{\bibitem[Ko88]{Ko88} \emph{U. Koschorke.} Link maps and the geometry of their invariants,
Manuscripta Math. 61:4 (1988) 383--415.}

\newcommand{\kono}{\bibitem[Ko91]{Ko91} \emph{U. Koschorke.} Link homotopy with many components,
Topology 30:2 (1991) 267--281.}

\newcommand{\kons}{\bibitem[Ko97]{Ko97} \emph{U. Koschorke.} A generalization of Milnor's $\mu$-invariants to higher-dimensional link maps, Topology 36:2 (1997) 301--324.}

\newcommand{\kps}{\bibitem[KPS]{KPS} * \emph{A. Kaibkhanov, D. Permyakov and A. Skopenkov.}
Realization of graphs with rotation, \url{http://www.turgor.ru/lktg/2005/3/index.htm}.}

\newcommand{\krzz}{\bibitem[Kr00]{Kr00} \emph{V. S. Krushkal.} Embedding obstructions and 4-dimensional thickenings of 2-complexes, Proc. Amer. Math. Soc. 128:12 (2000) 3683--3691. arXiv:math/0004058. }

\newcommand{\ksnn}{\bibitem[KS99]{KS99} * \emph{П. Кожевников и А. Скопенков.} Узкие деревья на плоскости, Мат. Образование. 1999. 2-3. С.~126-131.}

\newcommand{\kstz}{\bibitem[KS20]{KS20} \emph{R. Karasev and A. Skopenkov.}
Some `converses' to intrinsic linking theorems, Discr. Comp. Geom., 70:3 (2023), 921--930, arXiv:2008.02523.}


\newcommand{\ksto}{\bibitem[KS21]{KS21} * \emph{E. Kogan and A. Skopenkov.} A short proof of the Patak-Tancer theorem on non-embeddability of $k$-complexes in $2k$-manifolds,  arXiv:2106.14010.}

\newcommand{\kstoe}{\bibitem[KS21e]{KS21e} \emph{E. Kogan and A. Skopenkov.}
Embeddings of $k$-complexes in $2k$-manifolds and minimum rank of partial symmetric matrices, arXiv:2112.06636v2.}

\newcommand{\kstf}{\bibitem[KS24]{KS24} \emph{R. Karasev and A. Skopenkov.}
Short proofs of Tverberg-type theorems for cell complexes, Discr. Comp. Geom., (2025), arXiv:2405.05629.}


\newcommand{\kutt}{\bibitem[Ku23]{Ku23} \emph{W. K\"uhnel.} Generalized Heawood Numbers, The Electronic Journal of Combinatorics, 30:4 (2023) \#P4.17.}


\newcommand{\kuse}{\bibitem[Ku68]{Ku68} * \emph{К. Куратовский.} Топология. Т.~1,~2. М.: Мир, 1969.}

\newcommand{\kunfo}{\bibitem[Ku94]{Ku94} \emph{W. K\"uhnel.} Manifolds in the skeletons of convex polytopes, tightness, and generalized Heawood inequalities. In Polytopes: abstract, convex and computational (Scarborough, ON, 1993), volume 440 of NATO Adv. Sci. Inst. Ser. C Math. Phys. Sci., pp. 241--247. Kluwer
Acad. Publ., Dordrecht, 1994.}


\newcommand{\kunf}{\bibitem[Ku95]{Ku95} * \emph{W. K\"uhnel}, Tight Polyhedral Submanifolds and Tight Triangulations, Lecture Notes in Math. 1612, Springer, 1995.}

\newcommand{\kytz}{\bibitem[Ky20]{Ky20} \emph{J. Kyn{\v{c}}l,} Simple Realizability of Complete Abstract Topological Graphs Simplified, Discr. Comp. Geom., 64 (2020) 1--27, arXiv:1608.05867.}


\newcommand{\lazz}{\bibitem[La00]{La00} \emph{F. Lasheras.} An obstruction to 3-dimensional thickening,
Proc. Amer. Math. Soc. 2000. 128. P.~893-902.}

\newcommand{\lfma}{\bibitem[LF]{LF} \url{http://www.map.mpim-bonn.mpg.de/Linking_form}}

\newcommand{\lloe}{\bibitem[LL18]{LL18} \emph{A.S. Levine and T. Lidman.} Simply connected, spineless 4-manifolds, Forum of Math., Sigma, 7 (2019) e14, 1--11, arxiv:1803.01765.}

\newcommand{\lo}{\bibitem[Lo]{Lo} M.~de~Longueville. Notes on the topological Tverberg theorem.
Discrete Math.  247 (2002), no.~1--3, 271--297.
(The paper first appeared in
Discrete Math. 241 (2001) 207--233, but the original version suffered from serious publisher's typesetting errors.)}

\newcommand{\loot}{\bibitem[Lo13]{Lo13} \emph{M. de Longueville.} A course in topological combinatorics. Universitext. Springer, New York (2013).}

\newcommand{\lssn}{\bibitem[LS69]{LS69} \emph{W. B. R. Lickorish and L. C. Siebenmann.}
Regular neighborhoods and the stable range,  Trans. Amer. Math. Soc.. 1969. 139. P.~207-230.}

\newcommand{\lsne}{\bibitem[LS98]{LS98} \emph{L. Lovasz and A. Schrijver,}
A Borsuk theorem for antipodal links and a spectral characterization of linklessly embeddable graphs, Proc. Amer. Math. Soc. 126:5 (1998), 1275-1285.}

\newcommand{\ltof}{\bibitem[LT14]{LT14} \emph{E. Lindenstrauss and M. Tsukamoto,} Mean dimension and an embedding problem: an example, Israel J. Math. 199 (2014).}


\newcommand{\lyzf}{\bibitem[LY04]{LY04} * \emph{Y. Lin and A. Yang,} On 3-cutwidth critical graphs, Discrete Mathematics, 275 (2004), 339--346.}

\newcommand{\lz}{\bibitem[LZ]{LZ} * \emph{S. Lando and A. Zvonkin.} Embedded Graphs. Springer.}


\newcommand{\maez}{\bibitem[Ma80]{Ma80} * R. Mandelbaum, {\em Four-Dimensional Topology: An introduction},
Bull. Amer. Math. Soc. (N.S.) 2 (1980) 1-159.}

\newcommand{\mast}{\bibitem[Ma73]{Ma73} \emph{С. В. Матвеев.} Специальные остовы кусочно-линейных многообразий, Мат. Сборник. 1973. 92. С.~282-293.}

\newcommand{\maste}{\bibitem[Ma73]{Ma73} \emph{S. V. Matveev.} Special skeletons of PL manifolds (in Russian), Mat. Sbornik. 1973. 92. P.~282-293.}

\newcommand{\manz}{\bibitem[Ma90]{Ma90} \emph{W. S.  Massey.} Homotopy classification of 3-component links of codimension greater than 2, Topol.  Appl. 34 (1990) 269--300.}

\newcommand{\mans}{\bibitem[Ma97]{Ma97} \emph{Yu. Makarychev.} A short proof of Kuratowski's graph planarity criterion, J. of Graph Theory, 25 (1997), 129--131.}

\newcommand{\matns}{\bibitem[Mat97]{Mat97} \emph{J. Matou\v sek.} A Helly-type theorem for unions of convex sets. Discr. Comp. Geom., 18:1 (1997) 1-12.}

\newcommand{\mazt}{\bibitem[Ma03]{Ma03} * \emph{J.~Matou{\v{s}}ek.} Using the {B}orsuk-{U}lam theorem:
Lectures on topological methods in combinatorics and geometry. Springer Verlag, 2008.}


\newcommand{\mazf}{\bibitem[Ma05]{Ma05} \emph{V. Manturov.} A proof of the Vasiliev conjecture on the planarity of singular links, Izv. RAN 2005.}

\newcommand{\metn}{\bibitem[Me29]{Me29} \emph{K. Menger.} \"Uber pl\"attbare Dreiergraphen und Potenzen nicht pl\"attbarer Graphen, Ergebnisse Math. Kolloq., 2 (1929) 30--31.}

\newcommand{\mezf}{\bibitem[Me04]{Me04} \emph{S. Melikhov.} Sphere eversions and realization of mappings, Trudy MIAN 247 (2004) 159-181 (in Russian) arXiv:math.GT/0305158.}

\newcommand{\mezs}{\bibitem[Me06]{Me06} \emph{S. A. Melikhov}, The van Kampen obstruction and its relatives, 	
Proc. Steklov Inst. Math 266 (2009), 142-176 (= Trudy MIAN 266 (2009), 149-183), arXiv:math/0612082.}

\newcommand{\meoo}{\bibitem[Me11]{Me11} \emph{S. A. Melikhov}, Combinatorics of embeddings, arXiv:1103.5457.}

\newcommand{\meos}{\bibitem[Me17]{Me17} \emph{S. Melikhov,} Gauss type formulas for link map invariants, arXiv:1711.03530.}

\newcommand{\meoe}{\bibitem[Me18]{Me18} \emph{S. A. Melikhov,} A triple-point Whitney trick, J. Topol. Anal., 2018, 1--6. arXiv:2210.04016.}


\newcommand{\metz}{\bibitem[Me20]{Me20} \emph{S. A. Melikhov,} Topological isotopy and Cochran's derived invariants, in `Topology, Geometry, and Dynamics: Rokhlin Memorial', Contemporary Mathematics, 772, AMS, Providence, RI, 2021. arXiv:2011.01409.}

\newcommand{\mett}{\bibitem[Me22]{Me22} \emph{S. A. Melikhov,} Embeddability of joins and products of polyhedra, Topol. Methods in Nonlinear Analysis, 60:1 (2022), 185-201. arXiv:2210.04015.}

\newcommand{\miff}{\bibitem[Mi54]{Mi54} \emph{J. Milnor,} Link groups, Ann. of Math. 59 (1954), 177--195.}

\newcommand{\miso}{\bibitem[Mi61]{Mi61} \emph{J. Milnor,} A procedure for killing homotopy groups of differentiable manifolds, Proc. Sympos. Pure Math, Vol. III (1961), 39--55.}

\newcommand{\mins}{\bibitem[Mi97]{Mi97} \emph{P. Minc.} Embedding simplicial arcs into the plane, Topol. Proc. 1997. 22. 305--340.}


\newcommand{\adnsvr}{\bibitem[MNS]{MNS} * \emph{А. Мирошников, О. Никитенко и А. Скопенков.} Циклы в графах и в гиперграфах: в направлении теории гомологий, Мат. Просвещение, 35 (2025), 137-184, arXiv:2406.16705.}

\newcommand{\dmnse}{\bibitem[MNS]{MNS} * \emph{A. Miroshnikov, O. Nikitenko, A. Skopenkov.}
Cycles in graphs and in hypergraphs: towards homology theory (in Russian), Mat. Prosveschenie, 35 (2025), 137-184, arXiv:2406.16705.}

\newcommand{\moss}{\bibitem[Mo77]{Mo77} * \emph{E. E. Moise.} Geometric Topology in Dimensions 2 and 3 (GTM), Springer-Verlag, 1977.}

\newcommand{\moen}{\bibitem[Mo89]{Mo89} \textit{B. Mohar}. An obstruction to embedding graphs in
surfaces. Discrete Math. 78 (1989) 135--142.}

\newcommand{\moze}{\bibitem[Mo08]{Mo08} \textit{T. Moriyama}. An invariant of embeddings of 3-manifolds in 6-manifolds and Milnor's triple linking number, J. Math. Sci. Univ. Tokyo, 18 (2011), 193--237. arXiv:0806.3733.}


\newcommand{\mrst}{\bibitem[MRS+]{MRS+} \emph{A. de Mesmay, Y. Rieck, E. Sedgwick, M. Tancer,}
Embeddability in $\R^3$ is NP-hard. arXiv:1708.07734.}

\newcommand{\mesczs}{\bibitem[MS06]{MS06} \emph{S.A. Melikhov, E.V. Shchepin,} The telescope approach to embeddability of compacta. arXiv:math.GT/0612085.}

\newcommand{\msos}{\bibitem[MS17]{MS17}  \emph{T. Maciazek, A. Sawicki.} Homology groups for particles on one-connected graphs
J. Math. Phys. 58, 062103 (2017). arXiv:1606.03414.}

\newcommand{\mstwof}{\bibitem[MST+]{MST+} \emph{J. Matou\v sek, E. Sedgwick, M. Tancer, U. Wagner}, Embeddability in the 3-sphere is decidable, Journal of the ACM 65:1 (2018) 1--49, arXiv:1402.0815.}


\newcommand{\mtzo}{\bibitem[MT01]{MT01} * \emph{B. Mohar and C. Thomassen.} Graphs on Surfaces.
The John Hopkins University Press, 2001.}

\newcommand{\mtwoz}{\bibitem[MTW10]{MTW10} \emph{J. Matou\v sek, M. Tancer, U. Wagner.} A geometric proof of
the colored Tverberg theorem, Discr. and Comp. Geometry, 47:2 (2012), 245--265. arXiv:1008.5275.}


\newcommand{\mtwoo}{\bibitem[MTW]{MTW} \emph{J. Matou\v sek, M. Tancer, U. Wagner.}
Hardness of embedding simplicial complexes in $\R^d$, J. Eur. Math. Soc. 13:2 (2011), 259--295. arXiv:0807.0336.}


\newcommand{\musf}{\bibitem[Mu74]{Mu74} \emph{J. Muncres,} Elementary differential topology, Complement to the book by J. W. Milnor and J. D. Stasheff, {\it Characteristic Classes}, Ann. of Math. St. 76 (1974), Princeton Univ. Press, Princeton, NJ.}

\newcommand{\mwoe}{\bibitem[MW18]{MW18} * \emph{F. Manin, S. Weinberger.} Algorithmic aspects of immersibility and embeddability, Intern. Math. Res. Notices, rnae170. arXiv:1812.09413.}


\newcommand{\mwsn}{\bibitem[MW69]{MW69} * \emph{J. MacWilliams}. Orthogonal matrices over finite fields. Amer. Math. Monthly, 76 (1969) 152--164.}

\newcommand{\mwofo}{\bibitem[MW14]{MW14} \emph{I. Mabillard and U. Wagner.} Eliminating Tverberg Points, I. An Analogue of the Whitney Trick, Proc. of the 30th Annual Symp. on Comp. Geom. (SoCG'14), ACM, New York, 2014, pp. 171--180.}

\newcommand{\mwof}{\bibitem[MW15]{MW15} \emph{I. Mabillard and U. Wagner.}
Eliminating Higher-Multiplicity Intersections, I. A Whitney Trick for Tverberg-Type Problems. arXiv:1508.02349.}


\newcommand{\mwos}{\bibitem[MW16]{MW16} \emph{I. Mabillard and U. Wagner.} Eliminating Higher-Multiplicity Intersections, II. The Deleted Product Criterion in the $r$-Metastable Range. arXiv:1601.00876v2.}

\newcommand{\mwosd}{\bibitem[MW16']{MW16'} \emph{I. Mabillard and U. Wagner.} Eliminating Higher-Multiplicity Intersections, II. The Deleted Product Criterion in the r-Metastable Range,
Proceedings of the 32nd Annual Symposium on Computational Geometry (SoCG'16).}


\newcommand{\natz}{\bibitem[Na20]{Na20} * \emph{R. Naimi,} A brief survey on intrinsically knotted and linked graphs, arXiv:2006.07342.}

\newcommand{\neno}{\bibitem[Ne91]{Ne91} \emph{S. Negami.} Ramsey theorems for knots, links and spatial graphs,
Trans. Amer. Math. Soc., 324 (1991), 527--541.}



\newcommand{\nizz}{\bibitem[Ni00]{Ni00} \emph{R. Nikkuni.} The second skew-symmetric cohomology group and spatial embeddings of graphs, J. Knot Theory Ram. 9 (2000), 387--411.}

\newcommand{\nkon}{\bibitem[NKS]{NKS} * \emph{L. T. Nguyen, J. Kim, B. Shim.}
Low-Rank Matrix Completion: A Contemporary Survey. arXiv:1907.11705.}

\newcommand{\noss}{\bibitem[No76]{No76} * \emph{С. П. Новиков.} Топология-1. М.: Наука, 1976. (Итоги науки и техники. ВИНИТИ. Современные проблемы математики. Основные направления, 12).}

\newcommand{\nszn}{\bibitem[NS09]{NS09} \emph{I. Novik and E. Swartz,} Socles of Buchsbaum modules, complexes and posets, Adv. Math. 222 (2009), 2059-2084. arXiv:0711.0783.}

\newcommand{\nwns}{\bibitem[NW97]{NW97} \emph{A. Nabutovsky, S. Weinberger}. Algorithmic aspects of homeomorphism problems. arXiv:math/9707232.}


\newcommand{\omoe}{\bibitem[Om18]{Om18} * \emph{А. Омельченко,} Теория графов. М.: МЦНМО, 2018.}

\newcommand{\orszo}{\bibitem[ORS]{ORS} \emph{A. Onischenko, D. Repov\v s and A. Skopenkov.}
Resolutions of 2-polyhedra by fake surfaces and embeddings into $\R^4$, Contemp. Math.  288 (2001) 396--400.}

\newcommand{\ossf}{\bibitem[OS74]{OS74} \emph{R. P. Osborne and R. S. Stevens.} Group presentations
corresponding to spines of 3-manifolds, I, Amer. J.~Math. 1974. 96. P.~454-471; II, Amer. J.~Math. 1977. 234.
P.~213-243; III, Amer. J.~Math. 1977. 234 P.~245-251.}


\newcommand{\oz}{\bibitem[Oz]{Oz} \emph{M. \"Ozaydin,} Equivariant maps for the symmetric group, unpublished,
\url{http://minds.wisconsin.edu/handle/1793/63829}.}

\newcommand{\panof}{\bibitem[Pan15]{Pan15} \emph{K. Panagiotis.} A note on the topology of irreducible $SO(3)$-manifolds, 	arXiv:1508.06150.}

\newcommand{\paof}{\bibitem[Pa15]{Pa15} \emph{S. Parsa,} On links of vertices in simplicial $d$-complexes embeddable in the Euclidean $2d$-space, Discrete Comput. Geom. 59:3 (2018), 663--679.
This is arXiv:1512.05164v4 up to numbering of sections, theorems etc.; we refer to numbering in arxiv version.
Correction: Discrete Comput. Geom. 64:3 (2020) 227--228.}

\newcommand{\paoe}{\bibitem[Pa18]{Pa18} \emph{S. Parsa,} On links of vertices in simplicial $d$-complexes
embeddable in the euclidean $2d$-space, arXiv:1512.05164v6.}

\newcommand{\patz}{\bibitem[Pa20]{Pa20} \emph{S. Parsa,} On links of vertices in simplicial $d$-complexes
embeddable in the euclidean $2d$-space, arXiv:1512.05164v8.}


\newcommand{\patzl}{\bibitem[Pa20]{Pa20} \emph{S. Parsa,}
Correction to: On the Links of Vertices in Simplicial $d$-Complexes Embeddable in the Euclidean $2d$-Space,
Discrete Comput. Geom. 64:3 (2020) 227--228.}

\newcommand{\patza}{\bibitem[Pa20]{Pa20} \emph{S. Parsa,} On the Smith classes, the van Kampen obstruction and embeddability of $[3]*K$, arXiv:2001.06478.}

\newcommand{\patzb}{\bibitem[Pa20b]{Pa20b} \emph{S. Parsa,} On the embeddability of $[3]*K$, arXiv:2001.06506.}

\newcommand{\pato}{\bibitem[Pa21]{Pa21} \emph{S. Parsa,} Instability of the Smith index under joins and applications to embeddability, Trans. Amer. Math. Soc. 375 (2022), 7149--7185, arXiv:2103.02563.}

\newcommand{\pak}{\bibitem[Pa]{Pa} * \emph{I. Pak}, Lectures on Discrete and Polyhedral Geometry, \url{http://www.math.ucla.edu/~pak/geompol8.pdf}.}

\newcommand{\peze}{\bibitem[Pe08]{Pe08} \emph{Д. Пермяков.} Классификация погружений графов в плоскость,
Вестник МГУ, сер.1, 2008, N5, 55-56.}

\newcommand{\peos}{\bibitem[Pe16]{Pe16} \emph{Д. Пермяков.} Матем. сб., 207:6 (2016),  93--112.}

\newcommand{\pest}{\bibitem[Pe72]{Pe72} * \emph{B. B. Peterson.} The Geometry of Radon's Theorem, Amer. Math. Monthly 79 (1972), 949-963.}


\newcommand{\prnf}{\bibitem[Pr95]{Pr95} * \emph{V. V. Prasolov.} Intuitive topology. Amer. Math. Soc., Providence, R.I., 1995.}

\newcommand{\prnfr}{\bibitem[Pr95]{Pr95} * \emph{В. В. Прасолов.} Наглядная топология. М.: МЦНМО, 1995.}


\newcommand{\przs}{\bibitem[Pr06]{Pr06} * \emph{V. V. Prasolov.}
Elements of Combinatorial and Differential Topology, 2006, GSM 74, Amer. Math. Soc., Providence, RI.}

\newcommand{\przsru}{\bibitem[Pr04]{Pr04} * \emph{В. В. Прасолов.}
Элементы комбинаторной и дифференциальной топологии. М.: МЦНМО, 2004. \url{http://www.mccme.ru/prasolov}.}

\newcommand{\przse}{\bibitem[Pr07]{Pr07} * \emph{V. V. Prasolov.} Elements of homology theory. 2007, GSM 74, Amer. Math. Soc., Providence, RI.}


\newcommand{\przseru}{\bibitem[Pr06]{Pr06} * \emph{В. В. Прасолов.} Элементы теории гомологий. М.: МЦНМО, 2006.}


\newcommand{\psns}{\bibitem[PS96]{PS96} * \emph{V. V. Prasolov, A. B. Sossinsky } Knots, Links, Braids, and 3-manifolds. Amer. Math. Soc. Publ., Providence, R.I., 1996.}


\newcommand{\pszf}{\bibitem[PS05]{PS05} * \emph{В. В. Прасолов и М. Б. Скопенков.}
Рамсеевская теория зацеплений, Мат. Просвещение. 2005. 9. С.~108--115.}

\newcommand{\pszfen}{\bibitem[PS05]{PS05} * \emph{V. V. Prasolov and M.B. Skopenkov.}
Ramsey link theory, Mat, Prosvescheniye, 9 (2005), 108--115.}

\newcommand{\psoo}{\bibitem[PS11]{PS11} \emph{Y. Ponty and C. Saule.} A combinatorial framework for designing (pseudoknotted) RNA algorithms, Proc. of the 11th Intern. Workshop on Algorithms in Bioinformatics, WABI'11, 250--269.}


\newcommand{\pstz}{\bibitem[PS20]{PS20} \emph{S. Parsa and A. Skopenkov.} On embeddability of joins and their `factors', Topol. Appl., 326 (2023) 108409, arXiv:2003.12285.}



\newcommand{\psszn}{\bibitem[PSS]{PSS} \emph{M. J. Pelsmajer, M. Schaefer and D. Stasi.} Strong Hanani-Tutte on the projective plane. SIAM J. Discrete Math., 23:3 (2009) 1317--1323.}

\newcommand{\psszs}{\bibitem[PSS]{PSS} \emph{M. J. Pelsmajer, M. Schaefer, and D. \v Stefankovi\v c.}
Removing even crossings. J. Combin. Theory Ser. B, 97(4):489–500, 2007.}

\newcommand{\pton}{\bibitem[PT19]{PT19} \emph{P. Pat\'ak and M. Tancer.} Embeddings of $k$-complexes into $2k$-manifolds. Discrete Comput. Geom. 71 (2024), 960--991. arXiv:1904.02404v4.}

\newcommand{\pw}{\bibitem[PW]{PW} \emph{I. Pak, S. Wilson}, G\lowercase{EOMETRIC REALIZATIONS OF POLYHEDRAL COMPLEXES}, \linebreak \url{http://www.math.ucla.edu/~pak/papers/Fary-full31.pdf}.}


\newcommand{\razf}{\bibitem[RA05]{RA05} * \emph{J. L. Ram\'irez Alfons\'in.} Knots and links in spatial graphs: a survey. Discrete Math., 302 (2005), 225--242.}

\newcommand{\rep}{\bibitem[Rep]{Rep} Referee's report on the paper ``Some `converses' to intrinsic linking theorems', \url{https://www.mccme.ru/circles/oim/materials/ksreport.pdf}}

\newcommand{\rnoo}{\bibitem[RN11]{RN11} * \emph{R. L. Ricca, B. Nipoti.} Gauss' linking number revisited.
J. of Knot Theory and Its Ramif. 20:10 (2011) 1325--1343. \url{https://www.maths.ed.ac.uk/~v1ranick/papers/ricca.pdf} .}

\newcommand{\rrstz}{\bibitem[RRS]{RRS} * \emph{V. Retinskiy, A. Ryabichev and A. Skopenkov.}
Motivated exposition of the proof of the Tverberg Theorem (in Russian).
Mat. Prosveschenie, 27 (2021), 166--169. arXiv:2008.08361.}


\newcommand{\rssec}{\bibitem[RS68]{RS68} \emph{C. P. Rourke and B. J. Sanderson,} Block bundles II, Ann. of Math. (2), 87 (1968) 431--483.}

\newcommand{\rsst}{\bibitem[RS72]{RS72} * \emph{C. P. Rourke and B. J. Sanderson,}
\newblock Introduction to Piecewise-Linear Topology,
\newblock \emph{Ergebn.\ der Math.} 69, Springer-Verlag, Berlin, 1972.}

\newcommand{\rsstr}{\bibitem[RS72]{RS72} * \emph{К. П. Рурк и Б. Дж. Сандерсон.} Введение в кусочно-линейную топологию, Москва. Мир. 1974.}

\newcommand{\rsns}{\bibitem[RS96]{RS96} * \emph{D. Repov\v s and A. B. Skopenkov.}
Embeddability and isotopy of polyhedra in Euclidean spaces,
Proc. of the Steklov Inst. Math. 1996. 212. P.~173-188.}

\newcommand{\rsne}{\bibitem[RS98]{RS98} \emph{D. Repov\v s and A. B. Skopenkov.}
A deleted product criterion for approximability of a map by embeddings, Topol. Appl. 1998. 87 P.~1-19.}

\newcommand{\rsnn}{\bibitem[RS99]{RS99} * \emph{D. Repov\v s and A. B. Skopenkov.} New results on embeddings of polyhedra and manifolds into Euclidean spaces,
Russ. Math. Surv. 54:6 (1999), 1149--1196.}


\newcommand{\rsnnd}{\bibitem[RS99']{RS99'} * \emph{Д. Реповш и А. Скопенков.}
Кольца Борромео и препятствия к вложимости, Труды МИРАН. 1999. 225. С.~331-338.}

\newcommand{\rszz}{\bibitem[RS00]{RS00} \emph{D. Repov\v s and A. Skopenkov.} Cell-like resolutions of polyhedra by special ones,  Colloq. Math. 2000. 86:2. P. 231--237.}

\newcommand{\rszzd}{\bibitem[RS00']{RS00'} * \emph{Д. Реповш и А. Скопенков.} Характеристические классы для начинающих, Мат. Просвещение. 2000. 4. С.~151-176.}

\newcommand{\rszo}{\bibitem[RS01]{RS01} \emph{D. Repovs and A. Skopenkov.} On contractible $n$-dimensional compacta, non-embeddable into $\R^{2n}$, Proc. Amer. Math. Soc. 129 (2001) 627--628.}

\newcommand{\rszt}{\bibitem[RS02]{RS02} * \emph{Д. Реповш и А. Скопенков.} Теория препятствий для начинающих,
Мат. Просвещение. 2002. 6. C.~60-77.}

\newcommand{\rszf}{\bibitem[RS04]{RS04} \emph{N. Robertson and P. Seymour.} Graph Minors. XX. Wagner's conjecture, J. of Comb. Theory, B, 92:2 (2004) 325--357.}

\newcommand{\rssnf}{\bibitem[RSS]{RSS95} \emph{D. Repov\v s, A. B. Skopenkov  and E. V. \v S\v cepin.}
On uncountable collections of continua and their span, Colloq. Math. 1995. 69:2. P.~289-296.}

\newcommand{\rssnfd}{\bibitem[RSS']{RSS95'} \emph{D. Repov\v s, A. B. Skopenkov and E. V \v S\v cepin.}
On embeddability of $X\times I$ into Euclidean space, Houston J.~Math. 1995. 21. P.~199-204.}

\newcommand{\rssz}{\bibitem[RSS+]{RSSZ} * \emph{A. Rukhovich, A. Skopenkov, M. Skopenkov, A. Zimin},
Realizability of hypergraphs, \url{https://www.turgor.ru/lktg/2013/1/1-1en.pdf} .}


\newcommand{\rstnt}{\bibitem[RST']{RST93} \emph{N. Robertson, P. Seymour and R. Thomas}, Linkless embeddings of graphs in 3-space, Bull. of the Amer. Math. Soc., 21 (1993) 84--89.}

\newcommand{\rstno}{\bibitem[RST]{RST91} * \emph{N. Robertson, P. Seymour and R. Thomas}, A survey of
linkless embeddings, Graph Structure Theory (Seattle, WA, 1991), Contemp. Math. 147, (1993) 125--136.}


\newcommand{\rwzl}{\bibitem[RWZ+]{RWZ+} \emph{Y. Ren, C. Wen, S. Zhen, N. Lei, F. Luo, D.X. Gu},
Characteristic class of isotopy for surfaces, J. Syst. Sci. Complex. 33 (2020) 2139--2156.}


\newcommand{\saeo}{\bibitem[Sa81]{Sa81} \emph{H. Sachs.} On spatial representation of finite graphs,
in: Finite and infinite sets (Eger, 1981), 649--662, Colloq. Math. Soc. Janos Bolyai, 37, North-Holland, Amsterdam, 1984.}

\newcommand{\sano}{\bibitem[Sa91]{Sa91} \emph{K. S. Sarkaria.}
A one-dimensional Whitney trick and Kuratowski's graph planarity criterion, Israel J.~Math. 73 (1991), 79--89.}


\newcommand{\sanov}{\bibitem[Sa91g]{Sa91g} \emph{K. S. Sarkaria.} A generalized Van Kampen-Flores theorem, Proc. Amer. Math. Soc. 111 (1991), 559--565.}

\newcommand{\sant}{\bibitem[Sa92]{Sa92} \emph{K. S. Sarkaria.} Tverberg’s theorem via number fields. Israel J. Math., 79:317–320, 1992.}

\newcommand{\sann}{\bibitem[Sa99]{Sa99} O. Saeki {\em On punctured 3-manifolds in 5-sphere}, Hiroshima Math. J. 29 (1999) 255--272.}


\newcommand{\sazz}{\bibitem[Sa00]{Sa00} \emph{K. S. Sarkaria.} Tverberg partitions and Borsuk-Ulam theorems. Pacific J. Math., 196:1 (2000) 231--241.}

\newcommand{\sczf}{\bibitem[Sc04]{Sc04} \emph{T. Sch\"oneborn.} On the Topological Tverberg Theorem, arXiv:math/0405393.}


\newcommand{\scot}{\bibitem[Sc13]{Sc13} * \emph{M. Schaefer.} Hanani-Tutte and related results. In Geometry --- intuitive, discrete, and convex, Bolyai Soc. Math. Stud., 24 (2013), 259--299.
\url{http://ovid.cs.depaul.edu/documents/htsurvey.pdf} }


\newcommand{\sctz}{\bibitem[Sc20]{Sc20} \emph{M. Schaefer.} The Graph Crossing Number and
its Variants: A Survey. The Electr. J. of Comb. (2020), DS21, \url{https://www.combinatorics.org/files/Surveys/ds21/ds21v5-2020.pdf}}


\newcommand{\scef}{\bibitem[Sc84]{Sc84} \emph{E.~V.~\v S\v cepin.} Soft mappings of manifolds, Russian Math. Surveys, 39:5 (1984).}

\newcommand{\shfs}{\bibitem[Sh57]{Sh57} \emph{A. Shapiro,} Obstructions to the embedding of a complex in a Euclidean space, I, The first obstruction, Ann. Math. 66 (1957), 256--269.}


\newcommand{\shen}{\bibitem[Sh89]{Sh89} * \emph{Ю. А. Шашкин,} Неподвижные точки, М., Наука, 1989.}

\newcommand{\shoe}{\bibitem[Sh18]{Sh18} * \emph{S. Shlosman},  Topological Tverberg Theorem: the proofs and the counterexamples, Russian Math. Surveys, 73:2 (2018), 175–182. arXiv:1804.03120.}

\newcommand{\sisn}{\bibitem[Si69]{Si69} \emph{K. Sieklucki.} Realization of mappings, Fund. Math. 1969. 65. P.~325-343.}

\newcommand{\sios}{\bibitem[Si16]{Si16} \emph{S. Simon,} Average-Value Tverberg Partitions via Finite Fourier Analysis, Israel J. Math., 216 (2016), 891-904, arXiv:1501.04612.}
 

\newcommand{\sknf}{\bibitem[Sk94]{Sk94} \emph{А. Скопенков.} Геометрическое доказательство теоремы
Нойвирта об утолщаемости 2-мерных полиэдров, Math. Notes. 1995. 58:5. P.~1244-1247.}


\newcommand{\skns}{\bibitem[Sk97]{Sk97} \emph{A. Skopenkov,} On the deleted product criterion for embeddability of manifolds in $\R^m$, Comment. Math. Helv. 72 (1997), 543--555.}

\newcommand{\skne}{\bibitem[Sk98]{Sk98} \emph{A. B. Skopenkov.} On the deleted product criterion for embeddability in $\R^m$, Proc. Amer. Math. Soc., 126:8 (1998), 2467-2476.}

\newcommand{\skzz}{\bibitem[Sk00]{Sk00} \emph{A. Skopenkov,} On the generalized Massey--Rolfsen invariant for link maps, Fund. Math. 165 (2000), 1--15.}

\newcommand{\skzt}{\bibitem[Sk02]{Sk02} \emph{A. Skopenkov,} On the Haefliger-Hirsch-Wu invariants for embeddings and immersions, Comment. Math. Helv. 77 (2002), 78--124.}

\newcommand{\skzth}{\bibitem[Sk03]{Sk03} \emph{M. Skopenkov,} Embedding products of graphs into Euclidean spaces,
Fund. Math. 179 (2003),~191--198, arXiv:0808.1199.}

\newcommand{\skzthd}{\bibitem[Sk03']{Sk03'} \emph{M. Skopenkov,} On approximability by embeddings of cycles in the plane, Topol. Appl. 134 (2003),~1--22, arXiv:0808.1187.}

\newcommand{\skzf}{\bibitem[Sk05]{Sk05} * \emph{A. Skopenkov,}
On the Kuratowski graph planarity criterion, Mat. Prosveschenie, 9 (2005), 116-128. arXiv:0802.3820.}


\newcommand{\skzs}{\bibitem[Sk05i]{Sk05i} \emph{A. Skopenkov,} A new invariant and parametric connected sum of embeddings, Fund. Math. 197 (2007) 253--269. arxiv:math/0509621.}

\newcommand{\skzei}{\bibitem[Sk05]{Sk05} \emph{A.  Skopenkov,} A classification of smooth embeddings of
4-manifolds in 7-space, I, Topol. Appl., 157 (2010) 2094--2110. arXiv:math/0512594.}

\newcommand{\skze}{\bibitem[Sk06]{Sk06} * \emph{A. Skopenkov,} Embedding and knotting of manifolds in Euclidean spaces, London Math. Soc. Lect. Notes, 347 (2008) 248--342. arXiv:math/0604045.}

\newcommand{\skzsi}{\bibitem[Sk06']{Sk06'} \emph{A. Skopenkov,} A classification of smooth embeddings of 3-manifolds in 6-space, Math. Zeitschrift, 260:3 (2008) 647--672. arxiv:math/0603429.}

\newcommand{\skzsc}{\bibitem[Sk06c]{Sk06c} \emph{A. Skopenkov,} Classification of embeddings below the metastable dimension, arXiv:math/0607422.}

\newcommand{\skozp}{\bibitem[Sk08]{Sk08} \emph{A.  Skopenkov,} Embeddings of $k$-connected $n$-manifolds into
$\R^{2n-k-1}$. arxiv:math/0812.0263; earlier version published in Proc. Amer. Math. Soc., 138 (2010) 3377--3389.}

\newcommand{\skoz}{\bibitem[Sk10]{Sk10} * \emph{А. Скопенков,} Вложения в плоскость графов с вершинами степени 4,
Мат. просвещение, 21 (2017), arXiv:1008.4940.}

\newcommand{\skoo}{\bibitem[Sk11]{Sk11} \emph{M. Skopenkov,} When is the set of embeddings finite up to isotopy? Intern. J. Math. 26:7 (2015), 28 pp. arXiv:1106.1878.}

\newcommand{\skofo}{\bibitem[Sk14]{Sk14} \emph{A. Skopenkov,} How do autodiffeomorphisms act on embeddings, Proc. A of the Royal Society of Edinburgh, 148:4 (2018), 835--848. arXiv:1402.1853.}

\newcommand{\sks}{\bibitem[Sk14]{Sk14} * \emph{A. Skopenkov,} Realizability of hypergraphs and intrinsic linking  theory, Mat. Prosveschenie, 32 (2024), 125--159, arXiv:1402.0658.}

\newcommand{\sksr}{\bibitem[Sk14]{Sk14} * \emph{А. Скопенков,} Реализуемость гиперграфов и неотъемлемая зацепленность, Мат. просвещение, 32 (2024), 125--159. arXiv:1402.0658.}


\newcommand{\skof}{\bibitem[Sk15]{Sk15} * \emph{А. Скопенков,} Алгебраическая топология с геометрической точки зрения, Москва, МЦНМО, 2015 (1е издание).}

\newcommand{\skofe}{\bibitem[Sk15]{Sk15} * \emph{A. Skopenkov,} Algebraic Topology From Geometric Viewpoint (in Russian), MCCME, Moscow, 2015 (1st edition). }

\newcommand{\skofel}{\bibitem[Sk15e]{Sk15e} * \emph{А. Скопенков,} Алгебраическая топология
с геометрической точки зрения, эл. версия, \url{http://www.mccme.ru/circles/oim/home/combtop13.htm\#photo}}


\newcommand{\skoomp}{\bibitem[Sk11]{Sk11} A. Skopenkov, A simple proof of the Abel-Ruffini theorem (in Russian),
Mat. Prosveschenie, 15 (2011) 113-126, arXiv:1102.2100.}

\newcommand{\skofmp}{\bibitem[Sk15]{Sk15} A. Skopenkov, A short elementary proof of the Ruffini-Abel Theorem (in Russian),
Mat. Prosveschenie, 36 (2026) 95--113. Abridged English version is published in [Sk21m, \S8]; full English version: arXiv:1508.03317.}

\newcommand{\skotzr}{\bibitem[Sk20]{Sk20} * \emph{А. Скопенков,} Алгебраическая топология с геометрической точки зрения, Москва, МЦНМО, 2020 (2е издание).
Обновляемая версия части книги: \url{http://www.mccme.ru/circles/oim/obstruct.pdf}}

\newcommand{\skotz}{\bibitem[Sk20]{Sk20} * \emph{A. Skopenkov,} Algebraic Topology From Geometric Standpoint (in Russian), MCCME, Moscow, 2020 (2nd edition).
Update of a part: \url{http://www.mccme.ru/circles/oim/obstruct.pdf} .
Part of the English translation: \url{https://www.mccme.ru/circles/oim/obstructeng.pdf}.}


\newcommand{\skofp}{\bibitem[Sk15]{Sk15} \emph{A. Skopenkov,} Classification of knotted tori,
Proc. A of the Royal Soc. of Edinburgh, 150:2 (2020), 549-567. Full version: arXiv:1502.04470.}


\newcommand{\skos}{\bibitem[Sk16]{Sk16} * \emph{A. Skopenkov,} A user's guide to the topological Tverberg Conjecture, arXiv:1605.05141v5. Abridged earlier published version: Russian Math. Surveys, 73:2 (2018), 323--353.}


\newcommand{\skosd}{\bibitem[Sk16']{Sk16'} * \emph{A. Skopenkov,} Stability of intersections of graphs in the plane and the van Kampen obstruction, Topol. Appl. 240(2018) 259--269, arXiv:1609.03727.}


\newcommand{\skosc}{\bibitem[Sk16c]{Sk16c} * \emph{A. Skopenkov,}  Embeddings in Euclidean space: an introduction to their classification, to appear in Boll. Man. Atl. 
\url{http://www.map.mpim-bonn.mpg.de/} \url{Embeddings_in_Euclidean_space:_an_introduction_to_their_classification}}

\newcommand{\skosie}{\bibitem[Sk16e]{Sk16e} * \emph{A. Skopenkov,} Embeddings just below the stable range: classification, to appear in Boll. Man. Atl.
\url{http://www.map.mpim-bonn.mpg.de/Embeddings_just_below_the_stable_range:_classification}}

\newcommand{\skost}{\bibitem[Sk16t]{Sk16t} * \emph{A. Skopenkov,} 3-manifolds in 6-space, to appear in Boll. Man. Atl. \url{http://www.map.mpim-bonn.mpg.de/3-manifolds_in_6-space}.}

\newcommand{\skosf}{\bibitem[Sk16f]{Sk16f} * \emph{A. Skopenkov,} 4-manifolds in 7-space, to appear in Boll. Man. Atl. \url{http://www.map.mpim-bonn.mpg.de/4-manifolds_in_7-space}.}

\newcommand{\skosh}{\bibitem[Sk16h]{Sk16h} * \emph{A. Skopenkov,} High codimension links, to appear in Boll. Man. Atl.
\linebreak
\url{http://www.map.mpim-bonn.mpg.de/High_codimension_links}.}

\newcommand{\skosi}{\bibitem[Sk16i]{Sk16i} * \emph{A. Skopenkov,} Isotopy, submitted to Boll. Man. Atl.
\url{http://www.map.mpim-bonn.mpg.de/Isotopy}.}

\newcommand{\skosk}{\bibitem[Sk16k]{Sk16k} * \emph{A. Skopenkov,} Knotted tori,
\url{http://www.map.mpim-bonn.mpg.de/Knotted_tori}.}

\newcommand{\skoss}{\bibitem[Sk16s]{Sk16s} * \emph{A. Skopenkov,} Knots, i.e. embeddings of spheres,
\linebreak
\url{http://www.map.mpim-bonn.mpg.de/Knots,_i.e._embeddings_of_spheres}.}

\newcommand{\skose}{\bibitem[Sk17]{Sk17} \emph{A. Skopenkov,}
Eliminating higher-multiplicity intersections in the metastable dimension range. arXiv:1704.00143.}

\newcommand{\skosed}{\bibitem[Sk17v]{Sk17v} * \emph{A. Skopenkov,}
On van Kampen-Flores, Conway-Gordon-Sachs and Radon theorems,  arXiv:1704.00300.}

\newcommand{\sk}{\bibitem[Sk17o]{Sk17o} \emph{A. Skopenkov,} On the metastable Mabillard-Wagner conjecture.  arXiv:1702.04259.}

\newcommand{\skmos}{\bibitem[Sk17d]{Sk17d} \emph{M. Skopenkov}. Discrete field theory: symmetries and conservation laws, arXiv:1709.04788.}

\newcommand{\skoe}{\bibitem[Sk18]{Sk18} * \emph{A. Skopenkov.} Invariants of graph drawings in the plane.
Arnold Math. J., 6 (2020) 21--55; full version: arXiv:1805.10237.}


\newcommand{\skoer}{\bibitem[Sk18]{Sk18} * \emph{А. Скопенков,} Инварианты изображений графов на плоскости,
Мат. просвещение, 31 (2023), 74-127. arXiv:1805.10237.}


\newcommand{\sktthd}{\bibitem[Sk23']{Sk23'} * \emph{A. Skopenkov.} Invariants of graph drawings in the plane (in Russian). Mat. Prosveschenie, 31 (2023), 74-127. arXiv:1805.10237.}

\newcommand{\skoeo}{\bibitem[Sk18o]{Sk18o} * \emph{A. Skopenkov.} A short exposition of S. Parsa's theorems on intrinsic linking and non-realizability. Discr. Comp. Geom. 65:2 (2021), 584--585; full version:  arXiv:1808.08363.}


\newcommand{\skona}{\bibitem[Sk19]{Sk19} * \emph{A. Skopenkov,} A short exposition of the Levine-Lidman example of spineless 4-manifolds, arXiv:1911.07330.}

\newcommand{\sktze}{\bibitem[Sk21m]{Sk21m} * \emph{A. Skopenkov.} Mathematics via Problems. Part 1: Algebra. Amer. Math. Soc., Providence, 2021. Preliminary version: \url{https://www.mccme.ru/circles/oim/algebra_eng.pdf}}

\newcommand{\sktz}{\bibitem[Sk20u]{Sk20u} * \emph{A. Skopenkov.} A user's guide to basic knot and link theory,
in: Topology, Geometry, and Dynamics, Contemporary Mathematics, vol. 772, Amer. Math. Soc., Providence, RI, 2021, pp. 281--309.
Russian version: Mat. Prosveschenie 27 (2021), 128--165. arXiv:2001.01472.}

\newcommand{\sktzru}{\bibitem[Sk20u]{Sk20u} * \emph{А. Скопенков.} Основы теории узлов и зацеплений для пользователя, Мат. просвещение, 27 (2021), 128--165. arXiv:2001.01472.}

\newcommand{\sktzo}{\bibitem[Sk20o]{Sk20o} \emph{A. Skopenkov.} On some results of S. Abramyan and T. Panov, arXiv:2005.11152.}

\newcommand{\sktzr}{\bibitem[Sk20e]{Sk20e} * \emph{A. Skopenkov.}
Extendability of simplicial maps is undecidable, Discr. Comp. Geom., 69:1 (2023), 250--259, arXiv:2008.00492.}


\newcommand{\sktzd}{\bibitem[Sk21d]{Sk21d} * \emph{A. Skopenkov.}
On different reliability standards in current mathematical research, arXiv:2101.03745.
More often updated version: \url{https://www.mccme.ru/circles/oim/rese_inte.pdf}.}

\newcommand{\sktt}{\bibitem[Sk22]{Sk22} * \emph{A. Skopenkov.} Invariants of embeddings of 2-surfaces in 3-space,
arXiv:2201.10944.}

\newcommand{\skttn}{\bibitem[Sk22n]{Sk22n} * \emph{A. Skopenkov}, Netflix problem and realization of (hyper)graphs, \url{https://www.mccme.ru/circles/oim/home/netflix20sep.pdf}}

\newcommand{\sktth}{\bibitem[Sk23]{Sk23} \emph{A. Skopenkov.}  To S. Parsa's theorem on embeddability of joins, arXiv:2302.11537.}

\newcommand{\sktf}{\bibitem[Sk24]{Sk24} * \emph{A. Skopenkov.} Double and triple linking numbers in space (in Russian). Mat. Prosveschenie, 33 (2024), 87--132.}

\newcommand{\sktfr}{\bibitem[Sk24]{Sk24} * \emph{А. Скопенков.} Двойные и тройные коэффициенты зацепления в пространстве. Мат. просвещение, 33 (2024), 87--132.}

\newcommand{\sktfb}{\bibitem[Sk24]{Sk24} \emph{A. Skopenkov.} The band connected sum and the second Kirby move for higher-dimensional links, Stud. Sci. Math. Hung., 62:4 (2025) 320--335, arXiv:2406.15367.}


\newcommand{\sktfe}{\bibitem[Sk24]{Sk24} \emph{A. Skopenkov.}
Embeddings of $k$-complexes in $2k$-manifolds and minimum rank of partial symmetric matrices, arXiv:2112.06636v5.}

\newcommand{\skd}{\bibitem[Sk]{Sk} * \emph{А. Скопенков.} Алгебраическая топология с алгоритмической точки зрения, 
\url{http://www.mccme.ru/circles/oim/algor.pdf}.}

\newcommand{\skde}{\bibitem[Sk]{Sk} * \emph{A. Skopenkov.} Algebraic Topology From Algorithmic Standpoint, draft of a book, mostly in Russian,
\url{http://www.mccme.ru/circles/oim/algor.pdf}.}


\newcommand{\skon}{\bibitem[Skw]{Skw} * \emph{A. Skopenkov.} Whitney trick for eliminating multiple intersections, slides for talks at St Petersburg, Brno, Kiev, Moscow,  \url{https://www.mccme.ru/circles/oim/eliminat_talk.pdf}.}

\newcommand{\skl}{\bibitem[EEF]{EEF} * {\it Proposed by D. Eliseev, A. Enne, M. Fedorov, A. Glebov, N. Khoroshavkina, E. Morozov, A. Skopenkov, R. \v Zivaljevi\'c.}
A user's guide to knot and link theory, \url{https://www.turgor.ru/lktg/2019/3} .}

\newcommand{\skr}{\bibitem[Skr]{Skr} * \emph{A. Skopenkov.} Realizability of hypergraphs, slides for talks,  \url{https://www.mccme.ru/circles/oim/algor1_beamer.pdf}.}


\newcommand{\skt}{\bibitem[Skt]{Skt} * \emph{A. Skopenkov.} Transparent anonymous peer review,
\linebreak
\url{https://www.mccme.ru/circles/oim/home/transp_peer_review.htm} .}

\newcommand{\rslktg}{\bibitem[KRR+]{RRSl} * Towards higher-dimensional combinatorial geometry, presented by
E. Kogan, V. Retinskiy, E. Riabov and A. Skopenkov, \url{https://www.mccme.ru/circles/oim/multicomb.pdf} .}



\newcommand{\sm}{\bibitem[Sm]{Sm} S. Smirnov.}

\newcommand{\sper}{\bibitem[Sp]{Sp} * Sperner's lemma defeats the rental harmony problem, \url{https://www.youtube.com/watch?v=7s-YM-kcKME}.}

\newcommand{\sset}{\bibitem[SS83]{SS83} \emph{Е. В. Щепин, М. А. Штанько.} Спектральный критерий вложимости компактов в евклидовы пространства, Труды Ленинградской Международной Топологической конференции. Л.: Наука, 1983. С.~135-142.}

\newcommand{\ssnt}{\bibitem[SS92]{SS92} \emph{J.~Segal and S.~Spie\.z.} Quasi embeddings and embeddings of polyhedra in $\R^m$,  Topol. Appl., 45 (1992) 275--282.}

\newcommand{\sszt}{\bibitem[SS03]{SS03} \emph{F. W. Simmons and F. E. Su.}
Consensus-halving via theorems of Borsuk-Ulam and Tucker, Math. Social Sciences 45 (2003) 15–25. \url{https://www.math.hmc.edu/~su/papers.dir/tucker.pdf}.}

\newcommand{\ssot}{\bibitem[SS13]{SS13} \emph{M. Schaefer and D. \v Stefankovi\v c.} Block additivity of $\Z_2$-embeddings. In Graph drawing, volume 8242 of Lecture Notes in Comput. Sci., 185--195.
Springer, Cham, 2013. \url{http://ovid.cs.depaul.edu/documents/genus.pdf}}

\newcommand{\sstt}{\bibitem[SS23]{SS23} \emph{A. Skopenkov and O. Styrt,} Embeddability of joinpowers, and minimal rank of partial matrices, arXiv:2305.06339.}

\newcommand{\sssne}{\bibitem[SSS]{SSS} \emph{J. Segal, A. Skopenkov and S. Spie\. z.}
Embeddings of polyhedra in $\R^m$ and the deleted product obstruction, Topol. Appl., 85 (1998), 225-234.}

\newcommand{\sstnf}{\bibitem[SST95]{SST95} \emph{R. S. Simon, S. Spie\. z and H. Toru\'nczyk.}
T\lowercase{HE EXISTENCE OF EQUILIBRIA IN CERTAIN GAMES, SEPARATION FOR FAMILIES OF CONVEX FUNCTIONS
AND A THEOREM OF BORSUK-ULAM TYPE}, Israel J. Math 92 (1995) 1--21.}

\newcommand{\sstzt}{\bibitem[SST02]{SST02} \emph{R. S. Simon, S. Spie\. z and H. Toru\'nczyk.}
E\lowercase{QUILIBRIUM EXISTENCE AND TOPOLOGY IN SOME REPEATED GAMES WITH INCOMPLETE INFORMATION},
Trans. Amer. Math. Soc. 354:12 (2002) 5005-5026.}

\newcommand{\stez}{\bibitem[ST80]{ST80} * {\it H.~Seifert and W.~Threlfall.}
A textbook of topology, v~89 of {\em Pure and Applied Mathematics}.
Academic Press, New York-London, 1980.}


\newcommand{\stzs}{\bibitem[ST07]{ST07} * \emph{А. Скопенков и А. Телишев.}
И вновь о критерии Куратовского планарности графов, Мат. Просвещение, 11 (2007), 159--160.}

\newcommand{\stzse}{\bibitem[ST07]{ST07} * \emph{A. Skopenkov and A. Telishev}, Once again on the Kuratowski graph planarity criterion, Mat. Prosveschenie, 11 (2007), 159-160. arXiv:0802.3820.}

\newcommand{\stos}{\bibitem[ST17]{ST17} \emph{A. Skopenkov  and M. Tancer,}
Hardness of almost embedding simplicial complexes in $\R^d$, Discr. Comp. Geom., 61:2 (2019), 452--463. arXiv:1703.06305.}

\newcommand{\stno}{\bibitem[ST91]{ST91} \emph{S.~Spie\. z and H.~Toru\'nczyk}, Moving compacta in $\R^m$ apart,
Topol. Appl. 41 (1991), 193--204.}

\newcommand{\sttt}{\bibitem[St24]{St24} \emph{M. Starkov,} An example of an `unlinked' set of $2k+3$ points in $2k$-space, arXiv:2402.09002.}

\newcommand{\sunt}{\bibitem[Su]{Su} * \emph{Д. Судзуки.} Основы дзэн-буддизма. Наука дзэн --- ум дзэн. Киев: Преса Украiни. 1992.}

\newcommand{\stwh}{\bibitem[SW]{SW} * \url{http://www.map.mpim-bonn.mpg.de/Stiefel-Whitney_characteristic_classes}}

\newcommand{\sz}{\bibitem[SZ05]{SZ} \emph{T. Sch\"oneborn and G. Ziegler}, The Topological Tverberg Theorem and Winding Numbers, J. Comb. Theory, Ser. A, 112:1 (2005) 82--104, arXiv:math/0409081.}

\newcommand{\szno}{\bibitem[Sz91]{Sz91} \emph{A.~Sz\"ucs,} On the cobordism groups of immersions and embeddings,
Math. Proc. Camb. Phil. Soc., 109 (1991) 343--349.}


\newcommand{\ta}{\bibitem[Ta]{Ta} * Handbook of Graph Drawing and Visualization. ed. by R. Tamassia, CRC Press, 2016.}


\newcommand{\tanfo}{\bibitem[Ta94]{Ta94} \emph{K. Taniyama,} Cobordism, homotopy and homology of graphs in $\R^3$,
Topology 33:3 (1994), 509--523.}

\newcommand{\tanf}{\bibitem[Ta95]{Ta95} \emph{K. Taniyama,} Homology classification of spatial embeddings of a graph, Topol. Appl. 65 (1995) 205--228.}

\newcommand{\tazz}{\bibitem[Ta00]{Ta00} \emph{K. Taniyama,} Higher dimensional links in a simplicial complex embedded in a sphere, Pacific Jour. of Math. 194:2 (2000), 465-467.}

\newcommand{\theo}{\bibitem[Th81]{Th81} * \emph{C.~Thomassen,} Kuratowski's theorem, J.~Graph. Theory 5 (1981), 225--242.}

\newcommand{\tooo}{\bibitem[To11]{To11} \emph{Tonkonog D.} Embedding 3-manifolds with boundary into closed 3-manifolds, Topol. Appl. 158 (2011), 1157-1162. arXiv:1003.3029.}


\newcommand{\tsbzf}{\bibitem[TSB]{TSB} \emph{D. M. Thilikos, M. Serna and H. L. Bodlaender},
Cutwidth I: A linear time fixed parameter algorithm, J. of Algorithms, 56:1 (2005), 1--24.}


\newcommand{\tsbzfd}{\bibitem[TSB05']{TSB05'} \emph{D. M. Thilikos, M. Serna and H. L. Bodlaender},
Cutwidth II: , J. of Algorithms, 56:1 (2005), 25--49.}



\newcommand{\umse}{\bibitem[Um78]{Um78} \emph{B. Ummel.} The product of nonplanar complexes does not imbed in 4-space, Trans. Amer. Math. Soc., 242 (1978) 319--328.}




\newcommand{\vant}{\bibitem[Va92]{Va92} * \emph{V.~A.~Vassiliev.} Complements of discriminants of smooth maps: Topology and applications, Amer. Math. Soc., Providence, RI, 1992 (рус. перевод: В. А. Васильев, Топология дополнений к дискриминантам, Фазис, Москва, 1997).}

\newcommand{\val}{\bibitem[Val]{Val} * \url{https://en.wikipedia.org/wiki/Valknut}}


\newcommand{\vi}{\bibitem[Vi]{Vi} * \emph{O. Viro.}
Some integral calculus based on Euler characteristic, Lect. Notes in Math. 1346.}

\newcommand{\vizt}{\bibitem[Vi02]{Vi02} * \emph{Э. Б. Винберг.} Курс алгебры. Москва. Факториал Пресс. 2002.}

\newcommand{\vizteng}{\bibitem[Vi02]{Vi02} * \emph{E. B. Vinberg.} A Course in Algebra. Graduate Studies in Mathematics, vol. 56. 2003.}

\newcommand{\vinhzs}{\bibitem[VINH07]{VINH07} * \emph{О. Я. Виро, О. А. Иванов, Н. Ю. Нецветаев и В. М. Харламов.}
Элементарная топология, МЦНМО. 2007.}

\newcommand{\vktt}{\bibitem[vK32]{vK32} \emph{E.~R.~van~Kampen}, Komplexe in euklidischen R\"aumen, Abh. Math. Sem. Hamburg, 9 (1933) 72--78; Berichtigung dazu, 152--153.}

\newcommand{\kafo}{\bibitem[vK41]{vK41} \emph{E. R. van Kampen,} Remark on the address of S. S. Cairns,
in Lectures in Topology, 311--313, University of Michigan Press, Ann Arbor, MI, 1941.}

\newcommand{\vo}{\bibitem[Vo96]{vo96} \emph{A. Yu. Volovikov,} On a topological generalization of the Tverberg theorem. Math. Notes 59:3 (1996), 324--326.}

\newcommand{\vopns}{\bibitem[Vo96v]{Vo96v} \emph{A. Yu. Volovikov,} On the van Kampen-Flores Theorem.
Math. Notes 59:5 (1996), 477--481.}


\newcommand{\vznt}{\bibitem[VZ93]{VZ93} \emph{A. Vu\v ci\'c and R. T. \v Zivaljevi\'c}, Note on a conjecture of Sierksma, Discr. Comput. Geom. 9 (1993), 339-349.}

\newcommand{\vzzn}{\bibitem[VZ09]{VZ09} \emph{S. T. Vre\'cica and R. T. \v Zivaljevi\'c},  Chessboard complexes
indomitable, J. of Comb. Theory, Ser. A 118:7 (2011), 2157--2166. arXiv:0911.3512.}


\newcommand{\walst}{\bibitem[Wa62]{Wa62} \emph{C.~T.~C.~Wall}, Classification of $(n-1)$-connected $2n$-manifolds, Ann. of Math., 75 (1962) 163--189.}


\newcommand{\wallss}{\bibitem[Wa67]{Wa67} \emph{C.~T.~C.~Wall.} Classification problems in differential topology, IV, Thickenings, Topology 1966. 5. P. 73--94.}

\newcommand{\waldss}{\bibitem[Wa67m]{Wa67m} \emph{F. Waldhausen.} Eine Klasse von 3-dimensional Mannigfaltigkeiten, I. Invent. Math. 1967. 3. P.~308-333.}

\newcommand{\walsz}{\bibitem[Wa70]{Wa70} \emph{C. T. C. Wall,} Surgery on compact manifolds,
1970, Academic Press, London.}

\newcommand{\wess}{\bibitem[We67]{We67} \emph{C.~Weber.} Plongements de poly\`edres dans le domain metastable, Comment. Math. Helv. 42 (1967), 1--27.}

\newcommand{\whit}{\bibitem[Wl]{Wl} * \url{https://en.wikipedia.org/wiki/Whitehead_link}}

\newcommand{\winum}{\bibitem[Wn]{Wn} * \url{https://en.wikipedia.org/wiki/Winding_number}}

\newcommand{\wrss}{\bibitem[Wr77]{Wr77} \emph{P. Wright.} Covering 2-dimensional polyhedra by 3-manifolds spines.
Topology. 16 (1977), 435--439.}

\newcommand{\wufe}{\bibitem[Wu58]{Wu58} \emph{W. T. Wu.} On the realization of complexes in a euclidean space (in Chinese): I, Sci Sinica, 7 (1958) 251--297; II, Sci Sinica, 7 (1958) 365--387; III, Sci Sinica, 8 (1959) 133--150.}

\newcommand{\wufn}{\bibitem[Wu59]{Wu59} \emph{W.~T.~Wu.} On the isotopy of a finite complex in Euclidean space, I, II, Science Record, N.S. 3:8 (1959) 342--347, 348--351.}

\newcommand{\wusf}{\bibitem[Wu65]{Wu65} * \emph{W. T. Wu.} A Theory of Embedding, Immersion and Isotopy of Polytopes in an Euclidean Space. Peking: Science Press, 1965.}


\newcommand{\yann}{\bibitem[Ya99]{Ya99} \emph{Z. Yang.} Computing Equilibria and Fixed Points: The Solution of Nonlinear Inequalities, Kluwer, Springer Science + Business Media, 1990.}


\newcommand{\zesz}{\bibitem[Ze60]{Ze60} \emph{E. C. Zeeman}, Unknotting spheres in five dimensions, Bull. Amer. Math. Soc. 66 (1960) 198.
\linebreak
\url{https://www.ams.org/journals/bull/1960-66-03/S0002-9904-1960-10431-4/S0002-9904-1960-10431-4.pdf}}

\newcommand{\z}{\bibitem[Ze]{Z} * \emph{E. C. Zeeman}, A Brief History of Topology, UC Berkeley, October 27, 1993, On the occasion of Moe Hirsch's 60th birthday, \url{http://zakuski.utsa.edu/~gokhman/ecz/hirsch60.pdf}.}

\newcommand{\zioz}{\bibitem[Zi10]{Zi10} * \emph{D. \v Zivaljevi\'c}, Borromean and Brunnian Rings,
\url{http://www.rade-zivaljevic.appspot.com/borromean.html}.}

\newcommand{\zioo}{\bibitem[Zi11]{Zi11} * \emph{G. M. Ziegler}, 3N Colored Points in a Plane, Notices of the Amer. Math. Soc., 58:4 (2011), 550-557.}


\newcommand{\zot}{\bibitem[Zi13]{Z13} \emph{A. Zimin.} Alternative proofs of the Conway-Gordon-Sachs Theorems, arXiv:1311.2882.}

\newcommand{\zstf}{\bibitem[ZS25]{ZS25} А.А. Заславский и А.Б. Скопенков, Исследовательские задачи и Московская математическая конференция школьников, arXiv:2512.22191.}

\newcommand{\zss}{\bibitem[ZSS]{ZSS} * Элементы математики в задачах: через олимпиады и кружки к профессии.
Сборник под редакцией А. Заславского, А. Скопенкова и М. Скопенкова. М.: МЦНМО, 2018.
Обновляемая версия части книги: \url{http://www.mccme.ru/circles/oim/materials/sturm.pdf}.}

\newcommand{\zsse}{\bibitem[ZSS]{ZSS} Elements of mathematics via problems: from olympiades and math circles to a profession (in Russian), editors A. Zaslavsky, A. Skopenkov, and M. Skopenkov. MCCME, Moscow, 2018,
updated part of the book: \url{http://www.mccme.ru/circles/oim/sturm.pdf}.}


\newcommand{\zu}{\bibitem[Zu]{Zu} \emph{J. Zung.} A non-general-position Parity Lemma,
\url{http://www.turgor.ru/lktg/2013/1/parity.pdf}.}



\adnsv

\bibitem[AS]{AS} \emph{E. Alkin, A. Skopenkov,} A generalization of the Conway--Gordon--Sachs theorem to multiple linking, draft. 

\bibitem[AT21]{AT21} \emph{G. Avramidi, T. Tam Nguyen Phan,} Fungible obstructions to embedding 2-complexes, arXiv:2105.10984. 

\bibitem[FD05]{FD05} \emph{T. Fleming, A. Diesl,} Intrinsically linked graphs and even linking number, 	Algebr. Geom. Topol., 5 (2005), 1419--1432. 

\bibitem[FFN]{FFN} \emph{E. Flapan, W. Fletcher, R. Nikkuni,} Reduced Wu and generalized Simon invariants for spatial graphs, Math. Proc. Camb. Philos. Soc., 156:3 (2014), 521--544.



\hcon

\kstz

\dmnse

\bibitem[Mi]{Mi} \emph{A. Mizev,} Realization of sets of integers by linking numbers of embeddings of $K_6$ in 3-space, draft. 

\bibitem[MM01]{MM01} \emph{B. Mellor and P. Melvin.} A geometric interpretation of Milnor's triple linking numbers, Algebr. Geom. Topol., 3 (2003), 557--568, arXiv:math/0110001. 

\natz
\neno 
\nizz
\prnf
\psns
\pszfen
\rnoo

\bibitem[Si86]{Si86} \emph{J. Simon.} Topological chirality of certain molecules, Topology, 25 (1986), 229--235.

\skde

\bibitem[Sk']{Sk'} * \emph{A. Skopenkov.} 
Knotting of graphs in 3-space, and configuration spaces,
\url{https://old.mccme.ru//circles//oim/knogramodkno.pdf}.

\skos
\skosc
\skoe
\sktz

\sktze
\sktf

\bibitem[ST03]{ST03} \emph{R. Shinjo, K. Taniyama,} Homology classification of spatial graphs by linking numbers and Simon invariants, Topol. Appl. 134 (2003), 53--67.

\tanfo
\tanf
\zioz


\end{thebibliography}
\end{document}



\newpage
\title{Инварианты зацепленности пространственных графов}

\author{Э. Алкин, Ю. Хромин, А. Скопенков}

\thanks{\emph{Э. Алкин, А. Скопенков:} Московский физико-технический институт.
\newline
\emph{Ю. Хромин:} Физико-математическая школа №146 (Пермь).
\newline
\emph{А. Скопенков:} Независимый московский университет, \url{https://users.mccme.ru/skopenko}.
\newline
Этот текст представлен в качестве проекта Летней Конференции Турнира Городов 2026 года, см. \url{https://www.turgor.ru/lktg/index.php} .
Мы благодарны А. Мирошникову, А. Мизеву, О. Никитенко и Р. Никкуни за полезные обсуждения,
и О. Никитенко за подготовку некоторых рисунков.
Некоторые рисунки подготовлены с использованием статьи \cite{Ta95}.}
 
\date{}
 
\begin{abstract}
Мы излагаем определения коэффициентов зацепления и чисел Ву--Саймона для пространственных графов.  
Мы излагаем некоторые результаты и предлагаем нерешенные задачи. 
Изложение основано на плоских диаграммах, поэтому знания стереометрии не требуется. 
\end{abstract}

\maketitle

\tableofcontents

Основные определения даны в \S\ref{s:lkpd}.  
Основные результаты сформулированы в \S\S \ref{s:main}, \ref{s:maini} с использованием понятия \emph{коэффициента зацепления}, определенного (при помощи проекций) в \S\ref{s:lkn}. 
Наброски доказательств основных результатов даны в \S\S \ref{s:real}, \ref{s:sol}, \cite[\S9]{AKS}. 

\smallskip 
\textbf{Замечание} (о стиле этого текста). 

Основные идеи представлены на <<олимпиадных>> примерах: на простейших
частных случаях, свободных от технических деталей, и со сведением к необходимому минимуму научного языка.
Благодаря этому и текст становится доступным для начинающих, и удается быстро добраться до интересных сложных и важных результатов, методов и теоретических идей.

Таким образом, многие из приведенных идей работают для более общих случаев.
Мы не тратим время читателя на несложные обобщения.
Их легко придумать или найти в литературе.
Трудно
именно \emph{применить} общую теорию для ярких результатов, сформулированных вне этой теории, когда неизвестно заранее, какую теорию нужно применять (и вообще, можно ли хоть какую-то теорию применить).

\algor{

\emph{`Listeners are prepared to accept unstated (but hinted) generalizations much more than they are able, on the spur of the moment, to decode a precisely stated abstraction and to re-invent the special cases that motivated it in the first place.'}\footnote{Слушатели в гораздо большей степени готовы воспринимать не сформулированные (но упомянутые) обобщения, чем способны <<навскидку>> раскодировать явно сформулированную абстракцию и воспроизвести те конкретные примеры, которые изначально ее мотивировали.} \cite{Ha74}

\emph{`The preceding theorem generalizes to many situations. In fact, there may not actually be an optimal generalization in the sense that no matter what generalization is given, someone could produce a more general one. One of the advantages of being an understander rather than a theorem quoter is that one may be able to obtain approaches to a wide variety of theorems some of which may not even have been formulated yet.'}\footnote{Предшествующая теорема имеет многочисленные обобщения. В действительности, «наилучшего» обобщения может и не быть, в том смысле, что для каждого обобщения кто-то сможет найти еще более общее. Одним из преимуществ читателя, ориентированного на понимание сути, а не просто на цитирование теорем, является способность получать подходы к широкому спектру теорем, некоторые из которых могут быть даже еще не сформулированы.} \cite{Bi83}

}

Более простой материал приводится, чтобы сделать естественным и доступным более сложный.
Попытка начинать с более простого (например, с частных случаев) повышает самостоятельность --- а, значит, глубину и надежность --- освоения материала.
Проще самому доказать частный случай, самому продумать переход от частного к общему, чем  самому сразу доказать общий случай.
Самостоятельно придуманное надежнее запоминается и легче модифицируется.
Кроме того, обычно на частном случае проще отловить и исправить ошибки.
Подробнее см. \cite[\S11.2]{ZSS}.


Как правило, мы приводим \emph{формулировку} красивого или важного утверждения \emph{перед} последовательностью определений и результатов, составляющих его \emph{доказательство}.\algor{\footnote{Часто происходит обратное: формулировки красивых результатов и важных проблем, ради которых была придумана теория, приводятся только {\it после} продолжительного изучения этой теории (или не приводятся совсем). Это развивает представление о математике как науке, изучающей немотивированные понятия и теории. Такое представление принижает ценность математики.}}
В~таких случаях для доказательства утверждения требуется часть дальнейшего материала.
Об этом указано после формулировки утверждения.
Некоторые результаты приведены без доказательств, тогда даются ссылки вместо указаний. 
Если к утверждению (или задаче) не приведено ни доказательство, ни ссылка на него, то оно несложно.
Основные результаты называются <<теоремами>>, менее важные результаты "--- <<утверждениями>>, важные вспомогательные результаты "--- <<леммами>>.
\lktg{Если не оговорено противное, то все это (а также примеры) --- задачи для решения.}

В тексте есть небольшое количество задач (полезно рассматривать и утверждения как задачи).
\algor{Изучение путем решения задач характерно для серьезного обучения математике, см. \cite[\S1.1]{HC19}, \cite[\S1.2]{ZSS} и данные там ссылки.
Оно продолжает древнюю культурную традицию.\footnote{Например, послушники дзенских монастырей обучаются, размышляя над загадками, данными им наставниками.
Впрочем, эти загадки являются скорее наводящими на размышления парадоксами, а не задачами.
См. подробнее \cite{Su}.}
Задачи пронумерованы жирными цифрами, слово <<задача>> опускается.
In a lecture we often start with a question and only then give an answer, e.g.
`Может ли количество ... быть равным 7?'  instead of Assertion \ref{pla-kon}.a.}
Если условие задачи является формулировкой утверждения, то в~задаче требуется это утверждение доказать (и тогда в ссылках мы называем это утверждение утверждением, а не задачей).
\algor{
Краткие указания
к задачам
приводятся сразу, более длинные в конце параграфа.}
\emph{Загадкой} называется не сформулированный четко вопрос; здесь нужно придумать и~четкую формулировку, и~доказательство.
\algor{Если номер задачи помечен звездочкой,
то эта задача посложнее соседних, или менее важна для дальнейшего.}

Определения важных понятий даны \textbf{жирным шрифтом}, чтобы их было проще найти.
\algor{{\it Разделы и задачи, отмеченные звездочкой, а также замечания,} не используются в дальнейшем.}

\lktg{Мы приглашаем школьников представлять свои письменные решения на конференциях школьников, см., например, \cite{ZS25}.}

\section{Зацепления и плоские диаграммы}\label{s:lkpd}

\begin{figure}[H]\centering

\begin{minipage}{0.24\textwidth}
\centering
\includegraphics[scale=1.2]{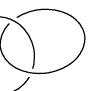}
\par (a)
\end{minipage}
\hfill
\begin{minipage}{0.24\textwidth}
\centering
\includegraphics[scale=1.2]{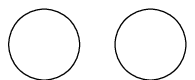}
\par (b)
\end{minipage}
\hfill
\begin{minipage}{0.24\textwidth}
\centering
\includegraphics[scale=0.6]{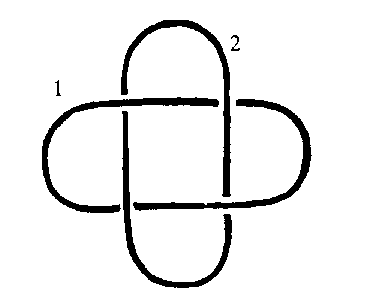}
\par (c)
\end{minipage}
\hfill
\begin{minipage}{0.24\textwidth}
\centering
\includegraphics[scale=0.4]{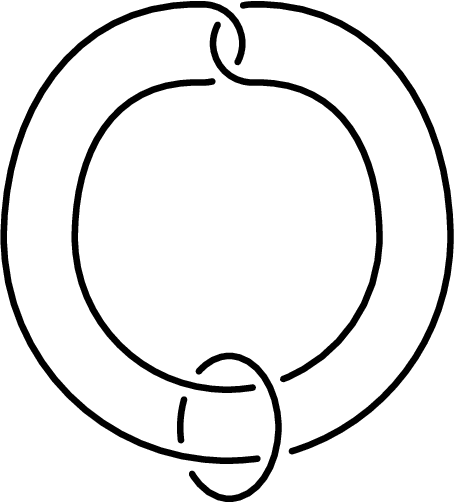}
\par (d)
\end{minipage}

\vspace{0.2cm}
\caption{Зацепление Хопфа (a), тривиальное зацепление (b) и еще два зацепления (c, d)}
\label{f:hopf}
\end{figure}

Можно представить себе {\it узел} как тонкую упругую нить, концы которой склеены; см. компоненты на рисунке~\ref{f:hopf}. 
Строго говоря, {\bf узел} — это пространственная замкнутая несамопересекающаяся ломаная.
{\bf Зацепление} (рисунок~\ref{f:hopf}) — это упорядоченный набор попарно непересекающихся узлов, которые называются {\it компонентами} зацепления.

Зацепления обычно изображаются с помощью их <<удачных>> проекций на плоскость, называемых {\it плоскими диаграммами} (рисунок~\ref{f:hopf};
под проекцией мы понимаем тень).
Представьте, что вы кладете веревку на стол и аккуратно записываете, как она пересекает саму себя (т.е. какая часть лежит сверху, а какая снизу).
{\it Плоская диаграмма зацепления} — это его проекция \emph{общего положения} на плоскость,  
вместе с указанием в каждом перекрестке, какой отрезок <<проходит снизу>>, а какой <<проходит сверху>> (рисунок~\ref{f:hopf}). 
 Проекции одного и того же зацепления на разные плоскости могут выглядеть совершенно по-разному.

В данном тексте мы строго работаем не с зацеплениями и их плоскими диаграммами, а с \emph{абстрактными плоскими диаграммами}. 
Здесь <<абстрактная>> означает <<не обязательно построенная на основе какого-либо зацепления>>, и в дальнейшем часто опускается.

Строго говоря, набор замкнутых ломаных на плоскости называется набором {\it общего положения}, если

(i) их вершины попарно различны,
 
(ii) любые два отрезка ломаных либо не пересекаются, либо пересекаются в общей концевой точке, либо пересекаются в общей внутренней точке, и
 
(iii) никакие три внутренности этих отрезков не имеют общей точки.
 
Общие внутренние точки, такие как указано в (ii), называются {\bf перекрестками}.

{\bf Плоская диаграмма (абстрактная)} состоит из 

$\bullet$ упорядоченного набора общего положения некоторых замкнутых ломаных на плоскости, 
  
$\bullet$ указания в каждом перекрестке, какой отрезок <<проходит снизу>>, а какой <<проходит сверху>> (рисунок~\ref{f:hopf}). 

Каждая из замкнутых ломаных называется {\bf компонентой} плоской диаграммы.

Отрезок ломаной, проходящий снизу, называется {\bf проходом}, а проходящий сверху называется {\bf переходом}.

\begin{assertion}\label{p:diag} *
Для любой абстрактной однокомпонентной плоской диаграммы существует узел, проецирующийся в эту диаграмму.
(Такой узел не обязательно единственный; см., однако, \cite[Утверждение 1.1.c]{Sk20u}.)
\end{assertion}

\section{Коэффициент зацепления по модулю 2}\label{s:lkn2}

Понятие \emph{коэффициента зацепления} двух замкнутых пространственных кривых было введено К. Гауссом в его исследованиях по электромагнетизму \cite{RN11}. 
Неформально говоря, коэффициент зацепления — это число оборотов, которое одна кривая делает вокруг другой. 
Не совсем ясно, как определить это число (даже неформально) для сложных пар кривых. 
Мы дадим строгое определение на языке  плоских диаграмм. 
  
{\bf Коэффициент зацепления по модулю 2}, $\lk_2\in\Z_2$, двухкомпонентной плоской диаграммы — это взятое по модулю $2$ число перекрестков, в которых первая компонента проходит сверху.

\begin{figure}[H]\centering
\begin{minipage}{0.2\textwidth}
\centering
\includegraphics[scale=1.2]{3-7.borr-link.eps}
\par (b)
\end{minipage}
\quad
\begin{minipage}{0.2\textwidth}
\centering
\includegraphics[scale=1.2]{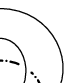}
\par (w)
\end{minipage}
\quad
\begin{minipage}{0.2\textwidth}
\centering
\includegraphics[scale=1.2]{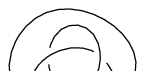}
\par (t)
\end{minipage}
\vspace{0.1cm}
\caption{Кольца Борромео (b), зацепление Уайтхеда (w) и трилистник (t).
\newline
Никакие два из колец Борромео не зацеплены, но все три вместе зацеплены (т.е. их нельзя продеформировать в три попарно непересекающихся шара, ср. \S\ref{s:isot})}
\label{f:borwhitre}
\end{figure}

\begin{pr}\label{a:lk2} 
(a,b,c) Для выбранного вами порядка компонент найдите коэффициенты зацепления по модулю 2 плоских диаграмм на рисунке~\ref{f:hopf}, двухкомпонентных поддиаграмм колец Борромео на рисунке~\ref{f:borwhitre}.b, и каждой пары непересекающихся циклов длины 3 на рисунке~\ref{f:ex-graph-emb}, справа.   
\end{pr}

\begin{figure}[H]\centering
\includegraphics[scale=0.6]{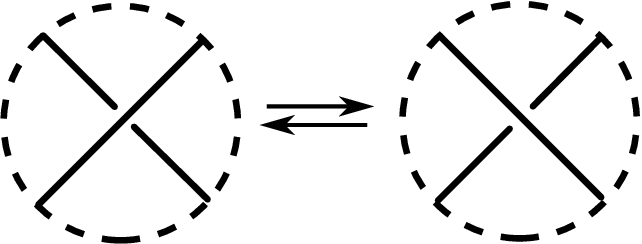}
\vspace{0.1cm}
\caption{Изменение перекрестка.
\emph{Изображены фрагменты плоских диаграмм в виде дисков, ограниченных пунктирными окружностями.
Вне дисков диаграммы идентичны.
Никакие другие ребра плоских диаграмм, кроме изображенных, не пересекают диски.
То же самое для рисунков~\ref{f:cross-ori},~\ref{f:triplic},~\ref{f:reid}.}}
\label{f:cross}
\end{figure}

{\bf Изменение перекрестка} — это замена прохода на переход или наоборот, см. рисунок~\ref{f:cross}.

\begin{assertion}\label{con-lk2} Для любых двухкомпонентных плоских диаграмм $D_+$ и $D_-$, отличающихся изменением перекрестка, 
образованного пересечением разных компонент, имеем \linebreak $\lk\phantom{}_2D_+-\lk\phantom{}_2D_- = 1\in\Z_2$. 
\end{assertion}

\begin{assertion}\label{a:lk2-pr}
Перестановка компонент сохраняет коэффициент зацепления по модулю 2.
\end{assertion}

В доказательстве (и далее) вы можете  использовать без доказательства следующую {\it лемму о четности:} 
любые две замкнутые ломаные общего положения на плоскости 
пересекаются в четном числе точек.
Обсуждение и доказательство см. в \S1.3 <<Число пересечений для ломаных на плоскости>> \cite{Sk18, Sk}. 

\textit{Обратная} к ломаной $l = L_1 L_2 \ldots L_m$ — это ломаная $\overline l := L_m L_{m-1} \ldots L_1$.
\textit{Конкатенация} ломаных $l = L_1\ldots L_m$ и $p = L_mP_2\ldots P_n$ — это ломаная $lp := L_1 \ldots L_m P_2 \ldots P_n$.  

\begin{lemma}\label{l:add2} Пусть $l_1,l_2,l_3$ — ломаные, соединяющие точку $A$ с точкой $B$, и не пересекающие ломаную $p$.
Тогда
$$\lk\phantom{}_2(l_1\overline{l_2}, p) + \lk\phantom{}_2(l_2\overline{l_3}, p) = \lk\phantom{}_2(l_1\overline{l_3}, p).$$
\end{lemma}

\section{Коэффициент зацепления}\label{s:lkn}

\begin{figure}[H]\centering
\includegraphics[scale=0.5]{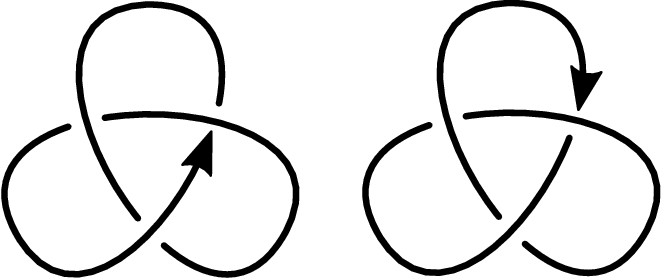}
\caption{Два трилистника с противоположными ориентациями}
\label{f:knot-ori}
\end{figure}

Вы знаете, что такое ориентированная ломаная, значит, вы знаете, что такое \emph{ориентированное зацепление} (рисунок~\ref{f:knot-ori}) и что такое \emph{ориентированная 
плоская диаграмма}.

\begin{figure}[H]\centering
\includegraphics[scale=0.8]{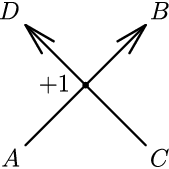}
\qquad 
\includegraphics[scale=0.8]{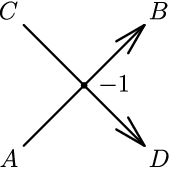}
\caption{Знак пары пересекающихся отрезков}
\label{f:sign}
\end{figure}

Пусть $(\overrightarrow{AB},\overrightarrow{CD})$ — упорядоченная пара пересекающихся векторов (ориентированных отрезков) на плоскости.
Определим {\bf знак} пары как $+1$, если $ABC$ ориентирован по часовой стрелке, и $-1$ в противном случае (рисунок~\ref{f:sign}).

{\bf Коэффициент зацепления} $\lk$ двухкомпонентной ориентированной плоской диаграммы 
— это сумма знаков во всех тех перекрестках, где первая компонента проходит над второй.
Здесь в каждом перекрестке {\it первый} ({\it второй}) вектор — это ориентированное ребро первой (второй) компоненты.

\begin{pr}\label{e:lk} (a,b,c) Для выбранных вами порядка и ориентации компонент найдите коэффициенты зацепления плоских диаграмм на рисунке~\ref{f:hopf}, двухкомпонентных поддиаграмм колец Борромео на рисунке~\ref{f:borwhitre}.b, и каждой пары циклов длины 3 на рисунке~\ref{f:ex-graph-emb}, справа.  

(d) Для любого целого $n$ существует ориентированная двухкомпонентная плоская диаграмма, коэффициент зацепления которой равен $n$. 

(e) Для любой тройки целых чисел $k,l,m$ существует ориентированная трехкомпонентная плоская диаграмма, попарные коэффициенты зацепления которой равны $k,l,m$. 
\end{pr}

\begin{figure}[H]\centering
\includegraphics[scale=0.4]{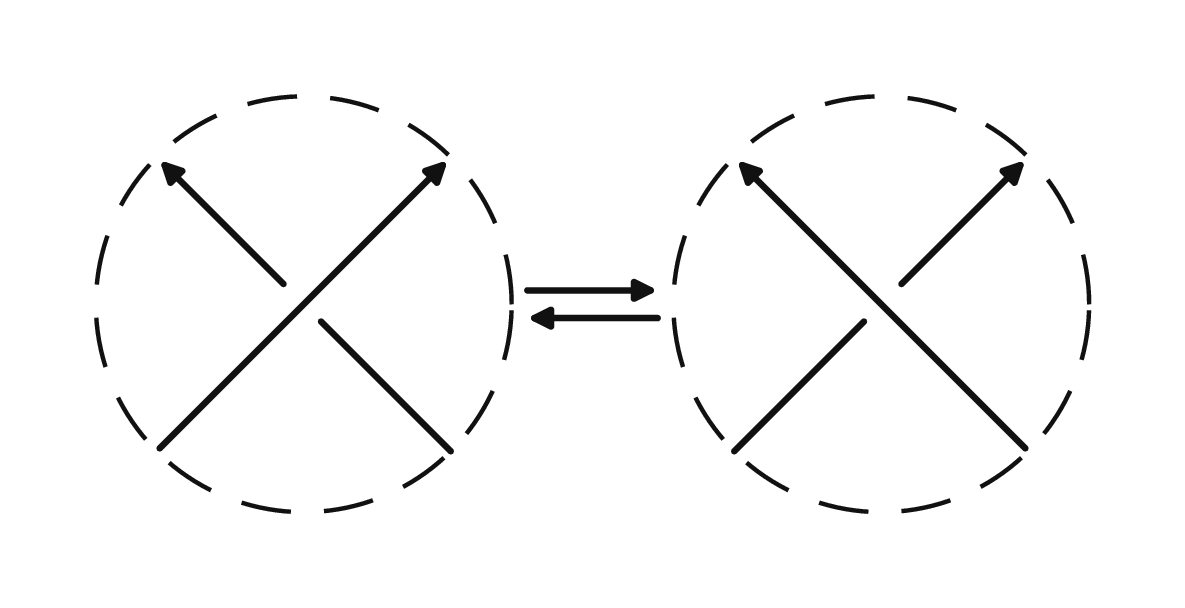}
\caption{Ориентированное изменение перекрестка}
\label{f:cross-ori}
\end{figure}

\begin{assertion}\label{con-lk} Для любых ориентированных двухкомпонентных плоских диаграмм $D_+$ и $D_-$ (соответственно слева и справа на рисунке~\ref{f:cross-ori}), 
отличающихся ориентированным изменением перекрестка, образованного пересечением разных компонент, имеем \linebreak $\lk D_+-\lk D_- = 1$. 
\end{assertion}


\begin{assertion}\label{p:lk-pr} 
(a) Изменение ориентации любой из компонент 
меняет знак коэффициента зацепления.

(b) Перестановка компонент сохраняет коэффициент зацепления.
\end{assertion}

Для п. (b) полезна следующая лемма. 

\begin{lemma}\label{l:half} Коэффициент зацепления равен $(n_{12}-n_{21})/2$, где
 
$\bullet$ $n_{12}$ — сумма знаков в тех перекрестках, в которых первая компонента проходит над второй, и

$\bullet$ $n_{21}$ — сумма знаков в тех перекрестках, в которых вторая компонента проходит над первой.
\end{lemma}

В доказательстве (и далее) вы можете  
использовать без доказательства следующую {\it лемму о тривиальности}: для любых двух замкнутых ориентированных ломаных общего положения на плоскости  
сумма знаков их точек пересечения равна нулю.
Обсуждение и доказательство см. в \S1.3 <<Число пересечений для ломаных на плоскости>> \cite{Sk18}, \cite{Sk}.

\begin{lemma}\label{l:add} Пусть $l_1,l_2,l_3$ — ломаные, соединяющие точку $A$ с точкой $B$, и не пересекающие замкнутую ломаную $p$.
Тогда
$$\lk(l_1\overline{l_2}, p) + \lk(l_2\overline{l_3}, p) = \lk(l_1\overline{l_3}, p).$$
\end{lemma}

\section{Основные результаты о наборах коэффициентов зацепления}\label{s:main}

\begin{figure}[H]
\centering
\includegraphics[scale=1.5]{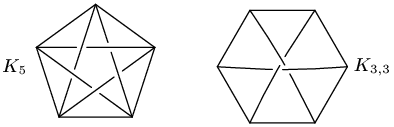}\qquad
\includegraphics[scale=0.8]{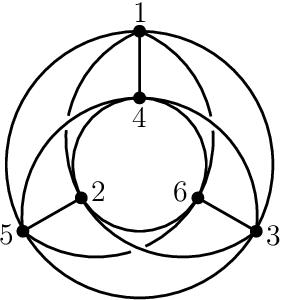}
\caption{Вложения графов 
в $\R^3$}
\label{f:ex-graph-emb}
\end{figure}

\emph{Вложение} $K\to\R^3$ графа $K$ в трехмерное пространство $\R^3$ --- это 
изображение без самопересечений графа $K$ (рисунок~\ref{f:ex-graph-emb}).
Мы рассматриваем такие вложения графов,
для которых ребра изображаются ломаными.
Простое строгое определение таких вложений дано в \cite[Remark 8.4]{AKS}.
Но вам не нужно знать его для решения задач, поскольку все результаты могут быть переформулированы и доказаны в терминах 
\emph{плоских диаграмм}, см. замечание~\ref{r:emvsdi}. 


Обозначим через

$\bullet$ $[n]$ множество $\{1, 2, \ldots, n\}$;

$\bullet$ $K_n$ полный граф с множеством вершин $[n]$;   

$\bullet$ $K_{m,n}$ полный двудольный граф с долями $[m]$ и $[n]'$ (мы обозначаем через $A'$ копию $A$).  

\begin{theorem}[Конвей--Гордон--Закс]\label{t:cgs}
Для любого вложения $K_6\to\R^3$ сумма коэффициентов зацепления по всем неупорядоченным парам непересекающихся циклов длины 3 в $K_6$ нечетна.
\end{theorem}

Этот результат не является частью данного проекта, но лишь служит мотивацией для теорем~\ref{t:div3}, \ref{t:odd}, 
утверждения~\ref{t:div3m} и гипотезы~\ref{t:div3mp}. 
 
\emph{Из этого и следующих разделов для решения в рамках проекта предназначены только теоремы~\ref{t:div3},~\ref{t:odd},~\ref{p:k5nontr}, 
 утверждение~\ref{t:div3m} и задача~\ref{a:crossing}.
Некоторые из них могут быть не очень простыми без подсказок, данных в следующих разделах.}

\begin{remark}[переформулировка в терминах плоских диаграмм]\label{r:emvsdi}
\emph{Плоская диаграмма вложения графа в $\R^3$} определяется   аналогично плоской диаграмме зацепления (используя понятие \emph{изображения графа на плоскости общего положения}, определенного в \cite[\S8]{AKS}).
\emph{Плоская диаграмма графа (абстрактная)} определена в \cite[\S8]{AKS} аналогично (абстрактной) плоской диаграмме (определенной в \S\ref{s:lkpd}). 
Теорема~\ref{t:cgs} переформулируется так:  
\emph{для любой плоской диаграммы графа $K_6$ сумма коэффициентов зацепления по всем неупорядоченным парам непересекающихся треугольников в $K_6$ нечетна.} 
Аналогично переформулируются все те результаты этого текста, которые сформулированы в терминах вложений графов в $\R^3$.   
\end{remark}
 
\begin{theorem}\label{t:div3}  
Для любого графа $K$ существует вложение $K\to\R^3$ такое, что коэффициент зацепления любых двух непересекающихся циклов длины 3 делится на 3. 
\end{theorem} 

Это было доказано Р. Никкуни в 2026 году. 
Для графов с менее чем 6 вершинами теорема~\ref{t:div3}, утверждение~\ref{t:div3m} и гипотеза~\ref{t:div3mp} не являются содержательными. 
В теореме~\ref{t:div3} (и, предположительно, в гипотезе~\ref{c:borr}) mod 3 можно заменить на mod $m$ для любого нечетного целого $m$. 
Но не на mod 2 (ввиду теоремы~\ref{t:cgs}).

\begin{assertion}\label{t:div3m}  
То же, что и в теореме~\ref{t:div3}, с заменой <<делится на 3>> на <<не равно нулю>>. 
\end{assertion}

\begin{conjecture}\label{t:div3mp}  
Для любого графа $K$ существует вложение $K\to\R^3$ такое, что для любых трех попарно непересекающихся циклов длины 3 
наибольший общий делитель их попарных коэффициентов зацепления равен 1.
\end{conjecture}

Аналог без  `длины 3' неверен ввиду теоремы Негами \cite{Ne91, PS05} 
(и ввиду \cite[Theorem 3.1]{FD05}, дающей конкретную константу для соответствующего частного случая теоремы Негами).  

\begin{problem}\label{p:k6} Для любого целого $k$ существует вложение $K_6\to\R^3$ такое, что коэффициент зацепления

$\bullet$ одной (неупорядоченной) пары непересекающихся циклов длины $3$
равен $2k+1$, а

$\bullet$ любой другой пары непересекающихся циклов длины $3$
равен нулю.
\end{problem}

Этот элементарный результат появился в недавних исследованиях \cite[Proposition 1.2]{KS20}
 и \cite[конец \S2.2]{AT21}, хотя мог быть известен и ранее.

\begin{theorem}\label{t:odd}
Возьмем любые $10$ целых чисел $n_{123,456},n_{124,356},\ldots$, соответствующих $10$ неупорядоченным разбиениям множества $[6]$ на два $3$-элементных подмножества.
Если сумма этих чисел нечетна, то существует вложение $K_6\to\R^3$ такое, что
коэффициент зацепления каждой пары $\{ijk,pqr\}$ непересекающихся циклов длины $3$ равен $n_{ijk,pqr}$.
Здесь мы считаем, что $i<j<k$ и $p<q<r$, и берем ориентации $(ijk)$ и $(pqr)$ на циклах. 
\end{theorem}

Это доказано в \cite{Mi}; это следует из лемм~\ref{l:move} и \ref{l:expr} ниже, ср. с теоремой~\ref{t:odd1}. 

\begin{problem}[нерешенная]\label{op:set_lk} Для каждого $n$ найдите все упорядоченные наборы целых чисел, реализуемые в качестве наборов коэффициентов зацепления вложения $K_n\to\R^3$.
\end{problem}


\emph{Общие подсказки.} Описывайте сложные конструкции с помощью простых модификаций простых конструкций! 
Большинство представленных геометрических задач сводятся к простым алгебраическим задачам.

\begin{figure}[H]\centering
\includegraphics[scale=.1]{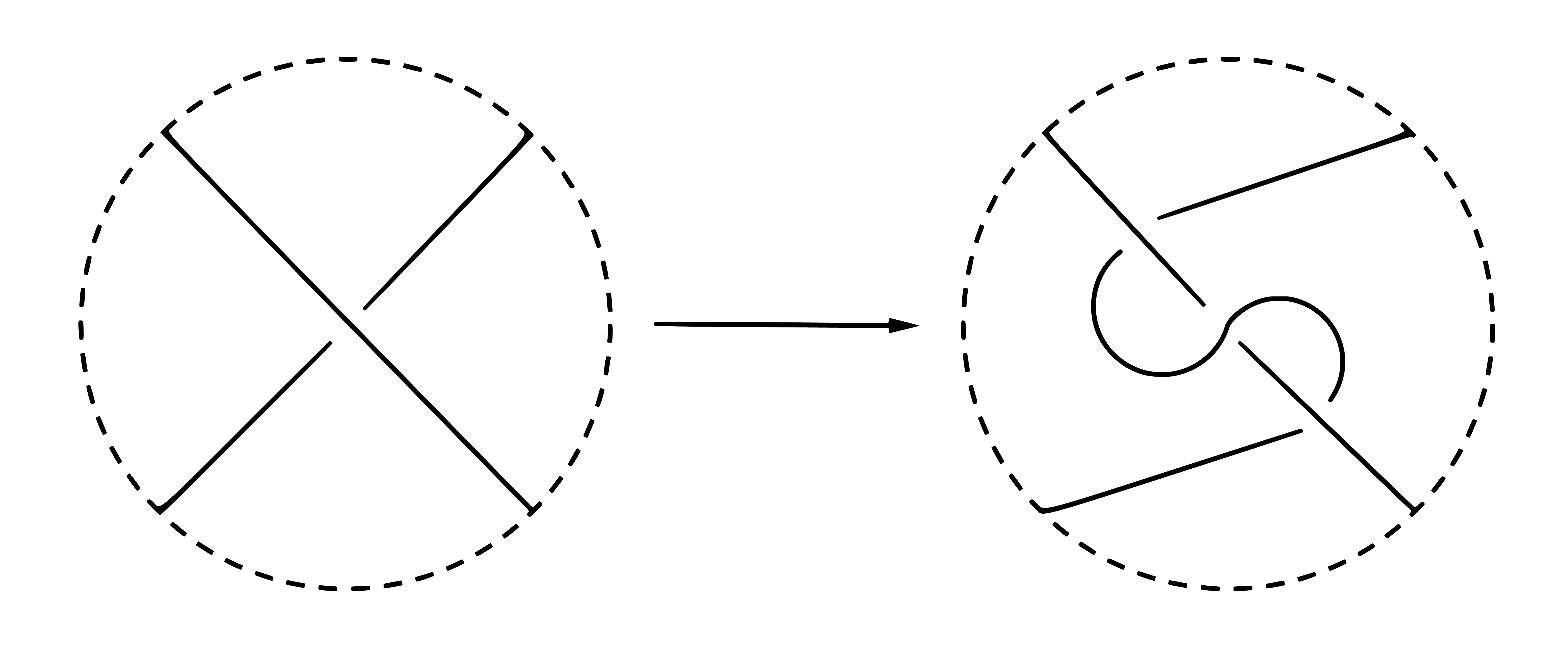}
\caption{Утроение перекрестка} 
 \label{f:triplic}
\end{figure}

\smallskip
\emph{Подсказка к теореме~\ref{t:div3}:} используйте рисунок \ref{f:triplic} и лемму \ref{l:half}.

\textbf{Обобщение теоремы~\ref{t:cgs} Конвея--Гордона--Закса на многократную зацепленность.}

Следующая гипотеза не является частью этого проекта, но служит мотивацией для теоремы~\ref{t:div3}. 
 
\begin{conjecture}\label{c:borr} 
Пусть $K_{11}\to\R^3$ --- вложение такое, что коэффициент зацепления любых двух непересекающихся циклов длины $3$ делится на $3$.  
Тогда существует <<борромеева по модулю 3>> тройка попарно непересекающихся циклов длины $3$ в $K_{11}$, т.е. тройка, у которой тройной коэффициент зацепления по модулю $3$ не равен нулю. 
\end{conjecture}

\textbf{Замечание.} 
По поводу колец Борромео, см. рисунок \ref{f:borwhitre}.b и, например, \cite{Zi10, Sk24}, \cite[\S4]{Sk} вместе с приведенными там ссылками. 
Целочисленный \emph{тройной коэффициент зацепления (Милнора)} определен для ориентированного трехкомпонентного зацепления, все попарные коэффициенты зацепления которого равны нулю, см., например,
\cite[\S2, Theorem~(3)]{MM01}, 
\cite{Sk24}, \cite[\S4.9]{Sk}. 
Аналогично определяется \emph{тройной коэффициент зацепления} по модулю $m$ для ориентированного трехкомпонентного зацепления, все попарные коэффициенты зацепления которого делятся на $m$.
Таким образом, для содержательности гипотезы~\ref{c:borr} необходима теорема~\ref{t:div3}. 
Согласно утверждению~\ref{t:div3m} и гипотезе~\ref{t:div3mp}, в гипотезе~\ref{c:borr} нельзя опустить предположение о <<попарной тривиальности по модулю $3$>>, даже если в заключении заменить <<по модулю 3>> на <<по модулю $m$ для некоторого $m>1$>>. 

Для достаточно больших $n$ и любого вложения $K_n\to\R^3$ в графе $K_n$ существуют три попарно непересекающихся цикла (не обязательно длины 3), образы которых образуют зацепление, изотопное кольцам Борромео.
Это следует из теоремы Негами \cite{Ne91, PS05}.
Однако аналог гипотезы~\ref{c:borr} с заменой $K_{11}$ на $K_n$ при достаточно больших $n$ нетривиален.

По поводу других, но близких проблем о неотъемлемой $n$-зацепленности см., например, обзор \cite{Na20} и приведенные там ссылки.

Предполагаемое доказательство и контрпример для $K_{10}$ вместо $K_{11}$ приведены в \cite{AS}.

Подробнее см. \cite[end of \S4]{AKS}.

\section{Основные результаты об изотопности и гомологичности}\label{s:maini}

\begin{figure}[H]\centering
\includegraphics[scale=0.8]{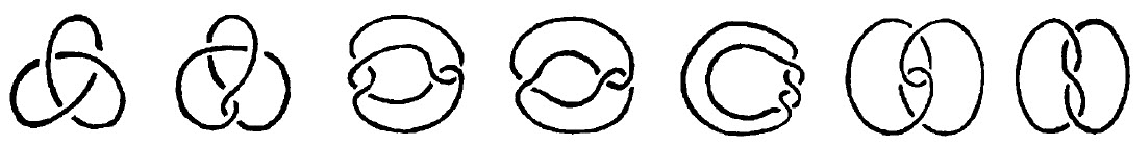}
\caption{Узлы, изотопные трилистнику} 
\label{f:trefoil}
\end{figure}

Под {\it изотопией} зацепления мы понимаем его непрерывную деформацию в трехмерном пространстве как набора тонких упругих нитей; на протяжении всей деформации не допускаются самопересечения (ни одной нити с собой, ни разных нитей).
См. рисунок~\ref{f:trefoil} и сведение к плоским диаграммам в \S\ref{s:isot}. 

\begin{pr}\label{a:crossing} Пожалуйста, докажите эти утверждения на интуитивном уровне (поскольку у вас нет строгого определения изотопии). 

(a) Все узлы (представленные своими плоскими диаграммами) на
рисунке~\ref{f:trefoil} изотопны друг другу.

(b) Зацепление Уайтхеда (рисунок~\ref{f:borwhitre}.w) изотопно зацеплению с рисунка~\ref{f:hopf}.d.

(c) Две плоские диаграммы зацепления получены проекциями на разные плоскости.
Тогда эти диаграммы получены проекциями на одну плоскость из двух (различных, но) изотопных зацеплений. 
\end{pr}

\emph{Подсказка к (a).} <<Вероятно, лучший способ решить эту задачу --- сделать модель узла трилистника из веревки или шнурка, 
а затем перемещать ее из одного положения в другое>> \cite[\S2]{Pr95}.

\begin{theorem}\label{p:k5nontr}  (a) Существуют неизотопные вложения $K_5\to\R^3$, сужения которых на любой простой цикл в $K_5$ изотопны. 

(b) (киральность \cite{Si86}) Никакое вложение $K_5\to\R^3$ не изотопно своему зеркальному образу. 

(a',b') То же для $K_{3,3}$ вместо $K_5$. 
\end{theorem}

Вложения из пп. (a,a') различаются целыми \emph{числом Ву--Саймона} и \emph{двудольным числом Ву--Саймона}, определенными в \cite[\S9]{AKS}. 
При помощи этих же инвариантов доказываются пп. (b,b').  

\begin{remark}[вложения с точностью до узлов]\label{r:emk}
\emph{Изотопическая} классификация узлов настолько сложна, что вряд ли может быть получена.
Если вложения $K\to\R^3$ изотопны, то их сужения на любой простой цикл в $K$ изотопны.
Таким образом, вряд ли удастся получить изотопическую классификацию вложений графов в $\R^3$.
Поэтому интересно определить и классифицировать <<вложения с точностью до узлов>>.
Существует несколько формализаций этой идеи: \emph{гомологичность} (определена в \cite[\S10]{AKS}), \emph{сингулярная гомотопия},
\emph{почти изотопия}, \emph{необъемлемая изотопия}; см. обзор \cite{Sk'}.
\end{remark}


\begin{theorem}\label{t:tani}
Два вложения графа (не обязательно связного) в пространство гомологичны тогда и только тогда, когда гомологичны их сужения на любой подграф, гомеоморфный $K_5,K_{3,3}$ или $K_3\sqcup K_3$.
Или, эквивалентно, когда все коэффициенты зацепления, числа Ву--Саймона и двудольные числа Ву--Саймона (сужений отображений $f,g$) равны.
\end{theorem}

Эта теорема не является частью данного проекта. 
Она доказана в \cite{ST03} с использованием гомологической классификации в терминах \emph{инварианта Ву--Саймона} \cite{Ta95} и
 вычисления некоторой \emph{группы гомологий} $H_2^{\text{skew-sym}}(K\times K-\diag;\Z)$
 \cite{Ni00}; см. подробнее в \cite[Remark~9.9]{AKS}.

\begin{problem}[нерешенная; ср. с задачей \ref{op:set_lk}]\label{op:set_lkwu}* 
(a) Для каждого $n$ найдите наборы целых чисел, реализуемые как наборы коэффициентов зацепления и чисел Ву--Саймона вложений $K_n\to\R^3$. (начните с $n=6$). 

(b) Аналогичный вопрос с добавлением двудольных чисел Ву--Саймона. 
\end{problem}


\section{Реализация коэффициентов зацепления: к доказательству теоремы~\ref{t:odd}}\label{s:real}

Переформулируем теорему~\ref{t:odd}. 
Обозначим через $X$ множество неупорядоченных разбиений
множества $[6]$ 
на два 3-элементных подмножества. 
{\bf Функция зацепленности} $l_f:X\to\Z$ вложения $f:K_6\to\R^3$ сопоставляет каждому разбиению $x=\{\{i,j,k\},\{p,q,r\}\}\in X$, где $i<j<k$ и $p<q<r$, коэффициент зацепления $f(ijk)$ и $f(pqr)$. 
(Множество всех функций $X \to \Z$ можно рассматривать как множество всех векторов с целыми координатами, индексированных разбиениями из $X$.)
 
\begin{theorem}[ср. с теоремой~\ref{t:odd}]\label{t:odd1} 
Для любой функции $\varphi: X \to \Z$ с нечетной суммой значений существует вложение $f$ с $l_f=\varphi$. 
\end{theorem}

\begin{figure}[H]\centering
\includegraphics[scale=0.6]{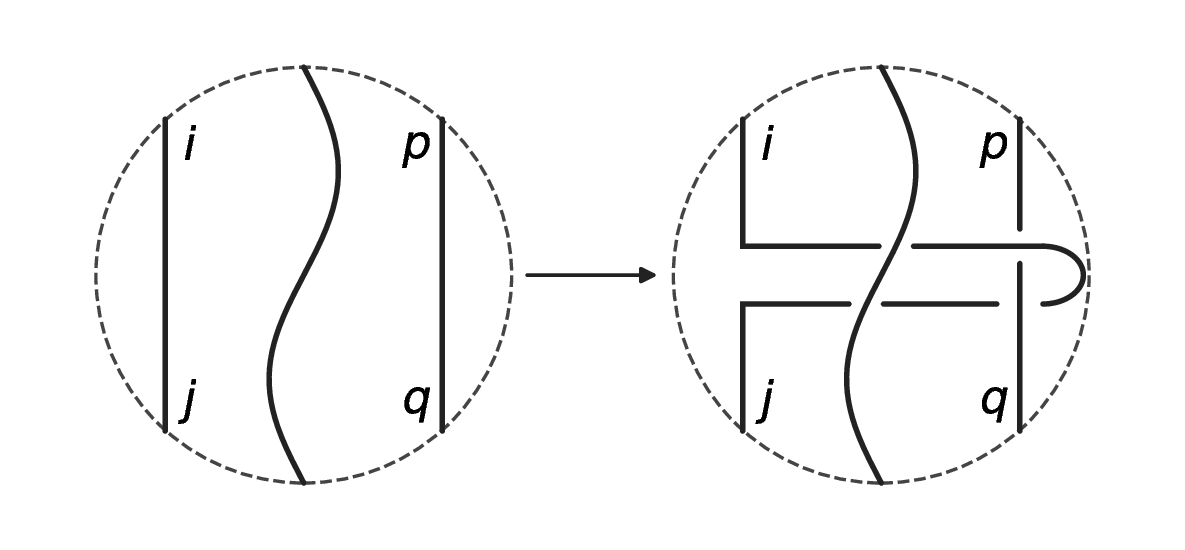}\qquad
\caption{Пальцевое движение одного ребра через другое}
\label{f:finger_gra3spa}
\end{figure}

\begin{pr}\label{p:move} (a) Как меняется $l_f$ при преобразовании вложения $f:K_6\to\R^3$, показанном на рисунке \ref{f:finger_gra3spa}? 

(b) (загадка) Тот же вопрос для вложения $f:K_n\to\R^3$. 
\end{pr}

Решив задачу \ref{p:move}.b, вы можете доказать утверждение~\ref{t:div3m}.

\emph{Интересной} парой назовем упорядоченную пару $(ij,pq)$ непересекающихся ориентированных ребер графа $K_6$ таких, что $i<j$ и $p<q$. 
Возьмем интересную пару $(ij,pq)$. 
Возьмем $a,b\in [6]$ такие, что $\{i,j,p,q,a,b\}=[6]$. 
Определим функцию $v_{ij,pq} : X\to\Z$ 
по формуле\footnote{Хотя формулы используют выбор $a$ и $b$ из $\{a,b\}$, они не зависят от этого выбора.}
$$v_{ij,pq}(\{i,j,a\},\{p,q,b\}) := \sgn(a-i)\sgn(a-j)\sgn(b-p)\sgn(b-q) \in \{ +1, -1 \},$$
$$v_{ij,pq}(\{i,j,b\},\{p,q,a\}) := \sgn(b-i)\sgn(b-j)\sgn(a-p)\sgn(a-q) \in \{ +1, -1 \},$$
и $v_{ij,pq}(x) := 0$ для всех остальных $x\in X$.

Например, $v_{15,23}(\{2,3,5\},\{1,4,6\})=0$, а 
$$v_{15,23}(\{2,3,6\},\{1,4,5\}) = v_{15,23}(\{1,5,4\},\{2,3,6\}) = 1\cdot(-1)\cdot1\cdot1 = -1.$$

\begin{pr}\label{p:intepair}  Как зависит функция $v_{ij,pq}$ от перестановки ребер $ij,pq$?
\end{pr}
 

\begin{lemma}\label{l:move} 
Для любого вложения $f:K_6\to\R^3$ и любой интересной пары $(ij,pq)$ существуют вложения $f_+,f_-:K_6\to\R^3$ такие, что $l_{f_\pm}-l_{f}=\pm v_{ij,pq}$.
\end{lemma}

Функция $v: X \to \Z$ называется \emph{реализуемой}, если $v$ равна некоторой целочисленной линейной комбинации функций $v_{ij,pq}$ по всем интересным парам $(ij,pq)$. 

\begin{lemma}[\cite{Mi}]\label{l:expr} Любая функция $X \to \Z$ с четной суммой значений реализуема.
\end{lemma}

Лемма следует из утверждений~\ref{a:coord}.ab. 

Для любого $x\in X$ определим \emph{индикаторную функцию}
$$\chi_x:X\to\Z\quad\text{следующим образом}\quad \chi_x(y)=\begin{cases}1 &\text{если }x=y \\ 0&\text{в противном случае}\end{cases}.$$

\begin{assertion}\label{a:coord}
(a) Для любых двух разбиений из $X$ либо сумма, либо разность их индикаторных функций реализуема.

(b) Удвоенная индикаторная функция любого разбиения из $X$ реализуема.
\end{assertion}


\textbf{Замечание.} Вот другое, более концептуальное определение (не используемое в этом проекте). 
\emph{Интересной} парой назовем упорядоченную пару $(ij,pq)$ непересекающихся ориентированных ребер графа $K_6$. 
Возьмем интересную пару $(ij,pq)$. 
Возьмем $a,b\in [6]$ такие, что $\{i,j,p,q,a,b\}=[6]$. 
Определим функцию $v_{ij,pq} : X\to\Z$ по формуле
$$v_{ij,pq}(\{i,j,a\},\{p,q,b\}) := \sgn(i,j,a)\sgn(p,q,b),\quad v_{ij,pq}(\{i,j,b\},\{p,q,a\}) := \sgn(i,j,b)\sgn(p,q,a),$$
и $v_{ij,pq}(x) := 0$ для всех остальных $x\in X$.
Здесь {\it знак} вектора из 3 целых чисел равен 
$$\sgn(k,l,m) := \sgn(k-l)\sgn(l-m)\sgn(m-k)\in\{+1,-1\}.$$ 
Это знак перестановки, упорядочивающей $k,l,m$ по возрастанию. 

Например, $v_{15,23}(\{2,3,5\},\{1,4,6\})=0$, а 
$$v_{15,23}(\{2,3,6\},\{1,4,5\}) = v_{15,23}(\{1,4,5\},\{2,3,6\}) = \sgn(1,5,4)\sgn(2,3,6) = (-1)\cdot1=-1.$$

\section{Некоторые ответы, подсказки и решения}\label{s:sol}

\begin{figure}[h]\centering
\includegraphics[scale=0.5]{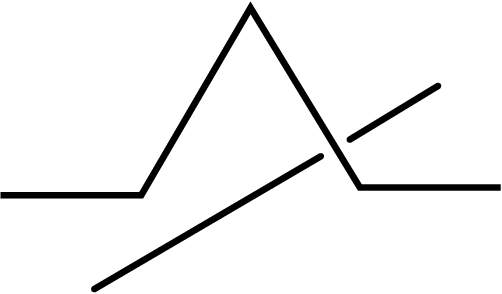}
\caption{Мостик над перекрестком}
\label{f:bridge}
\end{figure}

{\bf \ref{p:diag}.} См. рисунок \ref{f:bridge}.
На верхнем ребре каждого перекрестка плоской диаграммы выберем две близкие точки по разные стороны от перекрестка.
Заменим отрезок ломаной между двумя выбранными точками на <<мостик>>, поднимающийся над плоскостью диаграммы и соединяющий эти две точки. 
Заменим соответствующими мостиками отрезки ломаных во всех перекрестках.
Получим узел, проекция которого является объединением отрезков данной плоской диаграммы.

\smallskip
{\bf \ref{a:lk2}.} \textit{Ответы:} 
\qquad (a) 1 для зацепления Хопфа и 0 для остальных зацеплений; \qquad (b) 0;  

(c) 1 для пары циклов $(123)$ и $(456)$; 0 для остальных пар.
 
\smallskip
{\bf \ref{a:lk2-pr}.} Возьмем плоскую диаграмму зацепления. 
По лемме о четности (сформулированной после утверждения~\ref{a:lk2-pr})
число перекрестков, в которых первая компонента проходит над второй, имеет ту же
четность, что и число перекрестков, где вторая компонента проходит над первой.
Это завершает доказательство.

\smallskip
{\bf \ref{e:lk}.} \textit{Ответы:}
 \qquad (a) 1; 0; 2; 0 
(знаки могут заменяться на противоположные в зависимости от выбранной ориентации);

  \qquad(b) 0;  \qquad(c) +1 для пары циклов $(123)$ и $(456)$; 0 для остальных пар.

\smallskip
{\bf \ref{p:lk-pr}.} (a) Изменение ориентации любой из компонент меняет знак в каждом перекрестке на противоположный.

(b) Перестановка компонент меняет знак в каждом перекрестке на противоположный. 
Таким образом, п. (b) следует из леммы~\ref{l:half}.

\smallskip
{\bf \ref{l:half}.} Это следует из того, что по лемме о тривиальности (сформулированной после леммы) $n_{12}+n_{21}=0$.  
 
\smallskip
{\bf \ref{t:div3}.} Следуем тексту Р. Никкуни, с его любезного разрешения. 

Заменим каждый перекресток тремя перекрестками как показано на рисунке \ref{f:triplic}.
Возьмем любые два непересекающихся ориентированных цикла длины $3$ в графе. Тогда определен знак каждого перекрестка этих циклов, и в исходной, и в новой диаграмме.
Для каждого из перекрестков исходной диаграммы ребро, проходившее \emph{сверху}, после замены проходит \emph{сверху} в двух из новых перекрестков и \emph{снизу} в оставшемся новом перекрестке.
Знак каждого из этих двух новых перекрестков совпадает со знаком исходного перекрестка. Знак оставшегося нового перекрестка противоположен знаку исходного перекрестка.
Используем для исходной диаграммы обозначения $n_{12}$ и $n_{21}$ из леммы~\ref{l:half}, а для новой диаграммы аналогично определим $n'_{12}$ и $n'_{21}$.
Для каждого исходного перекрестка, вносившего вклад $\pm1$ в сумму $n_{12}-n_{21}$, соответствующие три новых перекрестка вносят вклад $\pm(1-(-1)+1)=\pm3$ в сумму $n'_{12}-n'_{21}$.
Тогда для любых двух непересекающихся циклов сумма $n'_{12}-n'_{21}$ делится на $3$. Значит, по лемме~\ref{l:half} полученная плоская диаграмма является требуемой.


\smallskip
{\bf \ref{t:div3m}.} \emph{Указание.} Выберем пару циклов с нулевым коэффициентом зацепления. Модифицируем плоскую диаграмму так, чтобы коэффициент зацепления в этой паре циклов стал ненулевым, а новых пар с нулевыми коэффициентами зацепления не появилось. Для этого многократно применим лемму~\ref{l:move}. При этом коэффициенты зацепления нужной и некоторых других пар циклов изменятся на $\pm N$ для произвольного наперед заданного целого $N$.

\smallskip
{\bf \ref{p:k6}.} Это следует из леммы~\ref{l:move}, утверждения~\ref{a:coord}.a и результата задачи~\ref{e:lk}.c.

\begin{figure}[h]\centering
\includegraphics[scale=0.25]{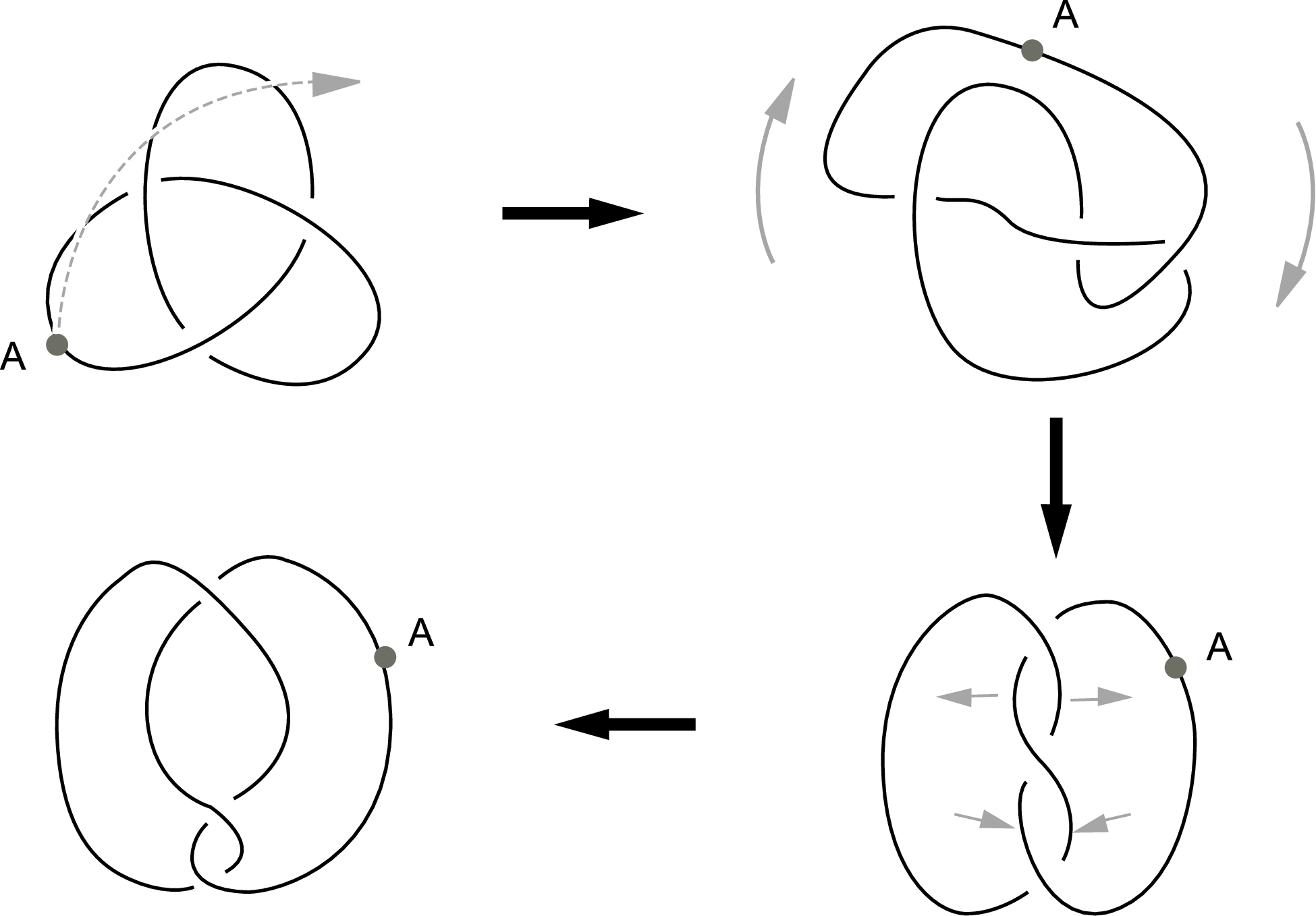}
\caption{Изотопия трилистника}
\label{f:itref}
\end{figure}

\smallskip
{\bf \ref{a:crossing}.} (a) Выполнить некоторые преобразования трилистника поможет рисунок \ref{f:itref} (выполненный Д. Кроо). 

\smallskip
{\bf \ref{p:move}.} (a) См. лемму~\ref{l:move}. 

\smallskip
{\bf \ref{p:intepair}.} Не изменяется.

\smallskip
{\bf \ref{l:move}.} Приведем упрощенное изложение доказательства из \cite{KS20}.

Приведем решение для вложения $f_-$. 
Неформально, мы получаем $f_-$, оборачивая ребро $f(ij)$ на $-1$ оборот вокруг ребра $f(pq)$, как на рисунке \ref{f:finger_gra3spa}.
Приведем строгую конструкцию.

На ребре $ij$ выберем близкие точки $I$ и $J$, между которыми нет перекрестков, так чтобы ориентация отрезка $IJ$ совпадала с ориентацией ребра $ij$.
На ребре $pq$ аналогично выберем точки $P$ и $Q$, так, чтобы ориентация отрезка $PQ$ совпадала с ориентацией ребра $pq$.
При помощи пальцевого движения заменим отрезок $IJ$ исходной диаграммы (т.е. диаграммы вложения $f$) на несамопересекающуюся ориентированную ломаную $l_1$ общего положения, так, чтобы она проходила через точки $P$ и $Q$, как показано на рисунке \ref{f:finger_gra3spa}. В частности, в перекрестке $P$ ломаная $l_1$ проходит сверху, а в перекрестке $Q$ и во всех остальных перекрестках --- снизу.

Рассмотрим в исходной диаграмме любые два непересекающиеся цикла $x$ и $y$ длины $3$.
Обозначим через $x_-$ и $y_-$ соответствующие им циклы в полученной диаграмме (т.е. диаграмме вложения $f_-$).
Обозначим через $\overline x$ разбиение множества вершин графа $K_6$ на два трехэлементных подмножества, в котором одно из подмножеств состоит из вершин цикла $x$.
Везде далее после введения ориентации на циклах $x$ и $y$ будем брать согласованные с ними ориентации на циклах $x_-$ и $y_-$.

Если ребро $ij$ не входит ни в один из циклов $x$ и $y$, то при любой ориентации на этих циклах выполняется $$\lk(x_-,y_-) = \lk(x,y).$$ Отсюда $l_{f_-}(\overline x)-l_f(\overline x)=0.$

Если ребро $ij$ входит в один из циклов $x$ и $y$, то не уменьшая общности оно входит в цикл $x$. Возьмем на цикле $x$ ориентацию такую же, как у ребра $ij$. Будем считать цикл $x$ первым.
Обозначим через $l_2$ ориентированную ломаную, являющуюся дополнением до ориентированного отрезка $IJ$ в цикле $x$.

Тогда если цикл $y$ не содержит ребро $pq$, то при любом выборе ориентации на нем, по лемме~\ref{l:add}
выполняется $$\lk(x_-,y_-)=\lk(l_1l_2,y)=\lk(l_1JI,y)+\lk(IJl_2,y)=0+\lk(x,y).$$ Отсюда $l_{f_-}(\overline x)-l_f(\overline x)=0$.

Если же цикл $y$ содержит ребро $pq$, то обозначим вершины циклов $x$ и $y$ так, что $x=ija, y=pqb$. 
Возьмем на цикле $y$ ориентацию такую же, как у ребра $pq$.
Тогда знак перекрестка в точке $P$ равен $-1$.
Следовательно, по лемме~\ref{l:add} выполняется $$\lk(x_-,y_-)=\lk(l_1l_2,y)=\lk(l_1JI,y)+\lk(IJl_2,y)=-1+\lk(x,y).$$
Отсюда $l_{f_-}(\overline x)-l_f(\overline x)=-\sgn(a-i)\sgn(a-j)\sgn(b-p)\sgn(b-q)$.

Из всего вышенаписанного вытекает формула $l_{f_-}-l_f=-v_{ij,pq}$.

Конструкция вложения $f_+$ получается из вышеприведенной конструкции изменением перекрестков в точках $P$ и $Q$ (проход меняется на переход и наоборот).  
 
\smallskip
{\bf \ref{a:coord}.} Следуем \cite{Mi}. 

(a) Для любых различных разбиений $x,y \in X$ найдутся $\alpha \in x,\beta\in y$ такие, что $|\alpha\cap\beta|=2$. 
Тогда существует интересная пара $(ij,pq)$ такая, что $\{i,j\}=\alpha\cap\beta$ и $\{p,q\}=\overline\alpha\cap\overline\beta$. 
Следовательно, $v_{ij,pq}(z)= \begin{cases}\pm1&\text{если }z=x\text{ или }z=y\\0&\text{в противном случае}\end{cases}$.
Значит, $\pm v_{ij,pq}=\chi_x\pm \chi_y$.

(b) Возьмем любое разбиение $x\in X$, и обозначим $s := \{\{1,2,3\},\{4,5,6\}\} \in X$. 
По п. (a), $\chi_x\pm\chi_s$ реализуема, следовательно, $2\chi_x\pm2\chi_s$ реализуема. 
Тогда $2\chi_s = v_{12,45}+v_{16,34}+v_{23,56}$, значит, $2\chi_s$ реализуема. 
Следовательно, $(2\chi_x\pm 2\chi_s)\mp2\chi_s = 2\chi_x$ также реализуема.

(По-иному, для любой перестановки $\sigma$ множества $[6]$ имеем $\pm2\chi_{\sigma(s)} = v_{\sigma(12,45)}+v_{\sigma(16,34)}+v_{\sigma(23,56)}$,
где мы расширяем определение $v_{ij,kl}$ с помощью замечания после утверждения \ref{a:coord}.)  


\section{Изотопность на языке плоских диаграмм}\label{s:isot}

\begin{figure}[H]\centering
\includegraphics[scale=0.45]{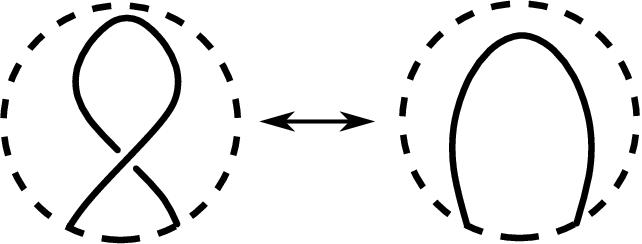}\qquad
\includegraphics[scale=0.45]{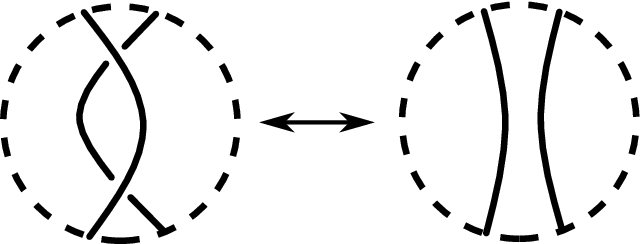}\qquad
\includegraphics[scale=0.45]{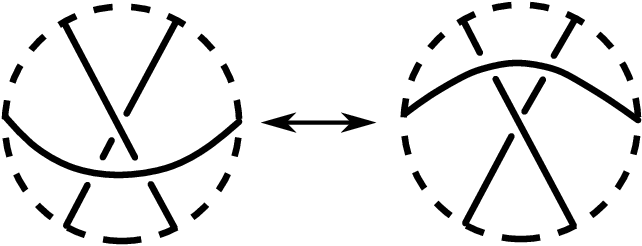}
\caption{Движения Райдемайстера}
\label{f:reid}
\end{figure}

Вместо зацеплений с точностью до изотопии в этом тексте мы изучаем 
плоские диаграммы с точностью до эквивалентности, порожденной {\bf движениями Райдемайстера} (рисунок~\ref{f:reid}) и {\it плоской изотопией}.
Под {\it плоской изотопией} 
плоской диаграммы мы понимаем ее непрерывную деформацию на плоскости без появления самопересечений. 
Это мотивируется следующей теоремой (которая не является частью данного проекта). 

\begin{theorem}[Райдемайстер; {\cite[\S1.7]{PS96}}]\label{l:reid} Два зацепления изотопны тогда и только тогда, когда некоторая плоская диаграмма первого зацепления может быть получена из некоторой плоской диаграммы второго зацепления с помощью движений Райдемайстера и плоской изотопии.
\end{theorem}

\begin{pr}\label{a:crossing-r} Для первой пары узлов на рисунке~\ref{f:trefoil} разложите вашу изотопию между ними на движения Райдемайстера, показанные на рисунке~\ref{f:reid}, и плоскую изотопию.
\end{pr}

\begin{lemma}\label{l:lk2} (a) Коэффициент зацепления по модулю 2

 \qquad (b) Коэффициент зацепления
 
 сохраняется при движениях Райдемайстера.
\end{lemma}

Эта лемма легко доказывается отдельно для каждого движения Райдемайстера.
Для п.~(b) проверьте, что знаки всех перекрестков не меняются.

Подробнее см., например, \cite[\S\S 1-3]{Sk20u}. 

\begin{figure}[H]\centering
\includegraphics[scale=1.]{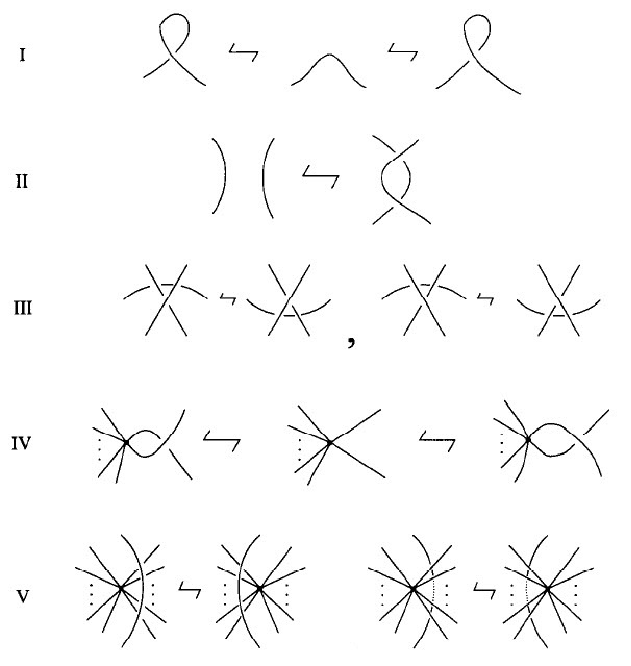}
\caption{Преобразования Райдемайстера для плоских диаграмм графов}
\label{f:reid-1-5}
\end{figure}

Аналогично теореме~\ref{l:reid},
\emph{два вложения графа в $\R^3$ изотопны тогда и только тогда, когда некоторая плоская диаграмма первого вложения (см. замечание~\ref{r:emvsdi}) может быть получена из некоторой плоской диаграммы второго вложения с помощью движений Райдемайстера для плоских диаграмм графов (рисунок~\ref{f:reid-1-5}) и плоской изотопии.} 
Понятие \emph{плоской диаграммы графа} определено в \cite[\S8]{AKS}.

{\it В этом списке звездочками отмечены книги, обзоры и популярные статьи.}